\documentclass[12pt]{report}
\usepackage{verbatim,amsmath,latexsym,amssymb,amsthm}
\usepackage{graphicx}
\usepackage{maple2e}

\theoremstyle{plain}
\newtheorem{thm}{Theorem}[section]

\newtheorem{prop}[thm]{Proposition}
\newtheorem{cor}[thm]{Corollary}

\theoremstyle{definition}
\newtheorem{alg}[thm]{Algorithm}
\newtheorem{df}[thm]{Definition}

\newtheorem{ex}[thm]{Example}

\newtheorem{rem}[thm]{Remark}

\def\zn{z^{(n)}}
\def\zk{z^{(k)}}
\def\zi{z^{(\infty)}}
\def\C{{\mathbb C}}
\def\R{{\mathbb R}}
\def\Z{{\mathbb Z}}
\def\bs{\backslash}
\def\CV{{\cal V}}
\def\CS{{\cal S}}
\def\CD{{\cal D}}

\def\CB{{\cal B}}
\def\CX{{\cal X}}
\def\CC{{\cal C}}
\def\tI{{\tilde I}}

\def\tF{{\tilde F}}

\def\tf{{\tilde f}}
\def\tz{{\tilde z}}
\def\q{{\tilde q}}
\def\x{{\bf x}}
\def\p{{\bf p}}
\def\GLC2{GL(2,\C)}
\def\Maple{{\sc Maple}}
\def\Id{{\rm I}}
\def\inv{{\grave \iota}}
\def\lD{{\bar \CD}}
\def\Eta{{\bar \eta}}
\def\Ji{J^\infty}
\def\as{\alpha}
\newcommand{\matrixb}[4]{\left(\begin{array}{cc} #1 &#2\\ #3& #4 \end{array}\right)}

\newlength{\originalbase}
\newcommand{\spacing}[1]{\setlength{\baselineskip}{#1\originalbase}}

\begin{document}                        


\spacing{1.66}


\begin{titlepage}
\vspace*{1cm}
\begin{centering}
\begin{Large}                         
	Inductive Approach to Cartan's Moving Frame Method with Applications to
Classical Invariant Theory.
\end{Large}                           
\vfill
	A Dissertation\\
	Presented to the Faculty of the Graduate School\\
	of\\
	University of Minnesota\\
	in Candidacy for the Degree of\\
	Doctor of Philosophy\\
\vfill
	by\\
	Irina A. Kogan\\
\vspace{1cm}
	Dissertation Advisor: Professor Peter J. Olver\\
\vspace{1cm}
	June 2000\\
\end{centering}
\end{titlepage}


\begin{titlepage}
    \begin{centering}
    \vspace*{3in}
        \copyright 2000 by Irina A. Kogan.\\
        All rights reserved.\\
    \end{centering}
\end{titlepage}
\pagenumbering{roman}
\begin{center}   
\Large {\bf Acknowledgments}
\end{center}

Let me express my deepest gratitude to my thesis advisor, Peter Olver, who guided my research with encouragement and support. 

I am also thankful to Mark Fels  for many interesting
discussions on Cartan's equivalence method  which has led  me to better
understanding of the subject; to Eric Schost and Gregoire Lecerf from  
Ecole Polytechnique in  France for their 
valuable suggestion on Gr\"obner basis computations;
and to  my undergraduate advisor, Nikolai I. Osetinski,  who first
introduced me to the subject of invariant theory.

I would like to thank the wonderful friends whom I met  in  graduate school. Their 
presence, as well as the friendly  and stimulating environment 
of the School of Mathematics, 
made my graduate studies  such a pleasant experience.

I cannot overvalue the love and support I have received from my parents, Svetlana and Aleksander Kogan; in particular, I would like to thank them for visiting me here in Minneapolis. I am thankful to my sister, Julia, and her family, Maksim, Ilya and Lev for being 
so close, in spite of the huge geographical distance between us. I am fortunate to have a wonderful 
son, Yasha, my constant source of energy and positive emotions. 

I devote this thesis to my grandfather, Ilya Guhn, with love and admiration. 

\newpage
\spacing{1.3}
\vskip-2cm

\begin{centering}
\end{centering}
\begin{center}   
\Large ABSTRACT
\end{center}

\begin{small}
This thesis is devoted to algorithmic aspects of the implementation of 
 Cartan's moving frame method to the problem of the equivalence of 
submanifolds under a Lie group action. We adopt a general definition of a moving frame 
as an equivariant map from the space of submanifolds to the group itself and 
introduce two algorithms, which simplify the construction of such maps. 
The first algorithm  is applicable when the group  factors as a product of 
two subgroups  $G=BA$,  allowing us to use 
 moving frames and differential invariants  for the groups $A$ and $B$ in order 
to construct
a moving frame and differential invariants for $G$. This approach not only simplifies the
computations, but also produces the relations among the
invariants of $G$ and its subgroups. We use the groups of the projective, the affine and the Euclidean transformations on the plane to illustrate the algorithm. 
We also introduce a recursive algorithm,  allowing, provided  the group action satisfies  certain conditions,
 to construct differential invariants order by order, at each step normalizing more and more of the group parameters, at the end obtaining a moving frame for the entire group. The development of this algorithm has been motivated by the  applications 
of the  moving frame method to the problems of the equivalence and symmetry of 
polynomials under linear changes of variables.  In  the complex or real case these 
problems  can be reduced and, in theory,
 completely solved  as the problem of the equivalence of 
submanifolds.  Its  solution  however involves algorithms based on the 
  Gr\"obner basis computations, which  due to their complexity, are not always feasible. 
Nevertheless, some interesting new  results were obtained, such as a classification  of 
ternary cubics and their groups 
of symmetries, and the  necessary and sufficient conditions for a homogeneous polynomial in three variables to be equivalent to $x^n+y^n+z^n.$   
\end{small}
\vspace*{30pt}


\spacing{1.66}
\tableofcontents
%
%
%
 

\makeatletter
 
\def\diagram{\m@th\leftwidth=\z@ \rightwidth=\z@ \topheight=\z@
\botheight=\z@ \setbox\@picbox\hbox\bgroup}
 
\def\enddiagram{\egroup\wd\@picbox\rightwidth\unitlength
\ht\@picbox\topheight\unitlength \dp\@picbox\botheight\unitlength
\hskip\leftwidth\unitlength\box\@picbox}
 
\def\bfig{\begin{diagram}}
\def\efig{\end{diagram}}
\newcount\wideness \newcount\leftwidth \newcount\rightwidth
\newcount\highness \newcount\topheight \newcount\botheight
 
\def\ratchet#1#2{\ifnum#1<#2 \global #1=#2 \fi}
 
\def\putbox(#1,#2)#3{%
\horsize{\wideness}{#3} \divide\wideness by 2
{\advance\wideness by #1 \ratchet{\rightwidth}{\wideness}}
{\advance\wideness by -#1 \ratchet{\leftwidth}{\wideness}}
\vertsize{\highness}{#3} \divide\highness by 2
{\advance\highness by #2 \ratchet{\topheight}{\highness}}
{\advance\highness by -#2 \ratchet{\botheight}{\highness}}
\put(#1,#2){\makebox(0,0){$#3$}}}
 
\def\putlbox(#1,#2)#3{%
\horsize{\wideness}{#3}
{\advance\wideness by #1 \ratchet{\rightwidth}{\wideness}}
{\ratchet{\leftwidth}{-#1}}
\vertsize{\highness}{#3} \divide\highness by 2
{\advance\highness by #2 \ratchet{\topheight}{\highness}}
{\advance\highness by -#2 \ratchet{\botheight}{\highness}}
\put(#1,#2){\makebox(0,0)[l]{$#3$}}}
 
\def\putrbox(#1,#2)#3{%
\horsize{\wideness}{#3}
{\ratchet{\rightwidth}{#1}}
{\advance\wideness by -#1 \ratchet{\leftwidth}{\wideness}}
\vertsize{\highness}{#3} \divide\highness by 2
{\advance\highness by #2 \ratchet{\topheight}{\highness}}
{\advance\highness by -#2 \ratchet{\botheight}{\highness}}
\put(#1,#2){\makebox(0,0)[r]{$#3$}}}

\def\adjust[#1]{} 
 
\newcount \coefa
\newcount \coefb
\newcount \coefc
\newcount\tempcounta
\newcount\tempcountb
\newcount\tempcountc
\newcount\tempcountd
\newcount\xext
\newcount\yext
\newcount\xoff
\newcount\yoff
\newcount\gap%
\newcount\arrowtypea
\newcount\arrowtypeb
\newcount\arrowtypec
\newcount\arrowtyped
\newcount\arrowtypee
\newcount\height
\newcount\width
\newcount\xpos
\newcount\ypos
\newcount\run
\newcount\rise
\newcount\arrowlength
\newcount\halflength
\newcount\arrowtype
\newdimen\tempdimen
\newdimen\xlen
\newdimen\ylen
\newsavebox{\tempboxa}%
\newsavebox{\tempboxb}%
\newsavebox{\tempboxc}%
 
\newdimen\w@dth
 
\def\setw@dth#1#2{\setbox\z@\hbox{\m@th$#1$}\w@dth=\wd\z@
\setbox\@ne\hbox{\m@th$#2$}\ifnum\w@dth<\wd\@ne \w@dth=\wd\@ne \fi
\advance\w@dth by 1.2em}
 
 
\def\t@^#1_#2{\allowbreak\def\n@one{#1}\def\n@two{#2}\mathrel
{\setw@dth{#1}{#2}
\mathop{\hbox to \w@dth{\rightarrowfill}}\limits
\ifx\n@one\empty\else ^{\box\z@}\fi
\ifx\n@two\empty\else _{\box\@ne}\fi}}
\def\t@@^#1{\@ifnextchar_{\t@^{#1}}{\t@^{#1}_{}}}
\def\to{\@ifnextchar^{\t@@}{\t@@^{}}}
 
\def\t@left^#1_#2{\def\n@one{#1}\def\n@two{#2}\mathrel{\setw@dth{#1}{#2}
\mathop{\hbox to \w@dth{\leftarrowfill}}\limits
\ifx\n@one\empty\else ^{\box\z@}\fi
\ifx\n@two\empty\else _{\box\@ne}\fi}}
\def\t@@left^#1{\@ifnextchar_{\t@left^{#1}}{\t@left^{#1}_{}}}
\def\toleft{\@ifnextchar^{\t@@left}{\t@@left^{}}}
 
\def\two@^#1_#2{\allowbreak
\def\n@one{#1}\def\n@two{#2}\mathrel{\setw@dth{#1}{#2}
\mathop{\vcenter{\lineskip\z@\baselineskip\z@
                 \hbox to \w@dth{\rightarrowfill}%
                 \hbox to \w@dth{\rightarrowfill}}%
       }\limits
\ifx\n@one\empty\else ^{\box\z@}\fi
\ifx\n@two\empty\else _{\box\@ne}\fi}}
\def\tw@@^#1{\@ifnextchar _{\two@^{#1}}{\two@^{#1}_{}}}
\def\two{\@ifnextchar ^{\tw@@}{\tw@@^{}}}
 
\def\tofr@^#1_#2{\def\n@one{#1}\def\n@two{#2}\mathrel{\setw@dth{#1}{#2}
\mathop{\vcenter{\hbox to \w@dth{\rightarrowfill}\kern-1.7ex
                 \hbox to \w@dth{\leftarrowfill}}%
       }\limits
\ifx\n@one\empty\else ^{\box\z@}\fi
\ifx\n@two\empty\else _{\box\@ne}\fi}}
\def\t@fr@^#1{\@ifnextchar_ {\tofr@^{#1}}{\tofr@^{#1}_{}}}
\def\tofro{\@ifnextchar^ {\t@fr@}{\t@fr@^{}}}
 
\def\epi{\mathop{\mathchar"221\mkern -12mu\mathchar"221}\limits}
\def\leftepi{\mathop{\mathchar"220\mkern -12mu\mathchar"220}\limits}
\def\mon{\mathop{\m@th\hbox to
      14.6\P@{\lasyb\char'51\hskip-2.1\P@$\arrext$\hss
$\mathord\rightarrow$}}\limits} 
\def\leftmono{\mathrel{\m@th\hbox to
14.6\P@{$\mathord\leftarrow$\hss$\arrext$\hskip-2.1\P@\lasyb\char'50%
}}\limits} 
\mathchardef\arrext="0200       

\setlength{\unitlength}{.01em}%
\def\settypes(#1,#2,#3){\arrowtypea#1 \arrowtypeb#2 \arrowtypec#3}
\def\settoheight#1#2{\setbox\@tempboxa\hbox{#2}#1\ht\@tempboxa\relax}%
\def\settodepth#1#2{\setbox\@tempboxa\hbox{#2}#1\dp\@tempboxa\relax}%
\def\settokens`#1`#2`#3`#4`{%
     \def\tokena{#1}\def\tokenb{#2}\def\tokenc{#3}\def\tokend{#4}}
\def\setsqparms[#1`#2`#3`#4;#5`#6]{%
\arrowtypea #1
\arrowtypeb #2
\arrowtypec #3
\arrowtyped #4
\width #5
\height #6
}
\def\setpos(#1,#2){\xpos=#1 \ypos#2}

\def\settriparms[#1`#2`#3;#4]{\settripairparms[#1`#2`#3`1`1;#4]}%
 
\def\settripairparms[#1`#2`#3`#4`#5;#6]{%
\arrowtypea #1
\arrowtypeb #2
\arrowtypec #3
\arrowtyped #4
\arrowtypee #5
\width #6
\height #6
}
 
\def\resetparms{\settripairparms[1`1`1`1`1;500]\width 500}
 
\resetparms
 
\def\mvector(#1,#2)#3{
\put(0,0){\vector(#1,#2){#3}}%
\put(0,0){\vector(#1,#2){26}}%
}
\def\evector(#1,#2)#3{{
\arrowlength #3
\put(0,0){\vector(#1,#2){\arrowlength}}%
\advance \arrowlength by-30
\put(0,0){\vector(#1,#2){\arrowlength}}%
}}
 
\def\horsize#1#2{%
\settowidth{\tempdimen}{$#2$}%
#1=\tempdimen
\divide #1 by\unitlength
}
 
\def\vertsize#1#2{%
\settoheight{\tempdimen}{$#2$}%
#1=\tempdimen
\settodepth{\tempdimen}{$#2$}%
\advance #1 by\tempdimen
\divide #1 by\unitlength
}
 
\def\putvector(#1,#2)(#3,#4)#5#6{{%
\ifnum3<\arrowtype
\putdashvector(#1,#2)(#3,#4)#5\arrowtype
\else
\ifnum\arrowtype<-3
\putdashvector(#1,#2)(#3,#4)#5\arrowtype
\else
\xpos=#1
\ypos=#2
\run=#3
\rise=#4
\arrowlength=#5
\ifnum \arrowtype<0
    \ifnum \run=0
        \advance \ypos by-\arrowlength
    \else
        \tempcounta \arrowlength
        \multiply \tempcounta by\rise
        \divide \tempcounta by\run
        \ifnum\run>0
            \advance \xpos by\arrowlength
            \advance \ypos by\tempcounta
        \else
            \advance \xpos by-\arrowlength
            \advance \ypos by-\tempcounta
        \fi
    \fi
    \multiply \arrowtype by-1
    \multiply \rise by-1
    \multiply \run by-1
\fi
\ifcase \arrowtype
\or \put(\xpos,\ypos){\vector(\run,\rise){\arrowlength}}%
\or \put(\xpos,\ypos){\mvector(\run,\rise)\arrowlength}%
\or \put(\xpos,\ypos){\evector(\run,\rise){\arrowlength}}%
\fi\fi\fi
}}
 
\def\putsplitvector(#1,#2)#3#4{
\xpos #1
\ypos #2
\arrowtype #4
\halflength #3
\arrowlength #3
\gap 140
\advance \halflength by-\gap
\divide \halflength by2
\ifnum\arrowtype>0
   \ifcase \arrowtype
   \or \put(\xpos,\ypos){\line(0,-1){\halflength}}%
       \advance\ypos by-\halflength
       \advance\ypos by-\gap
       \put(\xpos,\ypos){\vector(0,-1){\halflength}}%
   \or \put(\xpos,\ypos){\line(0,-1)\halflength}%
       \put(\xpos,\ypos){\vector(0,-1)3}%
       \advance\ypos by-\halflength
       \advance\ypos by-\gap
       \put(\xpos,\ypos){\vector(0,-1){\halflength}}%
   \or \put(\xpos,\ypos){\line(0,-1)\halflength}%
       \advance\ypos by-\halflength
       \advance\ypos by-\gap
       \put(\xpos,\ypos){\evector(0,-1){\halflength}}%
   \fi
\else \arrowtype=-\arrowtype
   \ifcase\arrowtype
   \or \advance \ypos by-\arrowlength
       \put(\xpos,\ypos){\line(0,1){\halflength}}%
       \advance\ypos by\halflength
       \advance\ypos by\gap
       \put(\xpos,\ypos){\vector(0,1){\halflength}}%
   \or \advance \ypos by-\arrowlength
       \put(\xpos,\ypos){\line(0,1)\halflength}%
       \put(\xpos,\ypos){\vector(0,1)3}%
       \advance\ypos by\halflength
       \advance\ypos by\gap
       \put(\xpos,\ypos){\vector(0,1){\halflength}}%
   \or \advance \ypos by-\arrowlength
       \put(\xpos,\ypos){\line(0,1)\halflength}%
       \advance\ypos by\halflength
       \advance\ypos by\gap
       \put(\xpos,\ypos){\evector(0,1){\halflength}}%
   \fi
\fi
}
 
\def\putmorphism(#1)(#2,#3)[#4`#5`#6]#7#8#9{{%
\run #2
\rise #3
\ifnum\rise=0
  \puthmorphism(#1)[#4`#5`#6]{#7}{#8}#9%
\else\ifnum\run=0
  \putvmorphism(#1)[#4`#5`#6]{#7}{#8}#9%
\else
\setpos(#1)%
\arrowlength #7
\arrowtype #8
\ifnum\run=0
\else\ifnum\rise=0
\else
\ifnum\run>0
    \coefa=1
\else
   \coefa=-1
\fi
\ifnum\arrowtype>0
   \coefb=0
   \coefc=-1
\else
   \coefb=\coefa
   \coefc=1
   \arrowtype=-\arrowtype
\fi
\width=2
\multiply \width by\run
\divide \width by\rise
\ifnum \width<0  \width=-\width\fi
\advance\width by60
\if l#9 \width=-\width\fi
\putbox(\xpos,\ypos){#4}
{\multiply \coefa by\arrowlength
\advance\xpos by\coefa
\multiply \coefa by\rise
\divide \coefa by\run
\advance \ypos by\coefa
\putbox(\xpos,\ypos){#5} }%
{\multiply \coefa by\arrowlength
\divide \coefa by2
\advance \xpos by\coefa
\advance \xpos by\width
\multiply \coefa by\rise
\divide \coefa by\run
\advance \ypos by\coefa
\if l#9%
   \putrbox(\xpos,\ypos){#6}%
\else\if r#9%
   \putlbox(\xpos,\ypos){#6}%
\fi\fi }%
{\multiply \rise by-\coefc
\multiply \run by-\coefc
\multiply \coefb by\arrowlength
\advance \xpos by\coefb
\multiply \coefb by\rise
\divide \coefb by\run
\advance \ypos by\coefb
\multiply \coefc by70
\advance \ypos by\coefc
\multiply \coefc by\run
\divide \coefc by\rise
\advance \xpos by\coefc
\multiply \coefa by140
\multiply \coefa by\run
\divide \coefa by\rise
\advance \arrowlength by\coefa
\ifcase\arrowtype
\or \put(\xpos,\ypos){\vector(\run,\rise){\arrowlength}}%
\or \put(\xpos,\ypos){\mvector(\run,\rise){\arrowlength}}%
\or \put(\xpos,\ypos){\evector(\run,\rise){\arrowlength}}%
\fi}\fi\fi\fi\fi}}

\newcount\numbdashes \newcount\lengthdash \newcount\increment
 
\def\howmanydashes{
\numbdashes=\arrowlength \lengthdash=40
\divide\numbdashes by \lengthdash
\lengthdash=\arrowlength
\divide\lengthdash by \numbdashes
\increment=\lengthdash
\multiply\lengthdash by 3
\divide\lengthdash by 5
}
 
\def\putdashvector(#1)(#2,#3)#4#5{%
\ifnum#3=0 \putdashhvector(#1){#4}#5
\else
\ifnum#2=0
\putdashvvector(#1){#4}#5\fi\fi}
 
\def\putdashhvector(#1,#2)#3#4{{%
\arrowlength=#3 \howmanydashes
\multiput(#1,#2)(\increment,0){\numbdashes}%
{\vrule height .4pt width \lengthdash\unitlength}
\arrowtype=#4 \xpos=#1
\ifnum\arrowtype<0 \advance\arrowtype by 7 \fi
\ifcase\arrowtype
\or \advance\xpos by 10
    \put(\xpos,#2){\vector(-1,0){\lengthdash}}
    \advance\xpos by 40
    \put(\xpos,#2){\vector(-1,0){\lengthdash}}
\or \advance \xpos by 10
    \put(\xpos,#2){\vector(-1,0){\lengthdash}}
    \advance\xpos by  \arrowlength
    \advance\xpos by  -50
    \put(\xpos,#2){\vector(-1,0){\lengthdash}}
\or \advance\xpos by 10
    \put(\xpos,#2){\vector(-1,0){\lengthdash}}
\or \advance\xpos by \arrowlength
    \advance\xpos by -\lengthdash
    \put(\xpos,#2){\vector(1,0){\lengthdash}}
\or {\advance\xpos by 10
    \put(\xpos,#2){\vector(1,0){\lengthdash}}}
    \advance\xpos by \arrowlength
    \advance\xpos by -\lengthdash
    \put(\xpos,#2){\vector(1,0){\lengthdash}}
\or \advance\xpos by \arrowlength
    \advance\xpos by -\lengthdash
    \put(\xpos,#2){\vector(1,0){\lengthdash}}
    \advance\xpos by -40
    \put(\xpos,#2){\vector(1,0){\lengthdash}}
   \fi
}}
 
\def\putdashvvector(#1,#2)#3#4{{%
\arrowlength=#3 \howmanydashes
\ypos=#2 \advance\ypos by -\arrowlength
\multiput(#1,#2)(0,\increment){\numbdashes}%
    {\vrule width .4pt height \lengthdash\unitlength}
\arrowtype=#4 \ypos=#2
\ifnum\arrowtype<0 \advance\arrowtype by 7 \fi
\ifcase\arrowtype
\or \advance\ypos by \arrowlength \advance\ypos by -40
    \put(#1,\ypos){\vector(0,1){\lengthdash}}
    \advance\ypos by -40
    \put(#1,\ypos){\vector(0,1){\lengthdash}}
\or \advance\ypos by 10
    \put(#1,\ypos){\vector(0,1){\lengthdash}}
    \advance\ypos by \arrowlength \advance\ypos by -40
    \put(#1,\ypos){\vector(0,1){\lengthdash}}
\or \advance\ypos by \arrowlength \advance\ypos by -40
    \put(#1,\ypos){\vector(0,1){\lengthdash}}
\or \advance\ypos by 10
    \put(#1,\ypos){\vector(0,-1){\lengthdash}}
\or \advance\ypos by 10
    \put(#1,\ypos){\vector(0,-1){\lengthdash}}
    \advance\ypos by \arrowlength \advance\ypos by -40
    \put(#1,\ypos){\vector(0,-1){\lengthdash}}
\or \advance\ypos by 10
    \put(#1,\ypos){\vector(0,-1){\lengthdash}}
    \advance\ypos by 40
    \put(#1,\ypos){\vector(0,-1){\lengthdash}}
\fi
}}
 
\def\puthmorphism(#1,#2)[#3`#4`#5]#6#7#8{{%
\xpos #1
\ypos #2
\width #6
\arrowlength #6
\arrowtype=#7
\putbox(\xpos,\ypos){#3\vphantom{#4}}%
{\advance \xpos by\arrowlength
\putbox(\xpos,\ypos){\vphantom{#3}#4}}%
\horsize{\tempcounta}{#3}%
\horsize{\tempcountb}{#4}%
\divide \tempcounta by2
\divide \tempcountb by2
\advance \tempcounta by30
\advance \tempcountb by30
\advance \xpos by\tempcounta
\advance \arrowlength by-\tempcounta
\advance \arrowlength by-\tempcountb
\putvector(\xpos,\ypos)(1,0)\arrowlength\arrowtype
\divide \arrowlength by2
\advance \xpos by\arrowlength
\vertsize{\tempcounta}{#5}%
\divide\tempcounta by2
\advance \tempcounta by20
\if a#8 %
   \advance \ypos by\tempcounta
   \putbox(\xpos,\ypos){#5}%
\else
   \advance \ypos by-\tempcounta
   \putbox(\xpos,\ypos){#5}%
\fi}}
 
\def\putvmorphism(#1,#2)[#3`#4`#5]#6#7#8{{%
\xpos #1
\ypos #2
\arrowlength #6
\arrowtype #7
\settowidth{\xlen}{$#5$}%
\putbox(\xpos,\ypos){#3}%
{\advance \ypos by-\arrowlength
\putbox(\xpos,\ypos){#4}}%
{\advance\arrowlength by-140
\advance \ypos by-70
\ifdim\xlen>0pt
   \if m#8%
      \putsplitvector(\xpos,\ypos)\arrowlength\arrowtype
   \else
   \putvector(\xpos,\ypos)(0,-1)\arrowlength\arrowtype
   \fi
\else
   \putvector(\xpos,\ypos)(0,-1)\arrowlength\arrowtype
\fi}%
\ifdim\xlen>0pt
   \divide \arrowlength by2
   \advance\ypos by-\arrowlength
   \if l#8%
      \advance \xpos by-40
      \putrbox(\xpos,\ypos){#5}%
   \else\if r#8%
      \advance \xpos by40
      \putlbox(\xpos,\ypos){#5}%
   \else
      \putbox(\xpos,\ypos){#5}%
   \fi\fi
\fi
}}
 
\def\putsquarep<#1>(#2)[#3;#4`#5`#6`#7]{{%
\setsqparms[#1]%
\setpos(#2)%
\settokens`#3`%
\puthmorphism(\xpos,\ypos)[\tokenc`\tokend`{#7}]{\width}{\arrowtyped}b%
\advance\ypos by \height
\puthmorphism(\xpos,\ypos)[\tokena`\tokenb`{#4}]{\width}{\arrowtypea}a%
\putvmorphism(\xpos,\ypos)[``{#5}]{\height}{\arrowtypeb}l%
\advance\xpos by \width
\putvmorphism(\xpos,\ypos)[``{#6}]{\height}{\arrowtypec}r%
}}
 
\def\putsquare{\@ifnextchar <{\putsquarep}{\putsquarep%
   <\arrowtypea`\arrowtypeb`\arrowtypec`\arrowtyped;\width`\height>}}
\def\square{\@ifnextchar< {\squarep}{\squarep
   <\arrowtypea`\arrowtypeb`\arrowtypec`\arrowtyped;\width`\height>}}
\def\squarep<#1>[#2`#3`#4`#5;#6`#7`#8`#9]{{
\setsqparms[#1]
\diagram
\putsquarep<\arrowtypea`\arrowtypeb`\arrowtypec`
\arrowtyped;\width`\height>
(0,0)[#2`#3`#4`{#5};#6`#7`#8`{#9}]
\enddiagram
}}                                                 
\def\putptrianglep<#1>(#2,#3)[#4`#5`#6;#7`#8`#9]{{%
\settriparms[#1]%
\xpos=#2 \ypos=#3
\advance\ypos by \height
\puthmorphism(\xpos,\ypos)[#4`#5`{#7}]{\height}{\arrowtypea}a%
\putvmorphism(\xpos,\ypos)[`#6`{#8}]{\height}{\arrowtypeb}l%
\advance\xpos by\height
\putmorphism(\xpos,\ypos)(-1,-1)[``{#9}]{\height}{\arrowtypec}r%
}}
 
\def\putptriangle{\@ifnextchar <{\putptrianglep}{\putptrianglep
   <\arrowtypea`\arrowtypeb`\arrowtypec;\height>}}
\def\ptriangle{\@ifnextchar <{\ptrianglep}{\ptrianglep
   <\arrowtypea`\arrowtypeb`\arrowtypec;\height>}}
\def\ptrianglep<#1>[#2`#3`#4;#5`#6`#7]{{
\settriparms[#1]
\diagram
\putptrianglep<\arrowtypea`\arrowtypeb`
\arrowtypec;\height>
(0,0)[#2`#3`#4;#5`#6`{#7}]
\enddiagram
}}                                            
 
\def\putqtrianglep<#1>(#2,#3)[#4`#5`#6;#7`#8`#9]{{%
\settriparms[#1]%
\xpos=#2 \ypos=#3
\advance\ypos by\height
\puthmorphism(\xpos,\ypos)[#4`#5`{#7}]{\height}{\arrowtypea}a%
\putmorphism(\xpos,\ypos)(1,-1)[``{#8}]{\height}{\arrowtypeb}l%
\advance\xpos by\height
\putvmorphism(\xpos,\ypos)[`#6`{#9}]{\height}{\arrowtypec}r%
}}
 
\def\putqtriangle{\@ifnextchar <{\putqtrianglep}{\putqtrianglep
   <\arrowtypea`\arrowtypeb`\arrowtypec;\height>}}
\def\qtriangle{\@ifnextchar <{\qtrianglep}{\qtrianglep
   <\arrowtypea`\arrowtypeb`\arrowtypec;\height>}}
\def\qtrianglep<#1>[#2`#3`#4;#5`#6`#7]{{
\settriparms[#1]
\width=\height                                
\diagram
\putqtrianglep<\arrowtypea`\arrowtypeb`
\arrowtypec;\height>
(0,0)[#2`#3`#4;#5`#6`{#7}]
\enddiagram
}}
 
\def\putdtrianglep<#1>(#2,#3)[#4`#5`#6;#7`#8`#9]{{%
\settriparms[#1]%
\xpos=#2 \ypos=#3
\puthmorphism(\xpos,\ypos)[#5`#6`{#9}]{\height}{\arrowtypec}b%
\advance\xpos by \height \advance\ypos by\height
\putmorphism(\xpos,\ypos)(-1,-1)[``{#7}]{\height}{\arrowtypea}l%
\putvmorphism(\xpos,\ypos)[#4``{#8}]{\height}{\arrowtypeb}r%
}}
 
\def\putdtriangle{\@ifnextchar <{\putdtrianglep}{\putdtrianglep
   <\arrowtypea`\arrowtypeb`\arrowtypec;\height>}}
\def\dtriangle{\@ifnextchar <{\dtrianglep}{\dtrianglep
   <\arrowtypea`\arrowtypeb`\arrowtypec;\height>}}
\def\dtrianglep<#1>[#2`#3`#4;#5`#6`#7]{{
\settriparms[#1]
\width=\height                                
\diagram
\putdtrianglep<\arrowtypea`\arrowtypeb`
\arrowtypec;\height>
(0,0)[#2`#3`#4;#5`#6`{#7}]
\enddiagram
}}
 
\def\putbtrianglep<#1>(#2,#3)[#4`#5`#6;#7`#8`#9]{{%
\settriparms[#1]%
\xpos=#2 \ypos=#3
\puthmorphism(\xpos,\ypos)[#5`#6`{#9}]{\height}{\arrowtypec}b%
\advance\ypos by\height
\putmorphism(\xpos,\ypos)(1,-1)[``{#8}]{\height}{\arrowtypeb}r%
\putvmorphism(\xpos,\ypos)[#4``{#7}]{\height}{\arrowtypea}l%
}}
 
\def\putbtriangle{\@ifnextchar <{\putbtrianglep}{\putbtrianglep
   <\arrowtypea`\arrowtypeb`\arrowtypec;\height>}}
\def\btriangle{\@ifnextchar <{\btrianglep}{\btrianglep
   <\arrowtypea`\arrowtypeb`\arrowtypec;\height>}}
\def\btrianglep<#1>[#2`#3`#4;#5`#6`#7]{{
\settriparms[#1]
\width=\height                               
\diagram
\putbtrianglep<\arrowtypea`\arrowtypeb`
\arrowtypec;\height>
(0,0)[#2`#3`#4;#5`#6`{#7}]
\enddiagram
}}
 
\def\putAtrianglep<#1>(#2,#3)[#4`#5`#6;#7`#8`#9]{{%
\settriparms[#1]%
\xpos=#2 \ypos=#3
{\multiply \height by2
\puthmorphism(\xpos,\ypos)[#5`#6`{#9}]{\height}{\arrowtypec}b}%
\advance\xpos by\height \advance\ypos by\height
\putmorphism(\xpos,\ypos)(-1,-1)[#4``{#7}]{\height}{\arrowtypea}l%
\putmorphism(\xpos,\ypos)(1,-1)[``{#8}]{\height}{\arrowtypeb}r%
}}
 
\def\putAtriangle{\@ifnextchar <{\putAtrianglep}{\putAtrianglep
   <\arrowtypea`\arrowtypeb`\arrowtypec;\height>}}
\def\Atriangle{\@ifnextchar <{\Atrianglep}{\Atrianglep
   <\arrowtypea`\arrowtypeb`\arrowtypec;\height>}}
\def\Atrianglep<#1>[#2`#3`#4;#5`#6`#7]{{
\settriparms[#1]
\width=\height                                     
\diagram
\putAtrianglep<\arrowtypea`\arrowtypeb`
\arrowtypec;\height>
(0,0)[#2`#3`#4;#5`#6`{#7}]
\enddiagram
}}
 
\def\putAtrianglepairp<#1>(#2)[#3;#4`#5`#6`#7`#8]{{%
\settripairparms[#1]%
\setpos(#2)%
\settokens`#3`%
\puthmorphism(\xpos,\ypos)[\tokenb`\tokenc`{#7}]{\height}{\arrowtyped}b%
\advance\xpos by\height
\puthmorphism(\xpos,\ypos)[\phantom{\tokenc}`\tokend`{#8}]%
{\height}{\arrowtypee}b%
\advance\ypos by\height
\putmorphism(\xpos,\ypos)(-1,-1)[\tokena``{#4}]{\height}{\arrowtypea}l%
\putvmorphism(\xpos,\ypos)[``{#5}]{\height}{\arrowtypeb}m%
\putmorphism(\xpos,\ypos)(1,-1)[``{#6}]{\height}{\arrowtypec}r%
}}
 
\def\putAtrianglepair{\@ifnextchar <{\putAtrianglepairp}{\putAtrianglepairp%
   <\arrowtypea`\arrowtypeb`\arrowtypec`\arrowtyped`\arrowtypee;\height>}}
\def\Atrianglepair{\@ifnextchar <{\Atrianglepairp}{\Atrianglepairp%
   <\arrowtypea`\arrowtypeb`\arrowtypec`\arrowtyped`\arrowtypee;\height>}}
 
\def\Atrianglepairp<#1>[#2;#3`#4`#5`#6`#7]{{
\settripairparms[#1]
\settokens`#2`
\width=\height                                
\diagram
\putAtrianglepairp                            
<\arrowtypea`\arrowtypeb`\arrowtypec`
\arrowtyped`\arrowtypee;\height>
(0,0)[{#2};#3`#4`#5`#6`{#7}]
\enddiagram
}}
 
\def\putVtrianglep<#1>(#2,#3)[#4`#5`#6;#7`#8`#9]{{%
\settriparms[#1]%
\xpos=#2 \ypos=#3
\advance\ypos by\height
{\multiply\height by2
\puthmorphism(\xpos,\ypos)[#4`#5`{#7}]{\height}{\arrowtypea}a}%
\putmorphism(\xpos,\ypos)(1,-1)[`#6`{#8}]{\height}{\arrowtypeb}l%
\advance\xpos by\height
\advance\xpos by\height
\putmorphism(\xpos,\ypos)(-1,-1)[``{#9}]{\height}{\arrowtypec}r%
}}
 
\def\putVtriangle{\@ifnextchar <{\putVtrianglep}{\putVtrianglep
   <\arrowtypea`\arrowtypeb`\arrowtypec;\height>}}
\def\Vtriangle{\@ifnextchar <{\Vtrianglep}{\Vtrianglep
   <\arrowtypea`\arrowtypeb`\arrowtypec;\height>}}
\def\Vtrianglep<#1>[#2`#3`#4;#5`#6`#7]{{
\settriparms[#1]
\width=\height                                 
\diagram
\putVtrianglep<\arrowtypea`\arrowtypeb`
\arrowtypec;\height>
(0,0)[#2`#3`#4;#5`#6`{#7}]
\enddiagram
}}
 
\def\putVtrianglepairp<#1>(#2)[#3;#4`#5`#6`#7`#8]{{
\settripairparms[#1]%
\setpos(#2)%
\settokens`#3`%
\advance\ypos by\height
\putmorphism(\xpos,\ypos)(1,-1)[`\tokend`{#6}]{\height}{\arrowtypec}l%
\puthmorphism(\xpos,\ypos)[\tokena`\tokenb`{#4}]{\height}{\arrowtypea}a%
\advance\xpos by\height
\puthmorphism(\xpos,\ypos)[\phantom{\tokenb}`\tokenc`{#5}]%
{\height}{\arrowtypeb}a%
\putvmorphism(\xpos,\ypos)[``{#7}]{\height}{\arrowtyped}m%
\advance\xpos by\height
\putmorphism(\xpos,\ypos)(-1,-1)[``{#8}]{\height}{\arrowtypee}r%
}}
 
\def\putVtrianglepair{\@ifnextchar <{\putVtrianglepairp}{\putVtrianglepairp%
    <\arrowtypea`\arrowtypeb`\arrowtypec`\arrowtyped`\arrowtypee;\height>}}
\def\Vtrianglepair{\@ifnextchar <{\Vtrianglepairp}{\Vtrianglepairp%
    <\arrowtypea`\arrowtypeb`\arrowtypec`\arrowtyped`\arrowtypee;\height>}}
\def\Vtrianglepairp<#1>[#2;#3`#4`#5`#6`#7]{{
\settripairparms[#1]
\settokens`#2`
\diagram
\putVtrianglepairp                             
<\arrowtypea`\arrowtypeb`\arrowtypec`
\arrowtyped`\arrowtypee;\height>
(0,0)[{#2};#3`#4`#5`#6`{#7}]
\enddiagram
}}

\def\putCtrianglep<#1>(#2,#3)[#4`#5`#6;#7`#8`#9]{{%
\settriparms[#1]%
\xpos=#2 \ypos=#3
\advance\ypos by\height
\putmorphism(\xpos,\ypos)(1,-1)[``{#9}]{\height}{\arrowtypec}l%
\advance\xpos by\height
\advance\ypos by\height
\putmorphism(\xpos,\ypos)(-1,-1)[#4`#5`{#7}]{\height}{\arrowtypea}l%
{\multiply\height by 2
\putvmorphism(\xpos,\ypos)[`#6`{#8}]{\height}{\arrowtypeb}r}%
}}
 
\def\putCtriangle{\@ifnextchar <{\putCtrianglep}{\putCtrianglep
    <\arrowtypea`\arrowtypeb`\arrowtypec;\height>}}
\def\Ctriangle{\@ifnextchar <{\Ctrianglep}{\Ctrianglep
    <\arrowtypea`\arrowtypeb`\arrowtypec;\height>}}
\def\Ctrianglep<#1>[#2`#3`#4;#5`#6`#7]{{
\settriparms[#1]
\width=\height                               
\diagram
\putCtrianglep<\arrowtypea`\arrowtypeb`
\arrowtypec;\height>
(0,0)[#2`#3`#4;#5`#6`{#7}]
\enddiagram
}}                                           
\def\putDtrianglep<#1>(#2,#3)[#4`#5`#6;#7`#8`#9]{{%
\settriparms[#1]%
\xpos=#2 \ypos=#3
\advance\xpos by\height \advance\ypos by\height
\putmorphism(\xpos,\ypos)(-1,-1)[``{#9}]{\height}{\arrowtypec}r%
\advance\xpos by-\height \advance\ypos by\height
\putmorphism(\xpos,\ypos)(1,-1)[`#5`{#8}]{\height}{\arrowtypeb}r%
{\multiply\height by 2
\putvmorphism(\xpos,\ypos)[#4`#6`{#7}]{\height}{\arrowtypea}l}%
}}
 
\def\putDtriangle{\@ifnextchar <{\putDtrianglep}{\putDtrianglep
    <\arrowtypea`\arrowtypeb`\arrowtypec;\height>}}
\def\Dtriangle{\@ifnextchar <{\Dtrianglep}{\Dtrianglep
   <\arrowtypea`\arrowtypeb`\arrowtypec;\height>}}
\def\Dtrianglep<#1>[#2`#3`#4;#5`#6`#7]{{
\settriparms[#1]
\width=\height                              
\diagram
\putDtrianglep<\arrowtypea`\arrowtypeb`
\arrowtypec;\height>
(0,0)[#2`#3`#4;#5`#6`{#7}]
\enddiagram
}}                                          
\def\setrecparms[#1`#2]{\width=#1 \height=#2}%
 
\def\recursep<#1`#2>[#3;#4`#5`#6`#7`#8]{{\m@th
\width=#1 \height=#2
\settokens`#3`
\settowidth{\tempdimen}{$\tokena$}
\ifdim\tempdimen=0pt
  \savebox{\tempboxa}{\hbox{$\tokenb$}}%
  \savebox{\tempboxb}{\hbox{$\tokend$}}%
  \savebox{\tempboxc}{\hbox{$#6$}}%
\else
  \savebox{\tempboxa}{\hbox{$\hbox{$\tokena$}\times\hbox{$\tokenb$}$}}%
  \savebox{\tempboxb}{\hbox{$\hbox{$\tokena$}\times\hbox{$\tokend$}$}}%
  \savebox{\tempboxc}{\hbox{$\hbox{$\tokena$}\times\hbox{$#6$}$}}%
\fi
\ypos=\height
\divide\ypos by 2
\xpos=\ypos
\advance\xpos by \width
\bfig
\putCtrianglep<-1`1`1;\ypos>(0,0)[`\tokenc`;#5`#6`{#7}]%
\puthmorphism(\ypos,0)[\tokend`\usebox{\tempboxb}`{#8}]{\width}{-1}b%
\puthmorphism(\ypos,\height)[\tokenb`\usebox{\tempboxa}`{#4}]{\width}{-1}a%
\advance\ypos by \width
\putvmorphism(\ypos,\height)[``\usebox{\tempboxc}]{\height}1r%
\efig
}}
 
\def\recurse{\@ifnextchar <{\recursep}{\recursep<\width`\height>}}
 
\def\puttwohmorphisms(#1,#2)[#3`#4;#5`#6]#7#8#9{{%
%
\puthmorphism(#1,#2)[#3`#4`]{#7}0a
\ypos=#2
\advance\ypos by 20
\puthmorphism(#1,\ypos)[\phantom{#3}`\phantom{#4}`#5]{#7}{#8}a
\advance\ypos by -40
\puthmorphism(#1,\ypos)[\phantom{#3}`\phantom{#4}`#6]{#7}{#9}b
}}
 
\def\puttwovmorphisms(#1,#2)[#3`#4;#5`#6]#7#8#9{{%
%
%
\putvmorphism(#1,#2)[#3`#4`]{#7}0a
\xpos=#1
\advance\xpos by -20
\putvmorphism(\xpos,#2)[\phantom{#3}`\phantom{#4}`#5]{#7}{#8}l
\advance\xpos by 40
\putvmorphism(\xpos,#2)[\phantom{#3}`\phantom{#4}`#6]{#7}{#9}r
}}
 
\def\puthcoequalizer(#1)[#2`#3`#4;#5`#6`#7]#8#9{{%
%
\setpos(#1)%
\puttwohmorphisms(\xpos,\ypos)[#2`#3;#5`#6]{#8}11%
\advance\xpos by #8
\puthmorphism(\xpos,\ypos)[\phantom{#3}`#4`#7]{#8}1{#9}
}}
 
\def\putvcoequalizer(#1)[#2`#3`#4;#5`#6`#7]#8#9{{%
%
%
\setpos(#1)%
\puttwovmorphisms(\xpos,\ypos)[#2`#3;#5`#6]{#8}11%
\advance\ypos by -#8
\putvmorphism(\xpos,\ypos)[\phantom{#3}`#4`#7]{#8}1{#9}
}}
 
\def\putthreehmorphisms(#1)[#2`#3;#4`#5`#6]#7(#8)#9{{%
\setpos(#1) \settypes(#8)
\if a#9 %
     \vertsize{\tempcounta}{#5}%
     \vertsize{\tempcountb}{#6}%
     \ifnum \tempcounta<\tempcountb \tempcounta=\tempcountb \fi
\else
     \vertsize{\tempcounta}{#4}%
     \vertsize{\tempcountb}{#5}%
     \ifnum \tempcounta<\tempcountb \tempcounta=\tempcountb \fi
\fi
\advance \tempcounta by 60
\puthmorphism(\xpos,\ypos)[#2`#3`#5]{#7}{\arrowtypeb}{#9}
\advance\ypos by \tempcounta
\puthmorphism(\xpos,\ypos)[\phantom{#2}`\phantom{#3}`#4]{#7}{\arrowtypea}{#9}
\advance\ypos by -\tempcounta \advance\ypos by -\tempcounta
\puthmorphism(\xpos,\ypos)[\phantom{#2}`\phantom{#3}`#6]{#7}{\arrowtypec}{#9}
}}
 
\def\setarrowtoks[#1`#2`#3`#4`#5`#6]{%
\def\toka{#1}
\def\tokb{#2}
\def\tokc{#3}
\def\tokd{#4}
\def\toke{#5}
\def\tokf{#6}
}
\def\hex{\@ifnextchar <{\hexp}{\hexp<1000`400>}}
\def\hexp<#1`#2>[#3`#4`#5`#6`#7`#8;#9]{%
\setarrowtoks[#9]
\yext=#2 \advance \yext by #2
\xext=#1 \advance\xext by \yext
\bfig
\putCtriangle<-1`0`1;#2>(0,0)[`#5`;\tokb``\tokd]
\xext=#1 \yext=#2 \advance \yext by #2
\putsquare<1`0`0`1;\xext`\yext>(#2,0)[#3`#4`#7`#8;\toka```\tokf]
\advance \xext by #2
\putDtriangle<0`1`-1;#2>(\xext,0)[`#6`;`\tokc`\toke]
\efig
}
\makeatother

\newpage
\pagenumbering{arabic}
\spacing{1.3}
\chapter*{Introduction.}
\addcontentsline{toc}{chapter}{\protect\numberline{}{Introduction}}
\pagestyle{myheadings}
\markboth{INTRODUCTION}{INTRODUCTION}
         
Elie Cartan's method of equivalence is a natural development of the Felix Klein  Erlangen program (1872), which describes geometry as the study of invariants of group actions on geometric objects.  Cartan formulated the problem of the equivalence of submanifolds under a group of transformations and  introduced the method of moving frames, which allows one to construct differential invariants under a group action \cite{C37}. The functional relations among the invariants provide a key to the classification of submanifolds under a prescribed group of transformations. Classically, a  moving frame is an equivariant map from the space of   submanifolds (or more rigorously, from the corresponding  jet bundle) to the  bundle of frames. Exterior  differentiation of this map produces an invariant coframe, which is used to construct a  number of differential invariants sufficient to solve the equivalence problem.  Considering  moving frame constructions on homogeneous spaces,  Griffiths \cite{Gr74} and Green \cite{Green78} observed  that a moving frame
can be viewed as an equivariant map from the space of  submanifolds to the group itself. As  pointed out in \cite{Gr74}, one of the  classical  moving  frames, the  Fr\'enet frame, is in fact  a map from the  space of curves to the Euclidean group.
Adopting this observation as a general definition of a  moving frame,  Fels and Olver  \cite{FO98},  \cite{FO99} generalized the method  for arbitrary, not necessarily 
transitive, 
finite-dimensional   Lie group actions on a manifold introducing, for the first time, a completely algorithmic  way for their construction.  

In the first chapter we give an overview of the Cartan's solution to the problem of the equivalence of submanifolds. We also describe a general algorithm for construction
 of the moving frames and differential invariants developed  
by Fels and Olver  \cite{FO99}.  Following this method, however, 
one might have to prolong 
the action to the jet spaces of high order before  obtaining any invariants, while 
 the earlier methods  \cite{Gr74}, \cite{Green78} allow one to construct invariants 
order 
by order. We combine the advantages of both approaches  in the recursive algorithm 
presented in Chapter~2. Not surprisingly, the construction of moving frames and 
differential invariants is simpler when the acting group has fewer parameters. 
Thus, it is desirable  to use  
the results obtained  for a subgroup $H\subset G$
to construct a moving frame and differential invariants for the entire group $G$.
The inductive algorithm  from Chapter~2 allows us to perform this for the groups 
that factor as a product. 
Using this algorithm, one obtains at the same time the relations among the invariants 
 of group $G$ and its subgroup $H$.  
An illustrative  example  is   induction from the Euclidean action on 
plane curves to the special affine
action, and then to the action of the entire  projective group. As a by-product,
one obtains 
the expression of the affine invariants  in terms of the Euclidean ones and the 
projective invariants in terms of the affine ones. The actions of all three  groups 
play an important role in computer image processing \cite{Fu94}, \cite{ST94}.

Equipped with all of these tools we approach the problem of the equivalence and symmetry of 
polynomials  under linear changes of variables in Chapter~3. 
In fact, applications to classical 
invariant theory have served as the initial motivation for the development of the algorithms from Chapter~2.  Two polynomials are said to be equivalent if there is a linear change of variables that transforms one into the other. The group of symmetries of a polynomial consists of all linear changes of 
variables that leave the polynomial  unchanged. It is desirable  
to describe the classes of equivalent polynomials and classify the corresponding symmetry groups. 
We concentrate on polynomials over complex numbers although  indicate how to adopt 
the results to real polynomials.
This problem has been traditionally approached by algebraic or algebraic-geometry tools
 \cite{GraceYoung03}, \cite{Gur64}, \cite{Kraft87}, \cite{VP89}.

Inspired by the ideas of Olver's book \cite{O99} we approach the problem using a
 differential geometric method of moving frames.
We consider the graph of a polynomial $F(\x)$ in $m$ variables as a submanifold in 
$(m+1)$-dimensional complex (or real) space therefore reducing the question 
to equivalence problems for  submanifolds.  
In the polynomial  case, differential invariants can be chosen to be
 rational functions,  and the polynomial  relations among them can be found via
 elimination algorithms based on Gr\"obner basis computations. Hence, the  problem
 of the equivalence and symmetry of polynomials can be completely solved, at least 
in theory. 
In practice, however, we are confronted with the complexity of Gr\"obner basis
 computations, which significantly limits our ability to solve specific problems. 
We start with the simplest case of polynomials in two variables, reproducing results from the paper by  Peter Olver and myself \cite{BO00}. We provide a {\sc Maple} code which determines the dimension of the symmetry group of a given polynomial and in the case
when  the symmetry group is finite, computes it explicitly. 
Computations become challenging, even  in the next case of polynomials in three variables. In order to construct a complete set of   differential invariants  
we apply the recursive 
algorithm from Chapter~2.  In some cases, the  
relations among the invariants are successfully computed via algorithms based on Gr\"obner basis computations, 
while in the other cases, this computation does not seem to be feasible. 
Nevertheless, some interesting new  results were obtained, such as  a classification  of ternary cubics and their groups 
of symmetries, and  necessary and sufficient conditions for a homogeneous polynomial in three variables to be equivalent to $x^n+y^n+z^n.$   

Cartan's method for solving the problem of equivalence for submanifolds 
was formulated in the category of smooth manifolds. Hence, its direct  applicability   
is  restricted  to polynomials over complex or real numbers.  It would be a
worthwhile  and interesting project to reformulate the method of moving frames 
in the algebraic-geometry language,  so that it can be applied to the problem of the equivalence 
and symmetry of algebraic varieties over fields of arbitrary  characteristics.

\pagestyle{headings}
\chapter{The Equivalence Problem for Submanifolds.}
     	
\label{prelim}
Two manifolds are said to be equivalent under a  transformation group if one can be 
mapped to the other by an element of the group. Symmetry can be considered as self-equivalence. The group of the symmetries, or the isotropy group, of a  submanifold consists of the transformations which map the submanifold to itself.
Given a transformation group one would like to find classes of equivalent
manifolds and
 to describe the  group of symmetries of a given submanifold.

A local  solution to this problem for a group acting on a homogeneous space was presented by E. Cartan in \cite{C37} and is known as  the method of moving frames. 
Classically, a  moving frame is an equivariant map from the space of   submanifolds (or more rigorously from the corresponding  jet bundle) to the  bundle of frames. Exterior  differentiation of this map produces an invariant coframe, which is used to construct a  number of differential invariants sufficient to solve the equivalence problem.

The method of moving frames  was  generalized  by Fels and Olver  \cite{FO99}  for arbitrary (not necessarily transitive) finite-dimensional   Lie group actions on a manifold.  It relies on the definition of  moving frame  as  an 
equivariant map from the space of  submanifolds to the group itself, which could be also found in Griffiths \cite{Gr74} and Green \cite{Green78}. 
As  pointed out by Griffiths in \cite{Gr74}, one of the  classical  moving  frames, the  Fr\'enet frame, is in fact  a map from the  space of curves to the Euclidean group.  

We remark that the problem of the equivalence for  submanifolds  is one of  many  problems which can be reformulated in terms of an exterior differential system and reduced to 
a question  of its integrability \cite{BCGGG91}. Other important examples include the  equivalence and symmetry problems  for differential equations and  for variational functionals. Cartan's method consists of rewriting the problem in terms  of differential forms  in a way that it becomes intrinsic with respect to group action. Then the exterior differentiation of this system  produces  the sufficient number of differential invariants. The equivalence problem can be solved by examining the functional relations among these invariants. Although not completely algorithmic this approach has led to solving numerous equivalence problems e.g., \cite{C35}, \cite{C53}, \cite{Ga89}, \cite{Kamran89} and  
\cite{O86} including  problems which involve  infinite dimensional pseudo-groups of transformation.

 In what follows we  consider the local  problem of equivalence and symmetry under 
finite dimensional Lie group of transformation on a smooth manifold.      
In the first section we review  some  basic  definitions and results about Lie group 
actions which we will use  later.  

\section {Lie Group Actions on Manifolds.}
\begin{df} A smooth map $w:G\times M\rightarrow M$ defines an {\it action} 
of a group $G$ on a  manifold $M$ if it satisfies the following properties:
\begin{equation}
\label{act}  w(e,z)=z,\qquad w(a,w(b,z))=w(ab,z),
\end{equation}  
for any $z\in M$, $a,b\in G$.
When it does not lead to confusion an action will be denoted as multiplication:
 $w(a,z)=a\cdot z$
\end{df}
We adopt the definition of local group actions from \cite{O86}.
\begin{df} 
A {\em local group of transformations} on $M$ is given by a (local) Lie group $G$ and an open subset  ${\cal U}\subset G\times M$ , such that
$\{e\}\times M\subset {\cal U} \subset G\times M,$ which is the domain of definition of the group action, and a smooth map $w\colon {\cal U}\rightarrow M$ such that if 
$(g,z)\in {\cal U}$ the so is 
$(g^{-1},w(g,z))$ and two properties (\ref{act}) are satisfied whenever $w$ is  defined.
\end{df}
The following definition extends the notion connectness  to local 
transformations:
\begin{df} A group of transformations $G$ acting on $M$ is called {\em connected} if the following requirements hold:
\begin{description}
\item[(i)] $G$ is a connected Lie group and $M$ is a connected manifold;
\item[(ii)] the domain ${\cal U} \subset G\times M$ of the group action is connected;
\item[(iii)] for each $z\in M$ the local group $G_z=\{g\in G|(g,z)\in {\cal U}\}$ is connected.
\end{description}
\end{df}
\begin{df}
The {\it orbit} ${\cal O}_z$ through a point $z\in M$ is the image of the
 map $w_z:G\rightarrow M$ given by $w_z(g)=w(g,z)$. 
Each point of $M$ belongs to a unique orbit. If there is only one orbit on $M$, then 
${\cal O}_z=M$ for all $z$ and the action is called {\it transitive}.   
\end{df}
The differential $dw_z:TG|_e \rightarrow TM|_z$  maps  
the Lie algebra of $G$ to the tangent space at the point $z$. Let $X\in {\frak g}=TG|_e$ 
then  ${\hat X}(z)= dw_z(X)$ is a smooth vector field on $M$ called an {\em 
infinitesimal generator} of  the  $G$-action:
\[ \exp(tX)\cdot z =\exp(t\hat X,z),\]
where $\exp(t\hat X,z)$ is the  flow of the vector field $\hat X$.

\begin{df} The {\em isotropy} group of a subset $S\subset M$ is 
\[G_S=\{g\in G|gS=S\}.\]
The {\em global isotropy subgroup} is the subgroup
\[G_S^*=\bigcap_{x\in S}G_x=\{g\in G|gS=S\}\]
consisting of those group elements which fix all points in $S$.
\end{df}
\begin{df} A group $G$ acts \textit{effectively} if different elements have
different actions, or equivalently  $G_M^*=\{e\}$. The action of $G$ is locally effective if $G_M^*$
is a discrete subgroup of $G$.   A group  $G$ acts \textit{effectively on subsets }
if the global isotropy subgroup of each open subset $U\subset M$ is trivial.
\end{df}
If a group does not act
effectively then we can replace its action with the action of the quotient
group $G/G_M$. This action is well defined, effective and essentially the
same as the action of $G.$  Clearly if $G$ acts \textit{effectively on subsets }then $G$
 acts effectively. The converse statement is true in analytic category, however as example 2.3 
in \cite{FO99} shows, it is not valid in general in smooth category.

\begin{df} A transformation group \textit{acts freely} if the isotropy subgroup of each
point is trivial and  \textit{locally} \textit{freely} if the
isotropy group of each point is discrete.
\end{df}
If the group acts (locally) freely then the  dimension of each orbit ${\cal O}_z$ is equal to 
the dimension of the group. In this case the map $w_z$ defines 
a (local)  diffeomorphism between $G$ and the orbit  ${\cal O}_z$.

An action of a Lie group $G$ on a manifold $M$ induces the  standard linear representation
on the space of smooth functions $F\colon M\rightarrow \R$:
\begin{equation}\label{gF}
(gF)(gz)=F(z),
\end{equation} 
where $g\in G,\, z\in M.$ 

\begin{df} A function $F(z)$ is {\em invariant} if it is a fixed point of the standard 
representation above, that is 
\begin {equation}\label{gF=F}
F(gz)=F(z).
\end{equation}
 We say that  $F$ is a {\it local invariant} if it is defined on an open  subset of $M$,
 and/or the equality (\ref {gF=F}) holds  only for $g$ in a
neighborhood of identity in $G$. 
\end{df}

\begin{df} The {\em symmetry group} $G_F$ of a function $F(z)$ is the isotropy
 group of $F$ under representation (\ref{gF}):
$$G_F=\{g\in G|F(gz)=F(z)\}$$
\end{df}

\begin{df}The action of a group $G$ on $M$ is called \textit{semi-regular} if all its
orbits have the same dimension. The action is called \textit{regular }if, in
addition, each point $z\in M$ has arbitrarily small neighborhood whose
intersection with each orbit is a connected subset thereof. 
\end{df}

Let $G$ act semi-regularly on an $m$-dimensional manifold $M$ and let $s$ be the 
dimension of 
the orbits, then  the infinitesimal generators of the $G$-action form an integrable 
distribution of the dimension $s$. The orbits  of $G$ are the integral manifolds for 
this distribution.
By  Frobenius' theorem, coordinates $(x^1,\dots,x^s,y^1,\dots,  y^{m-s})$ on a chart  $U\subset M$ can 
be chosen so that each orbit is a 
level set of the last  $m-s$ coordinates:  $y^i=c_i,\, i=1,\dots, m-s.$ The functions 
$y^i$ form a complete set of functionally independent local invariants.  Thus the number of functionally independent local invariants for  a semi-regular action of a Lie group 
equals to  the difference between the dimension of the manifold and the dimension of the orbits.  Let  $\CS$ be a level set  of the first $s$ coordinates $x^i=c_i,i=1,\dots,s$,
 then the submanifold $\CS$ has codimension $s$ and is transversal to each  orbit in $U$. Thus $\CS$ intersects each orbit in a discrete set of points. If the action is regular
 and $U$ is sufficiently small, then $\CS$ intersects each orbit only once.

\begin{df}\label{cs} Suppose $G$ acts semi-regularly on an $m$-dimensional manifold $M$ with $s$-dimensional orbits. A (local) {\em cross-section} is an $(m-s)$-dimensional manifold 
$\CS\subset M$ such that $\CS$ intersects each orbit transversally. The cross-section 
is {\em regular} if $\CS$ intersects each orbit at most once.
\end{df}  

We conclude this preliminary section with Lie's infinitesimal criterion of invariance: 
\begin{thm}\label{Lie} Let $G$ be a connected group of transformations acting on a manifold $M$. A function $I\colon M\rightarrow \R$ is invariant under $G$ if and only if for every 
infinitesimal of generator $X$ of the $G$-action:
\[ X[I] \equiv 0.\]
\end{thm}

\begin{rem}\label{local} {\em Throughout the thesis we will, without saying it explicitly, consider  
the  group actions, cross-sections   and invariants to be local. 
Nevertheless in order to shorten the  statements and formulas  
we will keep global notation.  For instance  we will  write  
$w:G\times M\rightarrow M$ instead of $w\colon {\cal U}\rightarrow M$ for 
a local action, or we will call
$\CS$  a regular cross-section on $M$ meaning that it is a regular cross-section 
in some open subset $U\subset M$.}
\end{rem}

\section {Jet Spaces and Differential Invariants.}

In the case when the group $G$ acts transitively on $M$, there are clearly no 
non-constant invariants. 
Nevertheless if we consider the transitive action of the group of
Euclidean motions  on the plane we can find important geometric invariant:
the curvature $\kappa=\frac {u_{xx}}{\sqrt{1+u_x^2}^3}$ of an embedded
curve $u=u(x)$. We notice that $\kappa$  depends not only on a point on the curve 
but also on the derivatives of $u$ with respect to $x$. This  is an example 
of a {\it differential invariant}, which is an ordinary invariant function  on the  
prolonged space (or jet space). 
 Differential invariants were used by Lie in his work on the symmetry reduction of differential equations. 
The formal definition of jet spaces was first given by Ehresmann.  
We are going to give a brief  description of the geometric structure of jet spaces, 
for more details see \cite{And92}, \cite{O86}, \cite{O95}.
\begin{df}\label{JkMp}
Given a smooth manifold $M$ of dimension $m$ and an integer $p<m$, the $k$-{\it th 
order jet bundle} $J^k=J^k(M,p)$ is a fiber bundle over $M$, such that
a fiber over a point $z\in M$ consists of the set of
equivalence classes  of $p$-dimensional submanifolds of $M$ with  $k$-$th$ order 
contact at $z$. In particular $J^0=M$.
\end{df}
\begin{rem} In the case when $M$ has itself a fiber  bundle structure $M\rightarrow B$
 with a  $p$-dimensional base, then its $k$-th  jet bundle $J^kM$ can be defined  by  
 sections $s\colon B\rightarrow M$ under the equivalence relations of $k$-$th$ order 
contact at $z\in M$. Then the  jet bundle $J^k=J^k(M,p)$ from 
Definition~\ref{JkMp}  can be called the {\em extended jet bundle}, since it also includes  
jets of  $p$-dimensional 
submanifolds of $M$ that are not transversal to the fibers over $B$. The fiber bundle 
$J^kM\rightarrow M$ is an open dense subset of the extended bundle $J^k(M,p)$.
\end{rem}

There is a natural projection $\pi^k_l:J^k\rightarrow J^l$.
The inverse sequence  of topological spaces 
$(J^k,\pi^l_k)$  determine  an inverse limit space $J^\infty=J^\infty(M,p)$ 
together with 
projection maps  $\pi^\infty_k:J^\infty \rightarrow J^k$. The space  $J^\infty$ 
is called  infinite jet bundle over $M$.  
In the same manner the tangent bundle $TJ^\infty$ can be defined as the inverse 
limit of topological spaces $TJ^k$ under the projections 
$(\pi^k_l)_*=d\pi^k_l$. For $l<k$ we identify functions and differential forms on $J^l$ 
with functions and 
forms on $J^k$ under the pull-back $(\pi^k_l)^*$. The smooth functions and forms 
on $J^\infty$
are defined as the direct limit of the space of smooth functions and forms on  on $J^k$.

Let $U$ be a chart of $M$ with coordinates $(x^1,\dots x^p,u^1,\dots, u^q)$ so that
$p+q=\dim M$.  We say that a submanifold $S\subset M$ is {\em transverse}
 with respect to this coordinates if the restriction of the forms $dx^1,\dots,dx^p$ to 
$S$ is  a coframe on $S$. Any   transverse submanifold  can be locally described  as a graph 
$u^\alpha=f^\alpha(x),\, \alpha=1,\dots,q$. Let $U^0\subset U$ be a union of 
transverse submanifolds of $U$
and $U^k\subset J^k$ be a subset such that $\pi^k_0(U^k)=U^0$.
 It is not difficult to see that $U^k$ is an  open subset 
of  $J^k$, which can be parameterized  by the  set of independent 
variables $\{x^1,...,x^p\}$,  the set $\{u^1,...,u^q\}$ of dependent variables  and coordinates $u_J^\alpha$ which correspond to
the derivatives of the dependent variables with
 respect to the 
independent ones, where the  subscript $J=(j_1,\dots, j_p)$ is 
a multi-index, such that $|J|=j_1+\dots + j_p\leq k,\, j_i\geq 0$.

\begin{df}
Let $S$ be a $p$-dimensional submanifold of $M$, then its $k$-th prolongation $j^k(S)$ is 
a $p$-dimensional submanifold  of  $J^k(M,p)$, defined by the $k$-th jets of $S$.
In local coordinates on an open set $U^k$ it is the graph of equations  
$$u^\alpha=f^\alpha(x),\qquad
 u^\alpha_J=\frac{\partial^k f^\alpha}{\partial^{j_1} x^1\ldots\partial^{j_p} x^p},$$
where $ \alpha =1,\dots,q$ and  $J=(j_1,\dots, j_p)$ are all possible  
multi-indices, such that $|J|\leq k$.
\end{df}
Although the  $k$-th prolongation $j^k(S)$ of $S$  is  a $p$-dimensional 
submanifold of $J^k(M,p)$, not every $p$-dimensional  submanifold of $J^k(p,M)$ is the 
prolongation of a submanifold in $M$. The exterior  differential forms which are 
identically zero when restricted to the  prolongation $j^k(S)$ of 
any submanifold $S\subset M$ for all
$k$  form a differential ideal on $J^\infty$, 
called the {\it contact ideal}. 
In local coordinates a  basis for the contact ideal can be written as:
 \[
\theta _J^\alpha =du_J^\alpha -\sum_{i=1}^pu_{J,i}^\alpha dx^i,\qquad \alpha
=1,...,q,\quad |J|\geq 0.
\]
Thus the  cotangent bundle on $J^\infty(U,p) $ splits into two
sub-bundles: the horizontal sub-bundle which is  spanned
by the forms $dx^1,...,dx^p,$ and the vertical sub-bundle spanned by the
contact forms. 
We emphasize that the  contact (vertical) sub-bundle has  intrinsic definition, 
independent of the choice of  coordinates, whence the choice of basic horizontal 
one-forms depends on the choice of independent coordinates.
 
The differential $d$ on $J^\infty $ also
splits into horizontal and vertical components, 
$$d=d_H+d_V.$$ 
The horizontal
differential is defined by 
$$d_HF=\sum_{i=1}^p\left( D_iF\right) dx^i,$$ 
where 
\[
D_i=\frac \partial {\partial x^i}+\sum_{\alpha =1}^q\sum_Ju_{J,i}\frac
\partial {\partial u_J^\alpha }. 
\]
The operators $D_i$ span a subspace of  {\em total vector fields} in the tangent bundle 
 $TJ^\infty$, which 
can be defined intrinsically as the set of vector fields annihilated 
by any contact form.

\noindent The vertical differential is defined by
$$d_VF=\sum_{\alpha,J}\frac{\partial F}{\partial u_J^\alpha }%
d\theta _J^\alpha. $$
The vector fields $ \frac{\partial }{\partial u_J^\alpha }$, which are annihilated by 
any horizontal form,  span the subspace of {\em vertical vector fields} in $TJ^\infty$. 
Each one-form on $J^\infty$ splits into horizontal and vertical component, inducing a 
bi-grading of the exterior differential forms  on $J^\infty$. The horizontal differential
$d_H$ increases the horizontal degree of a form and the vertical differential $d_V$ increases the vertical degree.
This splitting gives rise to a bicomplex of differential forms called  the variational 
bicomplex \cite{Tul82}, \cite{Tsu82}, \cite{Vin84}, and \cite{And92}, 
an important tool in the study of  geometry of differential equations and variational problems.

\begin{df}
 Let $X$ be  a vector field   on $M$, then  there is a unique  vector field $pr X$
 on $J^{\infty}$, called prolongation of $X$ such that
\begin{description}
\item[(i)] $X$ and $prX$ agree on functions on $M$,
\item [(ii)] $prX$ preserves the contact ideal: the Lie derivative of a contact form with respect to $prX$ is also a contact form.
\end {description}
The vector field $pr^kX=(\pi^{\infty}_k)_*prX$ is called  $k$-th  prolongation of $X$. 
\end{df}
In terms of local coordinates let
$$X=\sum_{i=1}^p\xi ^i(x,u)\frac \partial {\partial
x^i}+\sum_{\alpha =1}^q\varphi ^\alpha (x,u)\frac \partial {\partial
u^\alpha }$$
 be a vector field on $M$, then
its $k$-th order \textit{prolongation }is: 
\begin{equation}\label{pr(X)}
pr^k(X)=\sum_{i=1}^p\xi ^i(x,u)\frac \partial {\partial x^i}+\sum_{\alpha
=1}^q\sum_{|J|\leq k}\varphi _J^\alpha (x,u^{|J|})\frac \partial {\partial
u_J^\alpha }, 
\end{equation}
where
\[
\varphi _J^\alpha (x,u)=D_JQ^\alpha +\sum_{i=1}^p\xi ^i(x,u)u_{J,i}^\alpha. 
\]
and  $Q^\alpha $ denotes the characteristics of the vector field $X:$%
\[
Q^\alpha (x,u^{(1)})=\varphi ^\alpha (x,u)-\sum_{i=1}^p\xi ^i(x,u)u_i^\alpha
. 
\]
Any prolonged vector field can be decomposed into total and vertical components:
\begin{eqnarray*}
pr(X)=\sum_{i=1}^p\xi ^i(x,u)\frac \partial {\partial x^i}+\sum_{\alpha,J}
\left(D_JQ^\alpha+\sum_{i=1}^p\xi^i(x,u)u_{J,i}^\alpha\right)\frac \partial {\partial
u_J^\alpha }\\
=\sum_{i=1}^p\xi ^i(x,u)D_i+\sum_{\alpha,J}
D_JQ^\alpha\frac \partial {\partial
u_J^\alpha }.
\end{eqnarray*}

After we have defined the jet space we would like extend the action of $G$ on it so that
the extended action maps a prolongation of a submanifold of $M$  to the 
prolongation of its image.

\begin{df}
The {\it $k$-th order prolongation} of a smooth transformation $g\in G$ on $M$
is defined by the property 
\[g^{(k)}\cdot j^k(S)=j^k(g\cdot S)
\]
for any submanifold $S\in M$.
\end {df}
It follows from the definition that  the following diagram commutes 
$$\square<1`1`1`1;500`500>[{J^k}`{J^k}`{J^l}`{J^l};g^{(k)}`\pi^k_l`\pi^k_l`g^{(l)}]
$$

When it is clear from the context that the prolonged transformation is considered we 
will  omit the superscript $(k)$ over $g$. 

\begin{df}
A $k$-th order \textit{differential invariant} on $M$ under the action of $G$ is a
 function on $J^k(M,p)$ which is invariant under the $k$-th prolongation of the group
action.
\end{df}
\begin{rem} \label{lkinv} We consider  an  $l$-th order differential invariant for $l<k$ as a $k$-th order differential invariant
under the pull-back $(\pi^k_l)^*$. Thus a $k$-th order differential invariant might only
depend on derivatives of order strictly less than $k$.
\end{rem} 
If $X$ is an infinitesimal generator of $G$-action on $M$, then its $k$-th prolongation
(\ref{pr(X)})   generates the corresponding prolonged transformation on 
$J^k(M,p)$. 
From  Lie's criterion \ref{Lie} it follows that function $I\colon J^k\rightarrow \R$
 is a differential invariant  if and only if for every infinitesimal of generator $X$ 
of the $G$-action 
\[pr^k(X)[I] \equiv 0.\]
Theoretically this criterion can be used to find all differential
invariants. In practice however it is difficult to use since it requires
integration of a system of first order partial differential equation.
The advantage of the  method of  moving frames is that  it  requires only 
differentiation not integration.

Let ${\cal O}^k$ and ${\cal O}^l$ be the orbits of prolonged action on 
$J^k$ and $J^l$ respectively for $k>l$, then  $(\pi^k_l)({\cal O}^k)={\cal O}^l$ 
and hence the dimension of the orbits  
can only become larger when we prolong the action to the higher jet spaces.
 On the other hand the dimension of the orbits is bounded by the dimension of the group $G$. Hence  there is an order of prolongation $n$ at which the maximum possible dimension 
is attained on an open subset of $\CV^n\subset J^n$. If $\zk,\, k>n$ is a point in $J^k$ 
such that $\pi^k_n(\zk)\in {\CV^n}$ then the orbit of the prolonged action through $\zk$ 
also has the maximal dimension. We call such points {\em regular jets} and denote their 
union in $J^\infty$ as ${\cal V}$. 

\begin{df}
The minimal order at which the orbits reach  maximal dimension is called {\em the order 
of stabilization}. The subsets $\CV^k\subset J^k,\quad k=n,\dots,\infty $ 
which consists of the 
points $z^{(k)}$ such that the orbit through $\pi^k_n z^{(k)}$ has 
maximal dimension are called {\it regular}.
\label{defst}
\end{df}

The following result (Ovsiannikov  \cite{Ovs82}, Olver \cite{Ol00}) is crucial 
for moving frame construction. 
\begin{thm}
If the action of $G$ on $M$ is locally effective on subsets then the prolonged 
action is locally free on $\CV^k$, for $k\geq n$, where $n$ is the order of stabilization.
\label{ovsth}
\end{thm}
Since the dimension of the space grows with prolongation and the dimension of orbits stabilizes at order $n$, then  nontrivial local differential 
 invariants are guaranteed  to appear at least at the order $n+1$ by Frobenius' theorem.
Although there exists only a finite number of functionally independent
differential invariants at each order of prolongation, their total number is
infinite since one can prolong up to the infinite order. Fortunately one can 
obtain all invariants from a finite   {\em generating set of  invariants} 
by applying {\em invariant differential}
 operators. The latter can be constructed  as  dual vector fields  
to horizontal contact invariant forms, which are defined as follows: 
\begin{df}
A differential  one-form $\omega $
on $J^n$ is called {\em contact invariant} if for every $g\in G$ we have 
$\left(g^{(n)}\right)^*\omega =\omega +\theta_g $ for some contact form $\theta_g $.  
A set of $p$
linearly independent horizontal contact invariant   forms $\{\omega _1,...,\omega _p\}$
is called a {\em horizontal contact invariant coframe}.
\end{df}
 The freeness of the group action on $\CV^n \subset J^n$ guarantees the existence  of a
 contact invariant coframe $\omega^1,\dots,\omega^p$ on $\CV^n$ \cite{FO97}, \cite{Ovs82}. 
  The horizontal differential of a function $F$ can be written in terms of this coframe 
as: 
\[
d_HF=\sum_{i=1}^p\left( \mathcal{D}_iF\right) \omega ^i, 
\]
where the total vector fields  $\mathcal{D}_i$ have an important property: they commute with the prolonged action of $G$ and thus map
differential invariants to higher order differential invariants. Operators which possess
 this property  are called \textit{invariant differential operators.} 
The following theorem \cite{Ovs82},  \cite{FO99} asserts that  one can produce 
all differential invariants by applying a finite set of invariant differential 
operators to a finite set of generating  invariants.

\begin{thm} Suppose that $G$ is a transformation group and let $n$
be its order of stabilization. Then, in a neighborhood of any regular jet $%
z^{\left( n\right) }\in \mathcal{V}^n$, there exists a contact invariant
horizontal coframe $\{\omega _1,...,\omega _p\}$, and corresponding {\em invariant
differential operators} $\mathcal{D}_1,...,\mathcal{D}_p$. Moreover, there
exists a {\em generating system of differential invariants }$I_1,...,I_l$%
, of order at most $n+1$, such that, locally, every differential invariant
can written as a function of $I_1,...,I_l$ and their invariant derivatives:
\begin{equation}\label{recf}
I=H(\dots,\CD_{J}I_j,\dots ),
\end{equation}
where $\CD_{J}$ is a certain composition of invariant differential operators. 
\label{indiff}
\end{thm}

\begin{ex}
 Let us consider  the special Euclidean group $SE(2,\R)=SO(2,\R)\ltimes \R^2$ acting on
 curves in $\R^2$. The group acts freely and transitively on the first jet space $J^1$ and 
thus the lowest order  invariant, the Euclidean  curvature 
$$\kappa=\frac{u_{xx}}{\sqrt{1+u_x^2}^3}$$
appears on the  second jet space.
The contact invariant coframe consists of a single  form, infinitesimal  arc-length:
$$ds={\sqrt{1+u_x^2}}\,dx.$$ 
Higher order differential invariants can be obtained by taking the derivatives of 
$\kappa$ with respect to the arc-length, or in other words by applying invariant 
differential operator
$\frac 1 {\sqrt{1+u_x^2}}\, D_x$. From dimensional consideration it is clear that 
only one new functionally independent invariant appears at each order of prolongation, 
and thus
any differential invariant is a  function of 
$\kappa,\kappa_s=\frac {d\kappa} {ds},\kappa_{ss}=\frac {d \kappa_s} {ds}$, etc. 
\end{ex}
The generalized method of moving frame \cite{FO99} described  
in Section~\ref{smf} provides a consisting 
algorithm of constructing differential invariants, a  contact invariant coframe,
 invariant differential operators and recursion formulas (\ref{recf}). 
In the next section we explain the  role played by differential invariants in the solution of the equivalence and symmetry problems for submanifolds.

\section{Equivalence and Symmetry of Submanifolds. Signature Manifolds.}

\begin{df} Let $G$ be a group acting on a manifold $M.$ Two submanifolds $S$ and
$\bar{S}$ of dimension $p$
in $M$ are said to be locally {\it equivalent} if there exist a  point $z\in M$ and 
 an element $g\in G$ 
such that $\bar{S}=g\cdot S$  in a  neighborhood of the point  
$g\cdot z \in \bar{S}$. 
An element $g\in G$ is a local {\it symmetry} of $S$ if $S=g\cdot S$ at 
least in a neighborhood of  the  point $g\cdot x \in S$. 
\end{df}

Both the equivalence and the symmetry problems 
can be solved by the following construction. 
Let $n$ be the order of stabilization as in definition~\ref{defst}, 
then the group acts locally 
freely  on  ${\CV^n}\in J^n(M,p)$ and the differential invariants must appear at 
least on  $J^{n+1}$. Let   $\{I_1,...I_{N_k}\}$ be a {\em complete set} of $k$-th order differential invariants  for $k>n$, where  by complete set we  mean that any other 
invariant of order $k$ is a function of   $\{I_1,...I_{N_k}\}$.  
Let $\{\tI_1,...\tI_{N_k}\}$ be the restriction of $I$'s to $j^k(S)$.
The functions $\{\tI_1,...\tI_{N_k}\}$ 
define a map $\phi_k:S \rightarrow \R^{N_k}$.
Let $t_k=rank(\phi_k)$, then $t_k\leq p$ since  $\dim S=p$. 
It is not difficult to see that the rank $t_k$ does not depend on a choice of a 
complete system of invariants.
\begin{df}\label{regs} A submanifold $S\subset M$ is called {\em regular} if $j^n(S)\subset \CV^n$
and the rank $t_k,\, k\geq l$ does not vary on  $S$.
\end{df}
\begin{df}\label{sign} Let $S$ be a regular submanifold of $M$  then its 
{\it $k$-th order signature manifold} ${\cal C}^{k}(S)$ is
 an immersed submanifold  $Im(\phi_k)\subset \R^{N_k}$, where $\phi_k $ is defined as above.
\end{df}

Taking into account Remark~\ref{lkinv} we observe that 
$t_{n+1}\leq t_{n+2}\leq \dots \leq  t_{n+i}...\leq p$ and hence this sequence 
stabilizes. In fact  once $t_s=t_{s+1}$ for some $s$, all the subsequent ranks are equal.
\begin{prop}
If $t_s=t_{s+1}$ in  the sequence of ranks above then:
$t_s=t_{s+1}=t_{s+2}=\dots$
\label{rank}
\end{prop}
{\it Proof.} Given that $t_s=t_{s+1}$, we need to prove that $t_{s+1}=t_{s+2}$. 
Let $I$ be any invariant of order $s+2$, then from   theorem 
 \ref{indiff} it follows that there exist invariants $I_1,...,I_N$ of order $s+1$ such that $I=H(\CD_{i_1}I_1,\dots \CD_{i_N}I_N)$, where $H$ is some function. 
From the definition of total vector fields it follows that:
\begin{equation}\label{1}
\tI=H\left(\CD_{i_1}\tI_1,\dots, \CD_{i_N}\tI_N\right), 
\end{equation}
where  $\tI,\tI_1,...,\tI_N$ are restrictions of $I$'s to the jet $j^{s+2}(S)$.
On the other hand since  $t_{s+1}=t_s$  there exist $s$ order invariants $F_1,\dots F_t$ such that  each $\tI_i$ can be written  
as a function of the restrictions  $\tF_1,\dots \tF_t$  of $F$'s to   $j^{s}(S)$:
\begin{equation}\label{2}
\tI_i=H_i(\tF_1,\dots \tF_t).
\end{equation}
By substitution of (\ref{2}) in (\ref{1}) we obtain that 
\begin{equation}
\tI=H\left(\CD_{i_1}(H_1(\tF_1,\dots,\tF_t)),\dots, 
\CD_{i_N}(H_N(\tF_1,\dots,\tF_t))\right). 
\end{equation}
We note that  $\CD_{i_j}(H_j(\tF_1,\dots,\tF_t)),\, j=1,\dots,N$ are differential 
 invariants of order $s+1$ and thus
any $s+2$ 
order invariant restricted to $j^{s+2}(S)$
can be written as a function of the invariants of order $s+1$ restricted to $j^{s+1}(S)$.
\qed 

\begin{df}\label{dRankOrder} The minimal order $s$  such that $t_s=t_{s+1}$ is called {\it differential 
invariant order} of $S$ and the corresponding rank $t=t_s$ is called 
{\it differential invariant rank} of $S$.
\end{df}

\begin{rem}\label{rems} From Proposition~\ref{rank} it follows that
$t_n<t_{n+1}<\ldots <t_s=t_{s+1}\leq p$
and thus $s \leq n+p$, where $n$ is the order of stabilization. 
\end{rem}
From Proposition~\ref{rank} and recurrence formulas (\ref{recf}) it follows that the 
signature  manifold ${\cal C}^{s+1}(S)$   encodes all  functional relations 
among invariants restricted to $j^\infty(S)$.
Remark~\ref{signgeom} in the next section gives  a geometric description of the 
signature manifold. As the following two theorems \cite{O95}, \cite{FO99} show, the signature manifold of order $s+1$  plays a 
crucial role in  the
 solution of the equivalence and symmetry problems for submanifolds . 
\begin{thm} \label{equiv}Let $S,\bar{S}\in M$ be two regular $p$-dimensional submanifolds. 
Then $S$ and $\bar{S}$ are (locally) equivalent,  
$\bar{S}=gS$ if 
and only if they have the same differential order $s$ and their signature manifolds
(locally) coincide: ${\cal C}^{s+1}(S)={\cal C}^{s+1}(\bar{S})$
\end{thm}
\begin{thm}\label{symm} Let $S\subset M$ be a regular $p$-dimensional submanifold of differential 
invariant rank\,  $t$ with respect to the transformation group $G$. 
Then its isotropy group $G_S$ is a $(p-t)$-dimensional subgroup of $G$ acting locally freely on $S$. 
\end{thm}
\begin{cor} A submanifold $S\subset M$ 
has  a discrete symmetry group if  and only if its signature
 manifold has  maximal dimension $p$.
\end{cor}
\begin{rem}\label{rfinite}
The condition $rank(\phi)<p$ is closed and hence  a generic submanifold
 has a discrete group of symmetries.
\end {rem}

If the symmetry group is finite then its cardinality  can be found from the following theorem \cite{O99}.

\begin{thm} \label{sfinite} Let $\dim\, {\cal C}^{s+1}(S)=p$ and  $c\in {\cal C}^{s+1}(S)$ be a generic point. Then 
the cardinality of the symmetry group is equal to the cardinality of the preimage of $c$
 under $\phi$:\,  $|G_S|=|\phi^{-1} (c)|$ .
\end{thm}
\begin{ex}
 Let us return to the problem of equivalence and symmetry  of  curves in $\R^2$ under the action of  the special Euclidean group $SE(2,\R)=SO(2,\R)\ltimes \R^2$. In this case the order of stabilization $n=1$, $\CV=J^1$ and so  the differential order $s$  of any curve is not 
greater than two (see Remark~\ref{rems}). Thus the symmetry and equivalence problem for any submanifold can be solved by considering the signature manifold parameterized by the  curvature $\kappa=\frac{u_{xx}}{\sqrt{1+u_x^2}^3}$ and its derivative $\kappa_s=\frac {(1+u_x^2)u_{xxx}-3u_xu_{xx}^2}{(1+u_x^2)^3}$ with respect to the arc-length.
If $\kappa$,  restricted to the jet of a curve, is  constant then  
the signature manifold degenerates to  a point $(\kappa,0)$. The curvature is constant
 when the curve is either a circle or a line. The symmetry group of a line consists  
of translations parallel to this line and the symmetry group of a circle consists of 
rotations around the center of the circle.  Both groups are one-dimensional in
 agreement with  Theorem~\ref{symm}. The 
curvature of a generic curve is non-constant, and thus the signature manifold is one-dimensional. Let us construct 
the signature manifolds for two graphs $u=\sin(x)$ and $u=\cos(x)$, which are clearly equivalent under translation by $\frac{\pi}2$ in $x$-direction.
\begin{eqnarray*}
u=\cos(x)  & |&u=\sin(x)\\
           & |&		\\	  
\kappa=-\frac{\cos(x)}{(1+\sin^2(x))^{(3/2)}}&|&\kappa=\frac{\sin(x)}{(1+\cos^2(x))^{(3/2)}}\\
           &|&        \\
\kappa_s=\frac{\sin(x)(1+\cos^2(x))}{(1+\sin^2(x))^3}&|&
\kappa_s=-\frac{\cos (x)(1+\sin^2(x))}{(1+\cos^2(x))^3}\\
 & &
\end{eqnarray*}
\begin{center}
{\scalebox{0.4}{\includegraphics{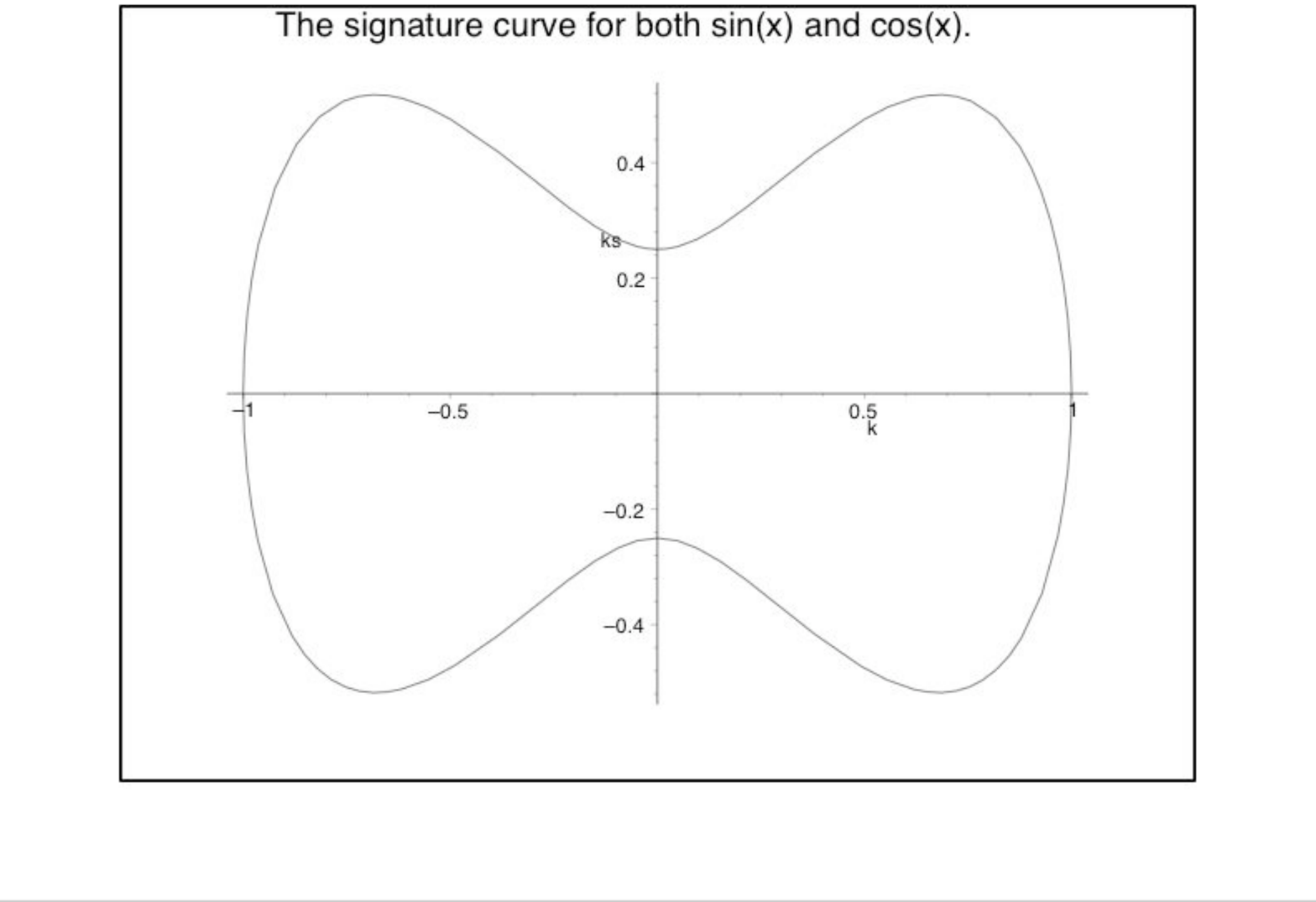}}}
\end{center}


As $x$ varies the signature curve will be traced over and over 
again with the period $2\pi$. This reflects the fact that the graphs of sine and cosine possess infinite discrete group of symmetries: translation by $2\pi k\,, k\in\Z$.
\label{ex1}
\end{ex}

\section {The Method of Moving Frames.}
\label{smf} 
In this section we describe the generalization of the Cartan's method of moving frame by Fels and Olver \cite{FO99}
which  provides an algorithm to construct a complete set of differential invariants of any order as well as a contact invariant coframe $\omega^i,\,i=1,\dots,p$, corresponding invariant differential operators ${\cal D}_i$ and recursion formulas (\ref{recf}).  
\begin{df}
Let a group $G$ act on a smooth manifold $N$. A {\it (right) moving frame} is defined 
as an equivariant map $\rho: N\rightarrow G$, where $G$ acts on itself  by  right multiplication. In other words the following diagram commutes:
$$\square<1`-1`-1`1;500`500>[{G}`{G}`{N}`{N};R_{g^{-1}}`\rho`\rho`g]
$$    
                                                                                                                                                                             
\end{df}
\begin{thm}
A local moving frame exists if and only if $G$ acts regularly and (locally) freely on $N$. 
\label{mf}
\end{thm}
{\it Proof.}
We sketch the proof of the sufficiency of these 
conditions and refer the reader to Theorem~4.4 in
\cite{FO99} for  the proof of their  necessity. 
Construction of a moving frame is equivalent to choosing a regular cross-section 
(see Definition~\ref{cs}). If an $r$-dimensional Lie group $G$  acts regularly and 
locally freely on $M$ then in a neighborhood $U$ of every point $z\in N$ there exists 
a cross-section $\CS$ that intersects each orbit at a unique point. We define a map
$ \rho:U\rightarrow G$   by the condition 
$$\rho(z)\cdot z \in \CS.$$
If $G$ acts freely on $U$ then the map  $\rho(z)\in G$ is well defined, if $G$ acts locally freely then the map on a neighborhood of the identity $V_e\subset G$ is well defined.
To show that $\rho$ is equivariant we observe that   
  $\rho(gz)\cdot gz=\rho(z)\cdot z$ for all $g\in G$ 
and so since the action of $V_e\subset G$ is free then $\rho(gz)=\rho(z)g^{-1}$.
In the future we will adopt  global notation for local objects and maps 
(see Remark~\ref{local}).
\qed
\begin{df} A moving frame on a submanifold $S\subset N$ is  the restriction of $\rho$
to $S$.
\end{df}
We notice that non-constant coordinates of $\rho(z)\cdot z$ provide a complete set of 
 functionally independent local invariants on $N$. In fact if $f(z)$ is any function on 
$N$ then its invariantization is defined as $\inv f(z)=f(\rho(z)\cdot z)$. In other words the function  $\inv f(z)$ is obtained by spreading the values of $f$ on $\CS$ along the orbits of 
$G$ and thus it is  invariant. 

This process can be described as normalization of  invariants on $\CB=G\times N$, which are called {\it lifted invariants}. Let the map $w:\CB\rightarrow N$ be defined by the 
group 
action, that is $w(g,z)=g\cdot z$, the map $\sigma:N\rightarrow \CB$ be defined by 
$\sigma(z)=(\rho(z),z)$ and the map $\iota: N\rightarrow N$ be defined by 
$\iota(z)=\rho(z)\cdot z$, then the following diagram commutes:
$$\Atriangle<-1`1`1;600>[{\CB=G\times N}`N`N;\sigma`w`\iota]
$$

We notice that the maps $w$ and $\iota$ are $G$-invariant and the map $\sigma$ is $G$-equivariant.  Let $z^i$ be the $i$-th coordinate function on $U\subset N$, then  
$(gz)^i=w^*z^i,$  are invariant functions on $\CB$,  under the action 
$g(h,z)=(hg^{-1},g\cdot z)$. Such functions are called {\em lifted invariants}.
Their normalization $I^i=\sigma^*w^*z^i$ include a complete set of functionally independent invariants
on $N$. The invariantization of an arbitrary   function on $N$ is defined by its 
pull-back under  the map  $\iota=\omega\circ\sigma: N\rightarrow N$. Indeed if $f$ is 
any function on $N$ then  $\inv f=\iota^*(f)=\sigma^*w^*f$ is an invariant function. 

It is clear that if $\eta$ is any form on $N$ then $\Eta=w^*\eta$ is an invariant form 
on $\CB$ and  $\iota^* \eta=\sigma^*(\Eta)$ is an invariant form on $N$. The pull back of any coframe on $N$ under $\iota$ produces a set of invariant one-forms of codimension 
$r$. We can complete this set to an invariant  coframe on $N$ by pulling back a left invariant coframe on $G$ under $\rho$. 
 However,  we 
define invariantization of differential forms in a different manner, for the reasons 
which become apparent later.  
We notice that cotangent bundle on  $\CB$ splits into two subspaces: the cotangent bundle to
$G$ spanned by the Maurer-Cartan forms and the cotangent bundle to the manifold $N$. 
The differential $d$ on $\CB$  splits accordingly $d=d_G+d_N$, leading to a well defined 
bicomplex of differential forms on $\CB$.
Since the action of $G$ on $\CB$ is a product of the actions on $G$ itself  and the action on $N$ then the action 
of $G$ on $\CB$  preserves this splitting. In particular, if $\Eta$ is a $G$-invariant form on 
$\CB$ then so are  $\pi_N(\Eta)$ and $\pi_G(\Eta)$, where $\pi_N$ and $\pi_G$ are the projection of the forms on $\CB$ on purely manifold component and on purely group component respectively.
If $\eta$ is any form on $N$ then 
$w^*(\eta)$ is an invariant form on $\CB$, and so is  $\pi_N\, w^*(\eta)$. 
Since $\sigma$ is a $G$-equivariant map from $N$ to $\CB$, then  
$\sigma^*\,\pi_N\,w^*(\eta)$ is an invariant form on $N$.
\begin{df}\label{inv} 
Let $\eta$ be a differential form on $N$, then its {\em invariantization} is defined 
by the formula:
\begin{equation}
\inv(\eta)=\sigma^*\, \pi_N\, w^*(\eta).
\end{equation}
\end{df}
In particular, if the forms $dz^i$ form a coframe on $N$ then the forms 
$\inv(dz^i)=\sigma^*\, \pi_N\, w^*\, dz^i=\sigma^*d_N(w^* z^i)$ form a 
$G$-invariant coframe on $N$.
Applied to zero-forms (functions) this definition coincide with the one which was given 
above.
It is easy to check that $\inv$ is a projection operator, that is it maps any invariant 
form to itself.

One can apply the process of invariantization to construct differential invariants, 
contact invariant forms and invariant differential operators on $J^\infty(M,p)$.
Due to theorem \ref{ovsth} there exists high enough order of prolongation $n$ such that $G$ acts locally freely on $\CV^n\subset J^n(M,p)$. Thus  there exists a local cross-section 
$\CS^n$ and a corresponding local moving frame
$\rho: {\cal V}^n \rightarrow G$. We can extend these cross-section and moving frame 
to any higher order regular set $\CV^k$  by defining 
$\CS^k=\{\zk|\pi^k_n\zk\in \CS^n\}$ and 
$\rho(\zk)=\rho\left(\pi^k_n(\zk)\right)$ for $k=n,\dots,\infty$.
 Let the map $w^k:\CB^k\rightarrow J^k$ be defined by the prolonged group action, 
that is 
$w(g,\zk)=g\cdot \zk$, the map $\sigma^k:J^k\rightarrow \CB^k$ be defined by 
$\sigma(\zk)=(\rho(\zk),\zk)$ and the map $\iota^k: J^k\rightarrow J^k$ be defined by 
$\iota(\zk)=\rho(\zk)\cdot \zk$, then  the following diagram commutes for all 
$k=n,\dots,\infty$:
$$\Atriangle<-1`1`1;600>[{\CB^k=G\times J^k}`{J^k}`{J^k};
\sigma^k`w^k`{\iota^k}]
$$
Later we will omit the superscript $k$ for the maps between jet spaces.

The process of invariantization can be defined as in Definition~\ref{inv}, 
where the role of $d_N$ is played by the jet differential $d_J$. In particular 
if 
$$x^i,\,\, i=1,\dots,p,\qquad  u^\alpha_J,\,\, \alpha=1,\dots,q$$
 are  coordinates on 
$J^\infty$, then the functions 
$$y^i=w^*(x^i),\,\, i=1,\dots,p,\qquad  v^\alpha_J=w^*(u^\alpha_J),\,\, \alpha=1,\dots,q$$
are lifted invariants on $\CB^\infty$ and the functions
$$\CX^i=\inv (x^i), \, i=1,\dots,p, \qquad I^\alpha_J= \inv (u^\alpha_J),\,\alpha=1,\dots,q$$
 form  a complete set of  differential invariants on $J^\infty$.
\begin{rem}   The set of invariants $\{ \CX^i, \dots, \CX^p, I^\alpha_J, \, \alpha=1,\dots,q$
 \} is complete in a sense that every other differential invariant can be expressed as a function of these invariants. However exactly $r$ of these
 invariants are functionally dependent on the others. If the cross-section $\CS^n$ 
is chosen as a level set of $r$ coordinates, then the invariantization  of these $r$ 
coordinates produce  constant functions, called {\em phantom} invariants.
\end{rem}
\begin{rem}\label{signgeom} The $k$-th order signature manifold of a submanifold $S$
(see Definition~\ref{sign}) can be described as the image of its 
$k$-th prolongation $j^k(S)$ under the projection $\iota^k$. Indeed,
for $k=n,\dots,\infty$ the map $\iota^k$ projects $J^k$ onto the subset $\CS^k$ of
 codimension $r=\dim G$, parameterized by a  complete set of  
functionally independent invariants.
Thus if $S\subset M$ is any submanifold then its $k$-th order signature manifold is the projection of the $k$-th jet of $S$ under  $\iota^k$:\, ${\cal C}^k(S)=\iota^k(j^k(S))$ onto the cross-section chosen to define the corresponding moving frame. The symmetry and equivalence theorems~\ref{equiv},~\ref{symm},~\ref{rfinite} have a nice  geometrical 
interpretation  in terms of this projection. For instance  the dimension of the signature manifold 
decreases when the jet of submanifold $j^{s+1}(S)$ is not transversal to the prolonged orbits,
and hence there are infinitesimal generators of the group action which are tangent to
$j^{s+1}(S)$. These infinitesimal generators give rise to the symmetry group of $S$.
\end{rem}
Let us return now to the process of invariantization. The invariantization of the basis form 
$dx^1,...,dx^p,\, \theta^\alpha_J $:
\begin{eqnarray*}\label{invfr}
\varpi^i&=&\inv(dx)=\sigma^*\, d_J\, w^*\,(x^i),\,\, i=1,\dots,p,\\
\vartheta^\alpha_J&=&\inv(\theta^\alpha_J)=
\sigma^*\,\pi_J\, w^*\, (\theta^\alpha_J),\,\, \alpha=1,\dots,q
\end{eqnarray*} 
produces an invariant coframe on $J^\infty$. We recall 
that  jet differential $d_J$ splits into vertical and horizontal components,
thus the differential on $\CB^\infty$ splits into three components $d=d_G+d_H+d_V.$
Definition~\ref{inv} was motivated by 
the fact that such invariantization preserves the contact ideal: invariantization
of a 
contact form is a contact form. The contact ideal is defined intrinsically, while the 
choice of horizontal forms depends on the choice of coordinates and is not preserved by 
invariantization, that is an invariantized horizontal form might gain a vertical 
component.  
By projecting the invariantization of a horizontal form to its purely horizontal 
part we obtain a contact invariant form.  
Let $\eta$ be any form on $J^\infty$, then $\Eta=w^*\eta$ is an invariant form on 
$\CB^\infty$. Let $\pi_H$   denote the projection of a form on its purely horizontal 
 component. Since the action of $G$ preserves the contact 
ideal 
then $ \pi_H(w^*\eta)$   is a  horizontal contact 
invariant forms on $\CB^\infty$, and  $\sigma^*\pi_Hw^*(\eta)$ is a horizontal 
contact invariant form on $J^\infty$. In particular  forms
$$dy^i=\pi_H\, w^*(dx^i)=d_H\, w^*(x^i),\,\, i=1,\dots,p$$
form a horizontal contact invariant coframe on $\CB^{\infty}$.
The dual vector fields $\lD_i$ produce 
a complete set of {\em lifted invariant differential operators}, such that 
 $v^\alpha_{J,i}=\lD_iv^\alpha_ J$. 
The forms 
\begin{equation}\label{cinvc}
\omega^i=\sigma^*d_H\,w^*\, (x^i)=\sigma^*(d_H\, y^i),\,\, i=1,\dots,p
\end{equation}
form a horizontal contact invariant coframe on $J^{\infty}$. 
In general we can call $\sigma^*\pi_Hw^*\eta$ the contact invariantization of a form 
$\eta$ on $\Ji$. 
The vector fields $\CD_i$ dual to the forms $\omega^i$ provide 
a complete set of invariant differential operators.
We notice that in contrast with invariant 
differential operators  $\CD_i$ on $J^\infty$,
the lifted operators $\lD_i$ commute.
The formulas which relate normalized invariants with invariants obtained by 
invariant differentiation are given in (\cite{FO99}, section~13).

\begin{ex}
\label{ex2}
Let us return to  the special Euclidean action on the plane as described in \ref{ex1}:
\begin{eqnarray*}
x &\mapsto& y=\cos(\alpha)x-\sin(\alpha)u+a,\\
u &\mapsto&  v=\sin(\alpha)x+\cos(\alpha)u+b. 
\end{eqnarray*}
The lifted contact invariant form $d_Hy=\tau dx$, where 
$\tau=\cos(\alpha)-\sin(\alpha)u_x$. 
The lifted invariant differential operator $\lD_x= \frac 1 \tau D_x$ and 
thus the lifted differential invariants are given by formulas:
\begin{eqnarray*}
v_1&=&\frac {\sin(\alpha)+\cos(\alpha)u_x}\tau,\\
v_2&=&\frac {u_{xx}}{\tau^3},\\
v_3&=&\frac {\tau u_{xxx}+3\sin(\alpha)u_{xx}^2}
{\tau^5},\\
v_4&=&\frac {\tau^2u_{xxxx}
+10\sin(\alpha)\tau u_{xx}u_{xxx}+15\sin^2(\alpha)u_{xx}^3}
{\tau^7}.
\end{eqnarray*}
We note that the lifted invariants define the prolongation of the group action: 
$u_x\mapsto v_1,\, u_{xx}\mapsto v_2$, etc.
 The moving frame can be defined on $J^1(\R^2,1)$ by choosing a cross-section
$\{x=0,\, u=0,\, u_x=0\}$, and so an equivariant   map  $J^1(\R^2,1)\rightarrow SE(2)$ can be
 found by 
solving  the equations:
\begin{eqnarray*}
y&=&\cos(\alpha)x-\sin(\alpha)u+a=0,\\
v&=&\sin(\alpha)x+\cos(\alpha)u+b=0,\\ 
v_1&=&\frac {\sin(\alpha)+\cos(\alpha)u_x}{\cos(\alpha)-\sin(\alpha)u_x}=0.
\end{eqnarray*}
Thus we obtain the  moving frame:
\begin{equation}
\alpha=-\arctan(u_x),\qquad a=-\frac{u_xu+x}{\sqrt{1+u_x^2}},\qquad
b=\frac{u_xx-u}{\sqrt{1+u_x^2}}.
\end{equation}

The corresponding element of the special Euclidean group can be written in a  matrix form: 
$$\rho_r=\left(
\begin{array}{ccc}
\frac 1{\sqrt{1+u_x^2}}&\frac {u_x}{\sqrt{1+u_x^2}} &-\frac {uu_x+x}{\sqrt{1+u_x^2}}\\
-\frac{u_x} {\sqrt{1+u_x^2}}&\frac {1}{\sqrt{1+u_x^2}} &\frac {xu_x-u}{\sqrt{1+u_x^2}}\\
0&0&1
\end{array}
\right)
$$
The differential invariants are obtained by normalization, that is substitution of the moving frame into the lifted invariants $v_k$:
\begin{eqnarray*}
I_2&=&\frac {u_{xx}}{(1+u_x^2) ^{3/2}}\\
I_3&=&\frac {(1+u_x^2)u_{xxx}-3u_xu_{xx}^2}{(1+u_x^2)^3}\\
I_4&=&\frac {(1+u_x^2)^2u_{xxxx}-10u_xu_{xx}u_{xxx}(1+u_x^2)+15u_x^2u_{xx}^3}
{(1+u_x^2) ^{9/2}}
\end{eqnarray*}
The contact invariant differential form is
$$\omega=\sigma^*(d_Hy)=\sqrt{1+u_x^2}\,dx=ds$$
We notice that $I_2=\kappa$, the Euclidean curvature, $I_3=\kappa_s=\frac{d\kappa}{ds}$ but 
$I_4=\kappa_{ss}+3\kappa^3$, according to recurrence formulas in \cite{FO99}.
As we have seen in Example~\ref{ex1} the first two invariants $\kappa$ and $\kappa_s$ are sufficient to solve the equivalence problem for curves in Euclidean space.
\begin{rem} The action of the special Euclidean group on $J^1$ is locally free, however
the isotropy group of each point $(x,u,u_x)$ contains one nontrivial transformation:
\begin{equation}\label{is}
\left(
\begin{array}{ccc}
-1&0&2x\\
0&-1&2u\\
0&0&1
\end{array}\right).
\end{equation}
This is reflected in the ambiguity of normalization for $\alpha$, which is defined up
 to the addition of $\pi n$. Thus the invariants $I_2, I_3,\dots$ are local. 
In particular $I_2=\kappa$ changes its sign under the transformation (\ref{is}). 
\end{rem}
We conclude this example with the discussion on how the procedure described above corresponds to the classical definition of the Euclidean curvature.  The Fr\'enet frame consists of the unit tangent $T$ and the  unit normal $N$ attached to each point on a curve $(x,u(x))$, with consistent orientation of the frame.   This produces  a map 
$\rho_l:J^1(\R^2,1) \rightarrow SE(2)$. Indeed, the pair $(T,N)$ defines
 a rotation matrix at each
 point, and a point $(x,u(x))$ defines a translation vector at each point. 
This map is equivariant 
with respect to the  action of the group on itself by {\it left} multiplication, and 
so it is called a {\it left} 
moving frame. The element of the Euclidean group assigned to each point can be written 
in  matrix form: 
$$\rho_l=\left(
\begin{array}{ccc}
\frac 1{\sqrt{1+u_x^2}}&-\frac {u_x}{\sqrt{1+u_x^2}} & x\\
\frac{u_x} {\sqrt{1+u_x^2}}&\frac {1}{\sqrt{1+u_x^2}} &u\\
0&0&1
\end{array}
\right)
$$
In this example as well as in general the {\it left} moving frame is the group  inverse of the corresponding  {\it right} moving frame: $\rho_l=(\rho_r)^{-1}$. 
In matrix form the Fr\'enet equations can be written as:
$$
(\frac{dT}{ds},\frac{dN}{ds},\frac{dX}{ds})=(T,N,X)\left(
\begin{array}{ccc}
0& -\kappa  & 1\\
\kappa &0 &0\\
0&0&0
\end{array}
\right),
$$
where $X=(x,u,1)^t$. We recall that $\rho_l=(T,N,X)$ and thus 
\[\rho_l^{-1}\frac d {ds} (\rho_l)\, ds=\left(
\begin{array}{ccc}
0& -\kappa ds & ds\\
\kappa ds&0 &0\\
0&0&0
\end{array}
\right)\]
is  a pull-back of the Maurer-Cartan forms on $SE(2)$ to the jet of a curve 
under the $G$-equivariant map 
$\rho_l$.  The elements of this matrix  
  are  contact invariant forms and their  ratio produce the second order 
differential invariant $\kappa$. 
\end{ex}

\chapter{Inductive Construction of Moving Frames.}
    \label{chInd}

In this chapter we present two  modifications of the moving frame method which
were motivated by its  practical implementation.
We call the first modification a recursive construction  because in contrast to the
algorithm from chapter one it allows to construct differential  invariants order
by order. At each step we normalize more and more of the group parameters at the end obtaining a moving frame for the group $G$.
 In the next
chapter we use this algorithm to construct a complete set of the differential
invariants for ternary cubics transformed by linear changes of variables.

The second modification can be used when the group $G$ factors as a product of
two subgroups: $G=AB,$ such that $A\cap B$ is discrete. In this case
invariants and moving frames for $A$ and $B$ can be used to construct invariants
and a moving frame for $G$. Such approach not only simplifies the
computations but also produces as a by-product the relations among the
invariants of $G$ and its subgroups. We notice that the group  of the
Euclidean motions on the plane is a factor of the group of the special
affine motions. In its turn the special affine motions is a factor of the group of projective transformations on the plane.  We use these  groups to illustrate
our algorithm.  As a consequence  we obtain the affine curvature in terms of the Euclidean and the projective curvature in terms of the special affine and also the relations among corresponding differential operators.

\section {Recursive Construction of Moving Frames.}

We  first sketch  the main idea of the algorithm.
Assume that $G$ acts on $M$ regularly but not freely. We notice that if
$\CS\subset M$ is a
cross-section to the orbits of $G$ and a smooth  map $\rho\colon M \rightarrow G$ is defined by the condition $\rho(z)\cdot z\in \CS$ for $z\in M$, then the non-constant coordinates of $\rho(z)\cdot z$ provide a complete set of zero order invariants. However in order to build
a moving frame recursively  we require that  $\CS$  satisfies a certain condition, namely
each point of $\CS$ has the same isotropy group. We postpone the discussion of the
existence of such a cross-section until later   and assume that $H_1$ is the  isotropy
group of each point in
$\CS$. In this case the map $\rho_0$ defined by the
 condition
$\rho_0(z)\cdot z\in \CS$ is  a $G$-equivariant map from $M$ to the right cosets
$H_1\backslash G$. We prolong the action of $G$ to the first order and define the set
$\CS^1=\{z^{(1)}|\pi_0^1(z^{(1)})\in \CS\}$. The set $\CS^1$ is invariant under the
action of $H_1$ and we assume that there is a cross-section $\CS^1_1\subset\CS^1$ with a
constant isotropy group $H_2$. We use this cross-section to define a map
$\rho_1:\CS^1\rightarrow H_2\backslash H_1$, by the condition
 $\rho_1(z^{(1)})\cdot z^{(1)}\in \CS^1_1$. The map
$\rho_1\left(\rho_0(z^{(1)})\cdot z^{(1)}\right)\rho_0(z^{(1)})$ is a
$G$-equivariant map
from $J^1$ to $H_2\backslash G$. The non-constant coordinates of
$\rho_1\left(\rho_0(z^{(1)})\cdot z^{(1)}\right)\rho_0(z^{(1)})\cdot z^{(1)}$ provide a
complete set of the first
order invariants. We continue this process by prolonging  the action
of $H_2$ to the next order, or we may prolong  by several orders at once if
we wish. The algorithm terminates at the order where the isotropy group becomes trivial.
For regular jets this happens
 at the order of stabilization  of the group. This procedure resembles in many ways
the algorithm presented by M. Green \cite{Green78}  for  constructing moving frames
for curves in homogeneous spaces, however taking advantage of the generalized
approach by Fels and Olver \cite{FO99},  we  can apply our algorithm
to construct a moving frame for submanifolds of any dimension under  more general
(not necessarily transitive) group actions.

Following
\cite{GOV93} we call a cross-section with a constant isotropy group a {\it slice. }
The obvious necessary condition for the existence of a slice in a
neighborhood $U$ of $z_0\in M$ is that the isotropy groups of any two
points in $U$ are conjugate by an element of $G$.
There is a simple counterexample of this phenomenon (\cite{VP89}, Example~1,\S~7).
\begin{ex} Let $\R^2$ act on $\R^2$ by
$$x \rightarrow x+au+b, \qquad
 u \rightarrow u
$$
The isotropy group of a  point $(x_0,u_0)$ is defined by the condition $au_0+b=0$. The orbits are lines parallel to the  $x$-axis. All the points that lie on the same orbit have equal isotropy groups. On the other hand  the isotropy groups of two points from different orbits are not equal  and they are not conjugate because the group is commutative.
\end{ex}

However, one can find  slices  for many common group actions. In
particular slices exist if the group action is proper
\cite{GOV93}.
\begin{df}
The action of $G$ is called {\it  proper} if the map $\theta: G\times M\rightarrow M\times M$ defined by $\theta(g,z)=(g\cdot z,z)$ is proper. In other words if  $K\subset M\times M$ is compact then so is  $\theta ^{-1}(K)\subset  G\times M$.
\end{df}
Since $\theta^{-1}(z,z)=(H,z)$ where $H$ is an isotropy group of
$z$, if the action is proper, then  the isotropy group of each
point is compact.

\begin{prop}
Let $G$ act regularly on a manifold $M$ and let $\CS$ be a slice.
Then the condition $\rho(z)\cdot z\in \CS$ for $z\in M$ defines a map
$[\rho] :M\rightarrow H\backslash G$
which is $G$-equivariant.
\label{Geq}
\end{prop}
{\it Proof.}
Let group elements $g_1$ and $g_2$ be such that $g_1\cdot z\in \CS$ and  $g_2\cdot z\in \CS$. Since each orbit intersects the slice $\CS$ at one point, then: $g_1\cdot z=g_2\cdot z$ and so
$ g_2^{-1}g_1$ belongs to the isotropy group $G_z$ of $z$. On the other hand since $H$ is the isotropy group of $g_2\cdot z$, then  $G_z=g_2^{-1}Hg_2$. Thus
$ g_2^{-1}g_1\in g_2^{-1 }Hg_2 $ and hence $ g_1\in Hg_2 $. We have proved that the  map $[\rho]$ is well defined.

To show the equivariance of $[\rho]$ we need to prove that
$[\rho] (g\cdot z)=[\rho](z)g^{-1}$.
Let us choose  $q\in [\rho](z)$  and  $\q\in [\rho](g\cdot z)$
By the construction of $[\rho] $
one has
\[
\q g\cdot z=q\cdot z \in {\CS}.
\]
It follows that $q^{-1}\q g\in G_z=q^{-1}Hq$, or equivalently
$$\q\in Hqg^{-1}.$$
Since $[\rho](z)=Hq$ and $[\rho] (gz)=H\q$ we have proved that
$[\rho] (g\cdot z)=[\rho](z)g^{-1}$.
\qed

We can extend $[\rho] $ to a $G$-equivariant map on $J^{k }$ by
$[\rho](\zk)=[\rho] \left(\pi_0^k(\zk )\right)$.
Locally we can choose  a section $s\colon H\backslash G\rightarrow G,$ such that
$s[H]=e$ and define the map
\begin{equation}\label{rhos}
\rho _s=s\circ [\rho]: J^k(M,p)\rightarrow G.
\end{equation}
The map  $\rho _s$ is not $G$-equivariant, but it is {\it $G$-equivariant up to the left action of $H$}, and so the functions
$f\left(\rho_s( \zk)\cdot \zk)\right)$, where $f$ is any function  on $J^k$ are invariant up
to the left action of $H$.
\begin{prop}\label{invupto}
There exists element $h\in H$ such that $$\rho_s(g\cdot \zk)=h\rho_s(\zk)g^{-1}.$$
\label{GHeq}
\end{prop}

{\it Proof}.
Since  $[\rho]$  is a $G$-equivariant map from  $J^k$ to $ H\backslash G$ then
$$\rho_s(g\cdot \zk)=s[\rho(g\cdot \zk)]=s[\rho(\zk)g^{-1}]=s[H \rho_s(\zk)g^{-1}]\in H\rho_s(\zk)
g^{-1}.$$
Thus there  exists  $h\in H$ such that
$\rho_s(g\cdot \zk)=h\rho_s(\zk)g^{-1}$.
\qed

Let $\CS^k=\{\zk|\pi^k_0(\zk)\in\CS\}$ be the subset of $J^k$ which projects to $\CS$.
In other words the  set $\CS^k$ is a pull-back of the fiber bundle
$J^k\rightarrow M$ under the
inclusion  $\CS\rightarrow M$.
By construction  $\CS^k$ is invariant under the prolongation of the $H$-action.

Let the map $w:\CB^k\rightarrow J^k$ be defined by the prolonged group action, that is
$w(g,\zk)=g\cdot \zk$, the map $\sigma_0\colon J^k\rightarrow \CB^k$ be defined by
$\sigma_0(\zk)=(\rho_s(\zk),\zk)$ and the map $\iota_0: J^k\rightarrow J^k$
be defined by $\iota_0(\zk)=\rho_s(\zk)\cdot \zk$ for $k=n,\dots,\infty$,
then we obtain a commutative diagram
similar to the one in Section \ref{smf} of Chapter~\ref{prelim}:
$$\Atriangle<-1`1`1;600>[{\CB^k=G\times J^k}`{J^k}`{J^k};
\sigma_0`w`{\iota_0}]
$$

If $f$ is any function on $J^k$ then due to the  proposition above the  function
$\iota_0^*(f)$ is invariant under $G$ up to an element of $H$.  We note  that $\iota_0$ projects  $ J^k$ to $ \CS^k$
 and thus   non-constant coordinates of  $\iota_0(\zk )$, restricted to $\CS^k$ are
transformed exactly  in the same way as coordinate functions on $\CS^k$.
More generally if $f$ is any function on $J^k$ then
$\iota_0^*(f)$ restricted to $\CS^k$ equals to $f$ restricted to  $\CS^k$
and so as functions on $\CS^k$ they are transformed exactly in the same way.
 This  trivial observation is used  in
 both algorithms and so we state it as a proposition.

\begin{prop}
\label{thesame}
Assume that  $\CS\subset N$ be a submanifold of $N$ invariant under the  action of  the
group $H$ and  there is a smooth projection  $\iota:N\rightarrow \CS$, that is
 $\iota(z)\in\CS$ for any point $z\in N$ and
$\iota(\tz)=\tz$ for any point $\tz\in \CS$. Let $f(z)$ be function on $N$ and let $\tf$ be its restriction to $\CS$. We define a function $F(z)$ by the formula
$F(z)=f\left(\iota (z)\right)$ and denote its restriction to $\CS $ as $\tF$.
Then $\tF=\tf$ and hence $h\cdot \tF=h\cdot \tf$ for all $h\in H$.
\end{prop}
\begin{rem} Note that the proposition above {\em does not} assert that
 $h\cdot F=h\cdot f$ on $N$, or equivalently, in general
$f\left(\iota(h\cdot z)\right)\neq f(h\cdot z)$.
\end{rem}

The following proposition asserts that $\CS^n\cap \CV^n
\neq\emptyset $ and hence one can construct a moving frame:
\begin{equation}\label{rhoH}
\rho_H:\CS^n\rightarrow H.
\end{equation}
\begin{prop}
Assume that $\CS$ is a cross-section for an
action of a group  $G$ on
a manifold $M$. Let  $n$ be the order of stabilization, $\CV^n$ be the regular set and
$\CS^n=\{\zn|\pi^n_0(\zn)\in\CS\}\subset J^n(M,p)$, then
$\CS^n\cap \CV^n\neq \emptyset$.
\end{prop}
{\em Proof.}
 The statement $\CS^n\cap \CV^n=\emptyset$ implies that
$\pi^n_0(\CV^n)\cap \CS=\emptyset$, that is there are no
$p$-dimensional
submanifolds of $M$  passing through any  point  $z\in\CS$ such that its $n$-th
prolongation
$j^n(S)$ at $z$ belongs to $\CV^n$.  But this  means that  for any $g\in  G$,
 there are no submanifolds through the point $g\cdot z$
that give rise to a regular $n$-th jet. Thus  $\{G\cdot \CS\}\cap
\pi^n_0(\CV^n)=\emptyset$. The set ${G\cdot \CS}$ is an open
subset of $M$ and hence  we arrive to a contradiction with the
assertion that the set $\CV^n$ is dense in $J^n$ if the action of
$G$
 \cite{Ovs82}. \qed
\begin{rem}
Although the map $\rho_H:\CS^n\rightarrow H$ is $H$-equivariant, its extension on $J^n$ defined by
$\rho_H\left(\iota_0(\zn)\right)$
is {\it not} an $H$-equivariant map from $J^n$ to $H$.
\end{rem}

\begin{prop}
The map $\rho_G\colon J^n(M,p)\rightarrow G$ defined by
\begin{equation}\label{rhoG}
\rho_G(\zn)=\rho_H\left(\iota_0(\zn)\right)\rho_s(\zn).
\end{equation}
is $G$-equivariant.
\end{prop}
{\it Proof.}
We recall that $\iota_0(\zn)=\rho_s(\zn)\cdot \zn$. From Proposition~\ref{GHeq} we know
that for an element $g \in G$ there is an element $h\in H$ such that
$$\rho_s(g\cdot \zn)=h\rho_s(\zn)g^{-1}.$$
and so
$$ \iota_0(g\cdot \zn)=h\rho_s(\zn)\cdot \zn=h\cdot \iota_0(\zn).$$
Thus
\begin{eqnarray*}
\rho_G(g\cdot \zn)&=&\rho_H\left(\iota_0(g\cdot \zn)\right)\rho_s(g\cdot \zn)\\
&=&\rho_H\left(h\cdot \iota_0(\zn)\right)h\rho_s(\zn)g^{-1}=
\rho_H\left(\iota_0(\zn)\right)\rho_s(\zn)g^{-1}.
\end{eqnarray*}
In the last equality we have used $H$-equivariance of the map $\rho_H$.
\qed

Since $\rho_G$ is a moving frame, then non-constant coordinates of
$\iota_G=\rho_G(\zn)\cdot \zn$ provide a complete set of $n$-th order differential
invariants.

Let $[g]$  denote the equivalence class of $g$ in $H\bs G$,
let $U_{H\backslash G}$ be a neighborhood of $H$ in ${H\backslash G}$ where
the local section
$s\colon {H\backslash G}\rightarrow G$ is defined, and let $e\in U_G\subset G$
be the preimage  of $U_{H\backslash G}$ under the canonical projection
$G\rightarrow H\bs G$.
Then the following two maps are local diffeomorphisms:
\begin {eqnarray*}
\phi\colon U_G\times J^\infty \rightarrow H\times U_{H\bs G}\times J^\infty&\colon&
\phi(g,\zi)=\left(g(s[g])^{-1}, [g],\zi\right),\\
\psi\colon  H\times U_{H\bs G}\times J^\infty\rightarrow &\colon&
\psi(h,[g],\zi)=\left(h\, s[g],\, \zi\right).
\end{eqnarray*}
We can summarize
the recursive  construction in the following commutative diagram:
\vskip5mm
$$\bfig
\putsquare<1`-1`1`1;1000`700>(0,0)[{H\times \Ji}`{H\times \Ji}`\Ji`\Ji;{\iota_0}
`\sigma_H`w_H`\iota_G]
\putAtriangle<-1`1`0;500>(0,700)[{H\times U_{H \bs G}\times \Ji}`
\phantom{H\times \Ji}
`\phantom{H\times \Ji};\sigma_0`w_0`]
\putmorphism(-700,1200)(1,0)[{U_G\times \Ji}``\phi]{800}1a
\putmorphism(1000,1200)(1,0)[`{G\times \Ji}`\psi]{800}1a
\putmorphism(-90,80)(-1,2)[\phantom{ \Ji}`\phantom{U_G\times \Ji}`\sigma_G]{480}{1}l
\putmorphism(1090,80)(1,2)[\phantom{ \Ji}`\phantom{U\times \Ji}`w_G]{480}{-1}r
\efig
$$
where the maps $w$ are defined by the group actions:
\begin{eqnarray*}
w_0(h, [g],\zi)&=&(h,s[g]\cdot \zi),\\
w_H(h, \zi)&=&h\cdot \zi,\\
w_G(g, \zi)&=&g\cdot \zi.
\end{eqnarray*}
The maps $\sigma$ are defined from the maps $\rho$ (see formulas (\ref{rhos},
\ref{rhoH}, \ref{rhoG}):
\begin{eqnarray*}
\sigma_H(\zi)&=&\left(\rho_H(\rho_s (\zi)\cdot \zi),\zi\right),\\
\sigma_0(h,\zi)&=&\left(h,[\rho_s(\zi)],\zi\right ),\\
\sigma_G(g,\zi)&=&\left (\rho_G(\zi),\zi\right ),
\end{eqnarray*}
and the maps  $\iota$ are projections:
\begin{eqnarray*}
\iota_0(h,\zi)&=&(h,\, \rho_s(\zi)\cdot \zi), \\
\iota_G(z)&=&\rho_H(\rho_s(\zi)\cdot\zi)\cdot \zi.
\end{eqnarray*}
\begin{ex}

Let the special rotation group  $G=SO(3,\R)$ act on $M=\R^3\backslash \{0\}\times \R$
by rotations on the independent variables $x,y,z$ and the trivial
action on the dependent variable $u$.

The action of  $SO(3,\R)$  on $\R^3$ is not free. The isotropy group of every  point on
 the  positive half of the $z$-axis consists of the rotations around the $z$-axis. On the other hand each orbit of  $SO(3,\R)$ intersects the positive half of  the $z$-axis at the unique point, and hence it can serve as a slice $\CS$ with the isotropy group
$H \backsimeq SO(2,\R)$

Our first step  is to  construct the
map $\rho_s :M\rightarrow G$ such that
 $\rho_s(p)\cdot p\in \CS$ for each  $p\in M$ and  $\rho_s (p_0)=\Id\,$ for any
$p_0\in \CS$.

Each coset $G\backslash H$ can be represented by the product of two rotations: $R_x(\theta) $
 with respect to the $x$-axis and  $R_y(\tau)$  with respect to the $y$-axis.
In matrix form this can be written as
\begin{eqnarray*}
&&\left(
\begin{array}{ccc}
\cos \tau & 0 & -\sin \tau \\
0& 1 & 0 \\
\sin \tau & 0 & \cos \tau
\end{array}
\right) \left(
\begin{array}{ccc}
1 & 0 & 0 \\
0& \cos \theta & -\sin \theta \\
0& \sin \theta & \cos \theta
\end{array}
\right) \\
&=&\left(
\begin{array}{ccc}
\cos \tau & -\sin \tau \sin \theta & -\sin \tau \cos \theta \\
0 & \cos \theta & -\sin \theta \\
\sin \tau & \cos \tau \sin \theta & \cos \tau \cos \theta
\end{array}
\right)
\end{eqnarray*}
We choose the  first rotation   $R_x(\theta) $
so that it  brings an arbitrary point $p=(x,y,z)$ to the upper $xz$-plane.
It can be achieved by choosing $\theta=\arctan(\frac y z)$, then
$\tilde{p}=R_x(\theta)\cdot p=(x,0,\sqrt{z^2+y^2}) $.  We choose   the rotation   $R_y(\tau)$
so that it  brings   $\tilde{p} $ to the $z$-axis.
We take  $  \tau =\arctan(\frac x{\sqrt{z^2+y^2}})$ and then
$R_y(\tau)\cdot \tilde{p}=(0,0,\sqrt{z^2+y^2+x^2})$ lies on the positive $z$-axis.
In the matrix form
\[
\rho_s(p)=R_y\, R_x=\left(
\begin{array}{ccc}
\frac{\sqrt{z^2+y^2}}r & -\frac{yx}{r\sqrt{z^2+y^2}} & -\frac{zx}{r\sqrt{z^2+y^2}} \\
0 & \frac z{\sqrt{z^2+y^2}} & -\frac y{\sqrt{z^2+y^2}} \\
\frac xr & \frac yr & \frac zr
\end{array}
\right) ,
\]
where $r=\sqrt{z^2+y^2+x^2}.$ The only non-constant   coordinate of
$\rho_s(p)\cdot p=(0,0,r)$
provides a zero order invariant, the radius.
We note that $\rho_s(0,0,z)=\Id$\, for $z>0$  and hence the conditions of
Proposition~\ref{thesame} are satisfied.

On the  next step we  prolong the transformation to the first jet
bundle $J^1$. We consider the set
$\CS^{1}=\{(u,x,y,z,u_x,u_y,u_z)|x=0,y=0,z>0\}\subset J^1,$ which projects to $\CS$ and is  invariant under the first prolongation of the $H$ action.
The first prolongation of the $\rho_s(p)$ transforms a point
$(u,x,y,z,u_x,u_y,u_z)\in J^1$ to a point $(U,X,Y,Z,U_x,U_y,U_z)$ in $\CS^{1}$
where

\begin{eqnarray*}
u\rightarrow U &=&u \\
x\rightarrow X &=&0 \\
y\rightarrow Y &=&0 \\
z\rightarrow Z &=&\sqrt{z^2+y^2+x^2}=r \\
u_x\rightarrow U_x &=&\frac{\sqrt{z^2+y^2}}ru_x-\frac{yx}{r\sqrt{z^2+y^2}}u_y-\frac{zx}{%
r\sqrt{z^2+y^2}}u_z; \\
u_y \rightarrow U_y &=&\frac z{\sqrt{z^2+y^2}}u_y-\frac y{\sqrt{z^2+y^2}}u_z; \\
u_z\rightarrow U_z &=&\frac xru_x+\frac yru_y+\frac zru_z.
\end{eqnarray*}
From  Proposition~\ref{thesame} we know that functions $\{U,Z,U_x,U_y,U_z\}$
are transformed on $\CS^1$ by $H$ by the same formulas  as $\{u,z,u_x,u_y,u_z\}$.
Let us represent an element of $H$ by the matrix:
\[
\left(
\begin{array}{ccc}
\cos \alpha & -\sin \alpha &  \\
\sin \alpha & \cos \alpha &  \\
&  & 1
\end{array}
\right) ,
\]
then
\begin{eqnarray*}
\bar{U} &=&U \\
\bar{Z} &=&Z; \\
\bar{U}_x &=&\cos \alpha\, U_x-\sin \alpha\, U_y; \\
\bar{U}_y &=&\sin \alpha\, U_x+\cos \alpha\, U_y; \\
\bar{U}_z &=&U_z.
\end{eqnarray*}
We observe that $H$ acts freely on $\CS^1$. We find a moving frame $\rho_H$
by setting $\bar{U}_x=0,$ and hence $\tan \alpha =\frac{U_x}{U_y}$.
Substitution of this normalization into $\bar{Z},\bar{U}_y,\bar{U%
}_z,$ produce invariants of $G=SO(3)$:
\[
U,\quad Z,\quad I_y=\sqrt{U_x^2+U_y^2},\quad I_z=U_z;
\]

In terms of coordinates on $J^1$ they can be written as:
\begin{eqnarray*}
U &=&u \\
Z &=&\sqrt{z^2+y^2+x^2}=r; \\
I_y &=&\sqrt{\left( yu_z-zu_y\right) ^2+\left( zu_x-xu_z\right) ^2+\left(
xu_y-yu_x\right) ^2}; \\
I_z &=&\frac xru_x+\frac yru_y+\frac zru_z;
\end{eqnarray*}
The corresponding  moving frame for the group $G$ is a product:
\[
\rho_H\rho_s=
\left(
\begin{array}{ccc}
\frac {U_y}{\sqrt{U_x^2+U_y^2}} & -\frac{U_x}{\sqrt{U_x^2+U_y^2}} &  \\
\frac {U_x}{\sqrt{U_x^2+U_y^2}} & \frac {U_y}{\sqrt{U_x^2+U_y^2}} &  \\
&  & 1
\end{array}
\right)
 \left(
\begin{array}{ccc}
\frac{\sqrt{z^2+y^2}}r & -\frac{yx}{r\sqrt{z^2+y^2}} &
-\frac{zx}{r\sqrt{z^2+y^2}} \\
0 & \frac z{\sqrt{z^2+y^2}} & -\frac y{\sqrt{z^2+y^2}} \\
\frac xr & \frac yr & \frac zr
\end{array}
\right).
\]

\end{ex}
The action of $H$ became free on  the first jet in the example above. If it is not,
we would have to either prolong the action up to the order where it is free or
to repeat our algorithm for $H$ acting on $\CS^1$. This can be done if there is a slice
$\CS^1_1\subset \CS^1$ with a constant isotropy group $H_1\subset H$.
The following procedure describes  order by order construction of invariants under the assumption that slices exist at each order.
 In our  notation the superscripts refer to the order of
prolongation and the subscripts refer to the induction step.
\begin{alg}\label{AlgHi}
On the zeroth step we consider the action  of the group $G=H_0$
such that there is a  slice ${\CS_0}\subset M$ with a constant isotropy group
$H_1$. We define the map $[\rho_0]:M\rightarrow H_1 \bs G$ such that
$\rho_0(z)\cdot z\in \CS_0$ and  $\rho_0(\CS_0)=e$, where  $\rho_0(z)$
is defined by a local  section from  $H_1 \bs G$ to $G$.
The non-constant coordinates of
\[
\iota_0 =\rho_0 (z)\cdot z,
\]
provide
a complete set of functionally independent zero order invariants.
Let
${\CS_0}^k=\{z^{\left( k\right) }\in J^k|\pi_0^k(\zk)\in \CS_0\}$ and
$\rho_0(\zk)=\rho_0(\pi_0^k(\zk))$,
then the map $\iota_0:J^k\rightarrow \CS_0^k$ defined by
\[
\iota_0(\zk) = \rho_0(\zk) \cdot \zk
\]
is invariant under the action of $G$ up to an element of $H_1$
as described in  Proposition~\ref{GHeq}.

Let $\CS^1_1\subset\CS^1_0$ be a slice for the action of $H_1$ on $\CS_0^1$
such that the isotropy group of each point equals  $H_2$
and let the map
$[\rho_1]: \CS_0^1 \rightarrow H_2\bs H_1$ be defined by the conditions
$\rho_1(\iota_0(z^{(1)}))\cdot \iota_0(z^{(1)})\in \CS_1^1$ and $\rho_1(\CS_1^1)=e$,
 where  $\rho_1(z)$
is defined by a local  section from  $H_2 \bs H_1$ to $H_1$.
Then the non-constant coordinates of
\[
\iota_1 =\rho_1 \left(\iota_0(z^{(1)})\right)\cdot \iota_0(z^{(1)})
\]
provide
a complete set of functionally independent first  order invariants  (this set
includes zero order invariants).
Let
${\CS_1}^k=\{z^{\left( k\right) }\in J^k|\pi_1^k(\zk)\in \CS_1^1\}$ and
$\rho_1(\zk)=\rho_1(\pi_1^k(\zk))$
then the map $\iota_1:J^k\rightarrow \CS_1^k$ is defined by
\[
\iota_1(\zk) = \rho_1(\iota_0(\zk)) \cdot \iota_0(\zk) .
\]
The group  product
\begin{equation}\label{pr}
\rho_1\left(\iota_0(z^{(1)})\right)\rho_0(z^{(1)})
\end{equation}
 defines a $G$-equivariant map from $J^1$ to the right cosets $H_2\bs G$.
If the isotropy group  $H_2={e}$ on
$\CS_0^1$,  the coordinates of $\iota_1(\zk)$ are invariant and the product
(\ref{pr}) is a moving frame for $G$.
Otherwise, we have to prolong the action of $H_2$ on $\CS^1_1$ to the second order
and repeat the algorithm. In the case when
 $\CS_0^n=\{\zn\in J^n(M)| \pi^n_0(\zn)\in\CS_0\}$ belongs to the regular set $\CV^n$, this process will terminate at the order of
stabilization $n$. Indeed on the $n$-th step we consider the action of the
isotropy group $H_{n}$ on the set $\CS^n_{n-1}$.   Any element  of Lie algebra of
$G$ generates a nontrivial transformation on $\CV^n$
and hence the action of $H_n$ on
$\CS^n_{n-1}\subset \CS_0^n$  is free.  We choose a cross-section
$\CS^n_{n}\subset \CS^n_{n-1}$ transversal to the action of $H_n$ and construct a
corresponding $H_n$-invariant  map $\rho_n:  \CS^n_{n-1} \rightarrow H_n$. The
  map $\rho_G=J^n\rightarrow G$ defined by
$$\rho_G=\rho_n(\iota_{n-1}(\zn))\rho_{n-1}(\iota_{n-2}\zn)\dots
\rho_1(\iota_0(\zn))\rho_0(\zn)$$
is a moving frame for $G$.
The complete set of $n$-th order differential invariants is given by the non-constant coordinates of
$$\iota_n=\rho_n(\iota_{n-1}(\zn)) \cdot \iota_{n-1}(\zn).$$
\end{alg}
\vskip 10mm
In the next  chapter we will return to  this algorithm to construct a complete set
of differential invariants for ternary cubics.

\section {A Moving Frame Construction for a Group that Factors as a Product.}

The moving frame construction is simpler for a group with a smaller number of
parameters. Thus  it is desirable to use a moving frame for a subgroup of $G$ to construct a moving frame for $G$.  We say that a group $G$ factors as a product of its subgroups $A$ and $B$ if $G=AB$, that is for any $g\in G$ there are $a\in A$ and $b\in B$ such that $g=ab$.
We reproduce two   useful statements   from \cite{GOV93}.
 \begin{thm} Let $G$ be a group, $A$ and $B$
are two subgroups. Then the following conditions are equivalent:

a) the reduction of the natural action of $G$ on $G/B$ to $A$ is transitive,

b) $G=AB$,

c) $G=BA$,

d) the reduction of the natural action of $G$ on $G/A$ to $B$ is transitive.
\label{product}
\end{thm}
\begin{cor}
The reduction of the natural action of $G$ on $G/B$ to $A$ is free and transitive
if and only if $G=AB$
(or $G=BA$) and $A\cap B=e$.
\end{cor}
\begin{rem}
If $G=AB$ and  $A\cap B=e$ then for each $g\in G$ there are {\it unique} elements $a\in A$ and $b\in B$ such that $g=ab$.
In this case the manifold $A\times B$ is diffeomorphic to $G$ and we will write
$A\times B \backsim G$
to denote  that two Lie groups are diffeomorphic as manifolds but  are  {\em not} necessarily
isomorphic as groups.  In the case when $A\cap B$ is discrete  then
$A\times B$ is locally diffeomorphic  to $G$.
\end{rem}

 The following theorem plays a central role in the construction of a moving frame for a
 product of two groups.

\begin{thm} Let  $A$ and $B$ act regularly on a manifold $M$ and assume that in a neighborhood   of a point $z_0$
the infinitesimal generators of the $A$-action are linearly
independent from the generators of the $B$-action.   Then locally
there exists a submanifold $\CS_A$ through the point $z_0$, which
is transverse to the  orbits of the  subgroup $A$ and is invariant
under the action of the subgroup $B$. \label{AB}
\end{thm}
{\it Proof}. Let $a$ be the dimension of the $A$-orbits, $b$ be
the dimension of the $B$-orbits
 on $U$ and $m=\dim M$. By  Frobenius' theorem we can locally rectify
 the orbits of $B$, that is we can introduce coordinates
$\{y_1,\dots,  y_b,x_1,\dots,x_{m-b}\}$  such that the orbits of
$B$ are defined by the equations $x_i=k_i, \, i=1,\dots, {m-b}$,
where $k_i$ are some constants. The orbits of $B$ are integral
manifolds for the distribution $\{\frac \partial {\partial
y_1}\dots \frac \partial {\partial y_b}\}$. The functions $x_i$
are invariant under the $B$-action. Let vector fields
$X_1,\dots,X_a$ and $Y_1,\dots,Y_b$ be a basis for infinitesimal
generators of the action of $A$ and $B$ respectively in a
neighborhood $U$ containing $ z_0$ .  The vector fields $Y_i,
i=1,\dots,b$ and $\frac \partial {\partial x_j},\, j=1\dots{m-b}$
are linearly independent by the choice of coordinates and their
union forms a basis in $TU$.
 We can choose $c=m-b-a$ vector fields
$\frac \partial {\partial x_{j_1}}\dots \frac \partial {\partial
x_{j_c}}$
  which  are linearly
independent from $X_1,\dots,X_a$ in $TU$. Let $\CS_A$ be an
integral manifold through the point $z_0$  for the involutive
distribution
 $\Delta=\{\frac \partial {\partial x_{j_1}}\dots \frac \partial
 {\partial x_{j_c}}, \frac \partial {\partial y_1}\dots \frac \partial {\partial y_b}\}$.
By construction  $\CS_A$ is invariant under the action of $B$ (it
is a union of the orbits of $B$). On the other the distribution
$\Delta$ is transversal to the infinitesimal generators
$X_1,\dots,X_a$  of the $A$-action, and so is transversal to the
orbits of $A$.
 \qed

With this result we construct a moving frame for a product of groups $A$ and be
$B$ as follows.
\begin{alg}\label{algAB}
Let  $G=BA$ and let $B\cap A$ be discrete.
 Since we are constructing a local moving frame,
 that is  a map to a neighborhood of the identity of the group we may assume that
$B\cap A=e$. Thus an element $g\in G$ can be written as a product $g=ba,\, a\in A,\, b\in B$
and $G\backsim B\times A$ in the category of smooth manifolds (but not as groups).
Let $n$ be the order of stabilization of the $G$-action. Since both $A$ and $B$ act
freely on $\CV^n\in J^n$ and their intersection is trivial then the infinitesimal
generators of the $A$-action and the $B$-action are linearly independent at each point
of $\CV^n$ and hence they satisfy the conditions of Theorem~\ref{AB}.
Thus  there is a
cross-section $\CS_A\subset \CV^n$ for the action of $A$ which is invariant under
the action of $B$. We use this cross-section to construct a moving frame $\rho_A$ for
$A$. The map  $\iota_A=\rho_A(\zn)\cdot \zn$ projects $\CV^n$ on the cross-section
 $\CS_A$, which is invariant under the action of $B$.  Moreover
 the action  of $B$ on $\CS_A$ is locally free and
hence we can choose a cross-section $\CS\subset \CS_A$ that defines   a moving frame
$\rho_B: \CS_A\rightarrow B$. We can extend $\rho_B$ to a map
$\tilde\rho_B:\CV^n\rightarrow B $, by the formula
\begin{equation}\label{rhoB}
 \tilde\rho_B=\rho_B\left(\rho_A(\zn)\cdot \zn\right).
\end{equation}
The map $ \tilde\rho_B$ is $A$-invariant but, in contrast to $\rho_B$,
it is not $B$-equivariant.

The cross-section $\CS$ is transversal to the orbits of $G$
and the map $\rho_G$ defined by
\begin{equation}\label{mfAB}
\rho_G(\zn)=\tilde \rho_B\left( \zn\right)\rho_A(\zn)
\end{equation}
satisfy the condition  $\rho_G(\zn)\zn\in \CS$, and hence is a moving frame for the
$G$-action. The maps
$$\iota_G^k(\zk)=\tilde \rho_B\left( \zn\right)\rho_A(\zk)\cdot \zk$$
define projections of $J^k$ onto $\CS^k$, for $k=n,\dots,\infty$ and so
 the non-constant coordinate functions of
$$\rho_B\left(\rho_A(\zk)\cdot \zk\right)\rho_A(\zk)\cdot \zk \label{eqAB}$$
provide a complete set of $k$-th order differential invariants for G.
\end{alg}
\begin{rem} We notice that the coordinates of  $\rho_A(\zk)\cdot \zk$ are invariant
under the $A$-action  and thus the formula above expresses   the invariants of the
 $G$-action in terms of the invariants of its subgroup $A$.
\end{rem}
We can summarize our construction in the following commutative diagram:
$$\bfig
\putsquare<1`-1`1`1;1000`700>(0,0)[{B\times \Ji}`{B\times \Ji}`\Ji`\Ji;{\tilde\iota_A}
`\sigma_B`w_B`\iota_G]
\putAtriangle<-1`1`0;500>(0,700)[{B\times A\times \Ji}`\phantom{B\times \Ji}
`\phantom{B\times \Ji};\sigma_A`w_A`]
\putmorphism(-700,1200)(1,0)[{G\times \Ji}`{\backsim}`]{500}0a
\putmorphism(1200,1200)(1,0)[{\backsim}`{G\times \Ji}`]{400}0a
\putmorphism(-90,80)(-1,2)[\phantom{ \Ji}`\phantom{G\times \Ji}`\sigma_G]{480}{1}l
\putmorphism(1090,80)(1,2)[\phantom{ \Ji}`\phantom{G\times \Ji}`w_G]{480}{-1}r
\efig
$$
where the maps $w$ are defined by the prolonged group action for $k=1,\dots,\infty$:
\begin{eqnarray*}
w_A(b,a,\zk)&=&(b,a\cdot \zk),\\
w_B(b,\zk)&=&b\cdot \zk,\\
w_G(g,\zk)&=&g\cdot \zk=w_B\circ w_A(b,a,\zk)\,\, \mbox{where}\,\,g=ba.
\end{eqnarray*}
The maps $\sigma$ are defined from moving frames for $A,B$ and $G$ for
$k=n,\dots,\infty$ (see formulas (\ref{rhoB}, \ref{mfAB})):
\begin{eqnarray*}
\sigma_A(b,\zk)&=&(b,\rho_A(\zk),\zk)\\
\sigma_B(\zk)&=&(\tilde\rho_B(\zk),\zk),\\
\sigma_G(\zk)&=&(\rho_G(\zk), \zk)=\sigma_A\circ \sigma_B(\zk).
\end{eqnarray*}
The maps $\iota$ are projections:
\begin{eqnarray*}
\tilde\iota_A(b,\zk)&=&(b,\rho_A(\zk)\cdot \zk)\,:B\times J^k\rightarrow B\times \CS_A\\
\iota_G(\zk)&=&\rho_G(\zk)\cdot \zk\,:J^k\rightarrow \CS
\end{eqnarray*}
We remind the reader that although all maps are written as global they might
be only defined  on
an open subset of $J^k$ and in neighborhoods of the identities of
the groups $A,B$ and $G$. The manifolds $B\times A\times J^k$ and $G\times J^k$
are diffeomorphic, and this diffeomorphism is $A$-equivariant.
 The maps $w_A$, $w_G$, $\tilde\iota_A$, and $\iota_G$  are $A$-invariant,
  whence the maps $\sigma_A$ and $\sigma_B$ are $A$-equivariant,
  with respect to the action defined by:
\begin{eqnarray*}
\tilde a\cdot (b,a,\zk)&=& (b,a\tilde a^{-1},\tilde a\cdot \zk),\\
\tilde a\cdot (g,\zk)&=& (g\tilde a^{-1},\tilde a\cdot \zk),\\
\tilde a\cdot (b,\zk)&=& (b,\tilde a\cdot \zk).
\end{eqnarray*}
We note that neither $\sigma_A$ nor $\sigma_B$ is $B$-equivariant, but their composition
is.
As  was shown in Section~\ref{smf} (see formula (\ref{cinvc})) the forms
\begin{equation}
\omega_G^i=\sigma_G^*\,d_H\,w_G^*(x^i),\,\, i=1,\dots,p
\end{equation}
form a horizontal contact invariant coframe on $J^{\infty}$. Since $A$ is a subgroup of
$G$ then  the forms $\omega_G^i$ retain their invariant properties under the action
of $A$.
On the other hand the moving frame $\rho_A$ provide us with another
horizontal coframe which is contact invariant under the action of $A$:
$$\omega_A^i=\sigma_A^*\,d_H\,w_A^*(x^i),\,\, i=1,\dots,p.$$
The two coframes are related by a linear transformation
$w_G^i=\sum{j=1}^p L^i_jw_A^j$, where $L^i_j$ are functions on
$J^\infty$ invariant  under the $A$-action. In fact,   $L^i_j$ can
be explicitly expressed in terms of the basis invariants of $A$,
indeed:
\begin{eqnarray}\label{ABcinvc}
\omega_G^i&=&\sigma_B^*\,\sigma_A^*\,d_H\,w_A^*\,w_B^*(x^i)\\
\nonumber &=&\sigma_B^*\sigma_A^*\pi_Hw_A^*d_H\chi^i(b_1,\ldots,b_l,x^1,\ldots,x^p,u_J^\alpha)=
\sigma_B^*\inv_A(d_H\chi^i),
\end{eqnarray}
where $\chi^i=w_B^*x^i$ is a function on $B\times \Ji$, written in local coordinates and
$\inv_A$ denotes (contact) invariantization with respect to the $A$-action.
The forms $\tilde\omega_A^i=\inv_A(d_H\chi)=\sigma_A^*\,\pi_H\,w_A^*\,d_H\chi^i$ provide a horizontal
coframe on $B\times \Ji$ which is contact invariant with respect to the action of $A$.
A  form $\tilde\omega_A^i$ is obtained from  $d_H\chi^i$ by replacing forms $dx^j$
with  $\omega_A^j$ and coordinate functions $x^1, ,\ldots,x^p,u_J^\alpha$
with the corresponding fundamental invariants
$\CX^{(A)1}, ,\ldots,\CX^{(A)p},I_J^{(A)\alpha}$.
The final pull-back $\sigma_B^*$ is equivalent to the replacement of parameters
$b_1,...,b_l$ with the corresponding coordinates of
$\rho_B(\rho_A(z^{(\infty)})\cdot z^{(\infty)})$. The latter are expressed in terms of invariants of the $A$-action.

In many situations the following reformulation of Theorem~\ref{product}
enables us to enlarge a moving frame for a transformation group $A$  to a moving frame
for a larger group containing $A$.

\begin{thm}
Let ${\cal O}\subset M $ be an orbit of  $G$  and  let $A$ be a subgroup which acts transitively on ${\cal O}$. Then $G=HA,$ where $H$ is the  isotropy group of a point in
${\cal O}$. If in addition $A$ acts locally freely on ${\cal O}$ then $A\cap H$ is
discrete.
\end {thm}

Let $A$ act regularly on $M$ and let $n_A$ be  the order of stabilization
for  $A$,  then the  action of $A$ is (locally) free  on a subset
$\CV_A\subset J^{n_A}(M,p)$.
Assume that the action of $A$ can be extended to the action of a group $G$,
containing $A$ so that there is a point $z_0\in \CV_A$ such that the orbits
of $A$ and $G$  through $z_0$ coincide. If this is the case then let $H$ be the
isotropy group of the  point $z_0$.
Due to the theorem above  $G=HA$ and  $A\cap H$ is
discrete and so Algorithm~\ref{algAB} can be applied.

The situation is especially favorable if the action of $A$ on the regular set
$\CV_A\subset  J^{n_A}(M,p)$ is transitive.  In this case we can extend a moving
frame for $A$ to
 a moving frame for any group  $G$ containing $A$.
Let $H$ be an isotropy group for a
point $z_0 $  in  $\CV_A$, then as above $G=HA$ and the intersection $A\cap H$ is
finite. Moreover the point $z_0$ is  a cross-section to the action of $G$  which is invariant under the action of $H$. We use $z_0$ to define a moving frame
$\rho_A: J^{n_A}\rightarrow A$. Let $n$ be the order of stabilization for the
$H$-action and let
 $\CS_A^n=\{\zn|\pi^n_0(\zn)=z_0\}$. The non-constant coordinates of
$\rho_A(\zn)\cdot \zn$ are invariant under the action of $A$ and,
restricted to $\CS_A^n$ they are  transformed by $H$ in the same
way as the coordinate functions. Since $H$ acts locally freely on
$\CS_A^n $ there is a local moving frame
$\rho_H:\CS_A^n\rightarrow H$. The   map $\rho_G: J^n \rightarrow
G$ defined by
$$\rho_G(\zn)=\rho_H\left(\rho_A(\zn)\cdot \zn\right)\rho_A(\zn)$$
is a local moving frame for the $G$-action and  the non-constant coordinate functions of
$$\rho_H\left(\rho_A(\zk)\cdot \zk\right)\rho_A(\zk)\cdot \zk$$
provide a complete set of $k$-th order differential invariants of $G$.

\section{Examples:  Euclidean, Affine and Projective Actions on the Plane.}

The group  of the
Euclidean motions on the plane is a factor of the group of the special
affine motions. In its turn the group of  special affine motions is a factor of the group of projective transformations on the plane. All three of these transformation groups
 play an important role in the computer image processing \cite{Fu94}, \cite{ST94}. Applying  the Inductive
Algorithm~\ref{algAB}  we express  projective invariants in terms of affine, and affine
invariants in terms of Euclidean. We also obtain the relations  among the Euclidean, affine and projective arc-lengths and the corresponding invariant differential operators.

\begin{ex} Let us  use the moving frame for the special Euclidean group $SE(2,\R)$
acting on curves in $\R^2$  obtained in Example \ref{ex2} to build  a
moving frame for the special affine group.
We recall that the moving frame for $SE(2,\R)$  has been obtained on the first
jet space by choosing a cross-section $\{x=0,u=0,u_x=0\}$.
The special Euclidean group acts transitively on $J^1$
and the first invariant, Euclidean curvature $\kappa$ appears on the second order of
prolongation. The normalization of $u_{xxx}$ and $u_{xxxx}$ yields a third and fourth
 order invariants
$I_3^e=\kappa_s$ and
$I_4^e=\kappa_{ss}+3\kappa^3$.

The special affine transformation $SA(2,\R)$ on the plane is a semi-direct product
of the special linear group  $SL(2,\R)$ and translations in $\R^2$. We prolong it
to the first jet space of curves on the plane and notice that the isotropy  group $B$ of
the point  $z_0^{(1)}=\{x=0,u=0,u_x=0\}$ is given by linear transformations
$$\left(
\begin{array}{cc}
\tau  & \lambda  \\
& \frac 1\tau
\end{array}
\right) .
$$
Thus $SA(2,\R)=B\cdot SE(2,\R)$ and $B\cap  SE(2,\R)$ is finite. In fact
$B\cap SE(2,\R)=\{\Id,-\Id\}$.

We prolong the action of $B$ up to the fourth order:
 \begin{eqnarray*}
x&\rightarrow& \tau x+\lambda u;\\
u&\rightarrow& \frac 1 \tau u;\\
u_x&\rightarrow &\frac {u_x}{\tau(\tau+\lambda u_x)};\\
u_{xx} &\rightarrow& \frac {u_{xx}}{(\tau+\lambda u_x)^3};\\
u_{xxx}&\rightarrow &\frac {(\tau+\lambda u_x)
 u_{xxx}-3\lambda u_{xx}^2}{(\tau+\lambda u_x)^5};\\
u_{xxxx}& \rightarrow& \frac{(\tau+\lambda u_x)^2 u_{xxxx}
-10(\tau+\lambda u_x)\lambda u_{xx}u_{xxx}+15\lambda^2 u_{xx}^3}{(\tau+\lambda u_x) ^7}.
\end{eqnarray*}
 We restrict this transformation to the set
$\CS_E^4=\{z^{(4)}|\pi^4_1(z^{(4)})=z_0^{(1)}=(0,0,0)\}$ parameterized by
 Euclidean invariants $I_2^e,I_3^e$ and $I_4^e$, which, restricted to $\CS^4_E$  are
transformed under $B$ exactly by the same formulas as coordinate function $u_{xx},u_{xxx}$ and $ u_{xxxx}$.
Thus
\begin{eqnarray*}
I_2^e&\rightarrow& \frac {I^e_2}{\tau^3};\\
I_3^e&\rightarrow& \frac {\tau I^e_3-3\lambda (I_2^e)^2}{\tau^5};\\
I_4^e &\rightarrow & \frac{\tau^2 I^e_4
-10\tau\lambda I^e_{2}I^e_{3}+15\lambda^2 (I^e_{2})^3}{\tau ^7}.
\end{eqnarray*}
We emphasize  that the transformation formulas above are valid only on $\CS_E^4$
but not on  $J^4$.
We normalize the first transformation to one and the second transformation to zero.
This corresponds to choosing a cross-section
$$z_0^{(4)}=\{x=0,u=0,u_x=0,u_{xx}=1,u_{xxx}=0\}$$
 to the orbits of $SA(2,\R)$ on  $J^4$. Then
$$\tau =( I^e_2)^{1/3} \mbox{ and } \lambda =\frac{I^e_{3}}{3(I^e_{2})^{5/3}}.$$
We substitute this normalization in the transformation for $I^e_4$ to obtain the fourth order special affine invariant:
\[
I^a_{4}=\frac{I^e_{2}I^e_{4}-\frac 5 3(I^e_{3})^2}{(I^e_{2})^{8/3}},
\]
which we call the affine curvature and denote as $\mu$. We recall that $I_2^e=\kappa$, $I_3^e=\kappa_s$ and  $I_2^e=\kappa_{ss}+3\kappa^3$ and so the  affine curvature can be written in terms of the Euclidean curvature and its derivatives as follows:
\[
\mu= \frac{\kappa (\kappa _{ss}+3\kappa ^3)-\frac 5 3\kappa _s^2}{\kappa
^{8/3}}.
\]
 The moving frame for the special affine group is the product of the matrices:
\[
\left(\begin{array}{ccc}
\kappa^{1/3}  &\frac 1 3 \frac{\kappa_s}{\kappa^{5/3}} &0 \\
0& \frac 1 {\kappa^{1/3}} &0\\
0&0 &1
\end{array}
\right)
\left(
\begin{array}{ccc}
\frac 1{\sqrt{1+u_x^2}}&\frac {u_x}{\sqrt{1+u_x^2}} &-\frac {uu_x+x}{\sqrt{1+u_x^2}}\\
-\frac{u_x} {\sqrt{1+u_x^2}}&\frac {1}{\sqrt{1+u_x^2}} &\frac {xu_x-u}{\sqrt{1+u_x^2}}\\
0&0&1
\end{array}
\right)
\]
Using formula (\ref{ABcinvc}) one can obtain an affine contact invariant
  horizontal form $d\alpha$ in terms of the Euclidean arc-length
$ds$:
$$d\alpha=\sigma_B^*\,\sigma_E^*\,\pi_H\,w_E^*\,d_H\,w_B^*\,(x),$$
where the  Euclidean invariantization of  $d_H\,w_B^*\,(x)=(\tau+\lambda u_x)\,dx$
equals to
$\tau\, ds$ and hence
\begin{equation}\label{kappamu}
d\alpha=\sigma^*_B(\tau\,ds)=(I_2^e)^{1/3}ds=\kappa^{1/3}ds.
\end{equation}
The form  $d\alpha$ is called the affine arc-length. Written in the
standard coordinates  $d\alpha=u_{xx}^{1/3}dx.$
The relation (\ref{kappamu}) between the affine and the Euclidean arc-lengths
provide a natural explanation for the affine curve evolution equation in \cite{ST94}.
The relation between invariant differential operators follows immediately:
$$\frac d {d\alpha}=\frac 1 {\kappa^{1/3}}\frac d {ds},$$
which enables us to obtain all higher order affine invariants in terms of the Euclidean
ones.

\end{ex}

\begin{ex}
Let us now use the moving frame for the special affine  group to build a moving frame for the  projective group $PSL(3,\R)$ locally acting on the plane by the transformations:
 \begin{eqnarray*}
x&\mapsto& \frac{\alpha  x+\beta u +\gamma}{\delta x+\epsilon u+\zeta};\\
u&\mapsto& \frac {\lambda  x+ \nu  u+\tau}{\delta x+\epsilon u+\zeta}.
\end {eqnarray*}
where the determinant of the  corresponding matrix equals to one.
The affine moving frame has been found by choosing the cross-section
  $$z_0^{(3)}=\{x=0,u=0,u_1=0,u_{2}=1,u_{3}=0\}\in J^3.$$
 The isotropy group $B$
of $z_0^{(3)}$ for the prolonged action of  $PSL(3,\R)$ consists of the transformations:
$$\left(
\begin{array}{ccc}
1  & ab &0 \\
0&a& 0\\
b & c&\frac 1 a
\end{array}
\right) .
$$
Thus $PSL(3,\R)=B\cdot SA(2,\R)$ and $B\cap  SA(2,\R)$ is finite.
Let
$$\CS_A^7=\{z^{(7)}|\pi^7_3(z^{(7)})=(0,0,0,1,0)=z_0^{(3)}\}.$$
 The affine invariants
$I^a_4=\mu,\, I^a_5,\, I^a_6, I^a_7$ can serve as coordinate
functions on $\CS_A^7$ which are  transformed under the action of
$B$ by exactly the same rules as coordinate functions
$u_4,u_5,u_6,u_7$. We have computed the prolongation of the action
of $B$ to the seventh order using {\sc Maple} and found that:

\begin{eqnarray*}
I_4^a&\rightarrow& \frac {I^a_4-3a^2b^2+6ac}{a^2};\\
I_5^a&\rightarrow& \frac {I^a_5}{a^3};\\
I_6^a &\rightarrow & \frac{I^a_6+3abI_5^a+30I^a_4(2ac-a^2b^2)+180a^2c(c-ab^2)+45a^2b^2}
{a ^4};\\
I_7^a &\rightarrow & \frac{I^a_7+7ab I^a_6+I_5^a(105ac-42b^2a^2)-35(I_4^a)^2ab}{a^5}.
\end{eqnarray*}

We normalize the  transformation for $I^a_5$ to one and the  transformations of
$I^a_4$  and $I^a_6$ to zero.
This corresponds to the cross-section
$$z_0^{(7)}=\{x=0,u=0,u_1=0,u_{2}=1,u_{3}=0,u_4=0,u_5=1,u_6=0 \}$$
to the orbits of $PSL(3,\R)$ on
$J^7$. Then
\begin{eqnarray*}
a &=&( I^a_5)^{1/3},\\
b & =& \frac {5(I^a_4)^2-I^a_{6}}{3(I^a_{5})^{4/3}},\\
c &=& \frac {(I_6^a)^2-10I^a_6(I^a_4)^2-3I^a_4(I^a_5)^2+25(I^a_4)^4}
{18(I_5^a)^{7/3}}.
\end {eqnarray*}
We substitute this normalization in the transformation for $I^a_7$ to obtain the seventh
 order projective  invariant:
\[
I^p_{7}= \frac {6I^a_7I^a_5-7(I_6^a)^2+70(I_4^a)^2I^a_6-105I^a_4(I^a_5)^2
-175(I^a_4)^4}{6(I^a_5)^{8/3}}
\]
which we call the projective  curvature and denote as $\eta$.
Using the recursion algorithm
from \cite{FO99} we can express
the higher order affine invariants in terms of $\mu$ and
its derivatives with respect to {\it affine} arc-length $d\alpha=u_{xx}^{1/3}dx$:
\begin{eqnarray*}
I_4^a=\mu,&\quad I_5^a=\mu_{\as},\\
I_6^a=\mu_{\as\as}+5\mu^2,&\quad I_7^a=\mu_{\as\as\as}+17\mu\mu_\as.
\end {eqnarray*}
This leads to the formula:
\[
\eta=\frac{-7\mu_{\as\as}^2+6\mu_{\as}\mu_{\as\as\as}-3\mu\mu_{\as}^2}{6\mu_{\as}^{8/3}}.
\]
The moving frame for the projective group  is the product of the matrices:
\begin{eqnarray*}
&& \left(\begin{array}{ccc}
1  & -\frac 1 3 \frac{\mu_{\as\as}}{\mu_\as} &0 \\
0& \mu_\as^{1/3} &0\\
-\frac 1 3 \frac{\mu_{\as\as}}{\mu_\as^{4/3}}&\frac 1 {18}
\frac{\mu_{\as\as}^2-3\mu\mu_\as^2}{\mu_\as^{7/3}} &\frac 1 { \mu_\as^{1/3}}
\end{array}
\right)
\left(\begin{array}{ccc}
\kappa^{1/3}  & \frac 1 3 \frac{\kappa_s}{\kappa^{5/3}} &0 \\
0& \frac 1 {\kappa^{1/3}} &0\\
0&0 &1
\end{array}
\right)\\
&\times &\left(
\begin{array}{ccc}
\frac 1{\sqrt{1+u_x^2}}&\frac {u_x}{\sqrt{1+u_x^2}} &-\frac {uu_x+x}{\sqrt{1+u_x^2}}\\
-\frac{u_x} {\sqrt{1+u_x^2}}&\frac {1}{\sqrt{1+u_x^2}} &\frac {xu_x-u}{\sqrt{1+u_x^2}}\\
0&0&1
\end{array}
\right).
\end{eqnarray*}
We can express the projective arc-length (that is a horizontal form which is
 contact invariant with respect to the projective action) in terms of the affine arc-length
$d\alpha$. We first lift the coordinate function $x$ to $B\times J^\infty$ by
$w_B^*\,(x)=\frac {x+abu}{bx+cu+\frac 1 a }$. The affine invariantization of
 $d_H\,w_B^*\,(x)$ produce a horizontal form $ad\alpha$ on $B\times J^\infty$
which is contact invariant with respect to the affine action. The projective arc-length
equals to
$$d\varrho=\sigma^*_B\,a d\alpha=(I_5^a)^{1/3} d\alpha=\mu_\as^{1/3}d\alpha.$$
The relation between invariant derivatives
$\frac d {d\varrho}=\frac 1 {\mu_\as^{1/3}}\frac d {d\alpha}$ allows us to obtain all higher order projective invariants in terms of the affine ones.
\end{ex}

\chapter{Application to  Classical Invariant Theory.}

\label{clinvth}
One of the central problems of  classical invariant theory is the equivalence and symmetry of multivariable polynomials  under linear changes of variables. We concentrate on the
 polynomials over complex numbers, but  we will  indicate how to adapt the results to real polynomials.
The standard action of the general linear group on $\C^m$ induces a representation
on the ring of polynomials $\C[\x]$:
\begin{equation}\label{AF}
\bar{F}(A\cdot \x)=F(\x),
\end{equation}
where
$A\in GL(m,\C)$,\, $\x\in \C^m$  and  $F\in \C[\x]$ is a polynomial in $m$ variables.
The set of  homogeneous polynomials of degree $n$ is mapped to itself under
transformation (\ref{AF}). Thus we will always restrict our attention to  homogeneous
polynomials of a certain degree, which are called {\it forms} in the classical invariant theory literature.
Polynomials of degree $n$ in $m$ variables form a linear space of dimension
${m+n-1 \choose n}$ isomorphic to the $n$-th symmetric tensor product of $\C^m$. The
coefficients of polynomials can serve as  coordinates  on this space and formula (\ref{AF}) induces a linear action of $GL(m,\C)$ on  the coefficients.
One can try to
classify polynomials   by computing {\em invariants} that are certain functions $H(\dots,a_{{i_1},\dots,{i_m}},\dots)$ of the coefficients,  invariant under the induced action. It turns out however that this action is not regular: not only the dimensions of the orbits vary but  also  some   of the orbits  are not closed and hence the orbits can
not be distinguished  by
continuous invariants. An algebro-geometric approach  to the classification of
 such orbits can be found in \cite{Kraft87}, \cite{VP89}.
We also note that even for the fixed number of variables the number of
invariants increases  when one increases the degree of polynomials.

We can bypass these difficulties by considering the graph of an $m$-variable
 polynomial $u=F(\x)$
 as a submanifold of  $\C^{m+1}$ (or  $\R^{m+1}$). The first $m$ coordinates
represent  independent variables and they are transformed by the
the general linear group $GL(m,\C)$ in the standard way, the last coordinate is
 considered to be dependent variable and the action on it is trivial.
The described action is regular on the open subset where not  all of the first $m$
coordinates are zero. We classify the orbits by constructing the corresponding  signature manifolds, parameterized by a certain number of differential invariants.
One of the advantages of this approach is that the set of
differential invariants parameterizing the
signature manifold depends only on the number of variables, but not on the degree of
polynomials. We also note that differential invariants restricted to the graphs of polynomials are   functions
$H(\x,\,\dots,a_{{i_1},\dots,{i_m}},\dots)$ depending on both the coefficients and the
 variables that are  invariant under simultaneous action of the general linear group on
the variables and the coefficients of  polynomials. Polynomial or rational functions with such properties  were widely used in
the classical invariant theory \cite{GraceYoung03}, \cite{Gur64} and were called
{\em (absolute) covariants}, whence covariants that depend  on the
coefficients $a_{{i_1},\dots,{i_m}}$ only  were called  {\em (absolute) invariants}
The simplest example of a covariant is the polynomial $F(\x)$ itself
(see formula (\ref{AF})).  Rational absolute covariants can be obtained as the ratio
of {\em relative  covariants} that are not strictly invariant
but might be multiplied by a certain power of the determinant of matrix
$A\in GL(\C,n) $:
$$H(A\cdot \x,\,\dots,A\cdot a_{{i_1},\dots,{i_m}},\dots)=(\det A)^k \, H( \x,\,\dots,a_{{i_1},\dots,{i_m}},\dots ).$$
The exponent $k$ is called the {\em weight} of the covariant.
In the terminology of classical invariant theory the completeness of
the fundamental  set of polynomial invariants (covariants) means that any
other  invariant (covariant)
 can be expressed as a {\em polynomial} function of the fundamental ones.
The finiteness of the fundamental set for the actions of linear reductive groups
was proved by Hilbert in 1890. This crucial result became a turning point
from  classical
computational approach in the  invariant theory to the modern
algebraic geometry approach.

There are several classical methods to obtain  a complete set of fundamental
 covariants.
Two powerful
 methods, known as {\em omega process (or transvection)} and {\em symbolic method} were formulated in the second part of the last century by the German school of invariant
theory led by Aronhold, Clebesh and Gordan. Both of the methods
are based on application of  certain differential operators (see \cite{O99} for outline of these methods and more historical remarks).  Classical processes can be also used to obtain joint covariants of several forms under a simultaneous linear transformation.

Despite of the enormous amount of results  obtained by classical approaches
 many equivalence and symmetry problems remain unsolved  even for the case of
polynomials over the real or complex numbers. The formulation  of the problem as
 a  problem of equivalence of submanifolds, a novel approach introduced by
Olver \cite{O99}, produces new results even in the most studied case of
{\em binary forms}, or  homogeneous polynomials in two variables.
In Section~\ref{bf} the results for binary forms from my paper with Peter Olver
\cite{BO00} are reproduced and then in the next sections
 the same approach is extended  to {\em ternary  forms}, that is homogeneous polynomials
in three variables. We note that differential invariants in the case of polynomials can
be chosen to be rational functions in the variables and the coefficients and so the
signature  manifold construction reduces to the problem of eliminating parameters
from rational expressions. This problem  can be solved  using  Gr\"obner basis algorithms
 \cite{CLO96}, \cite{Sturmfels93} and \cite{Froberg97}.
Thus theoretically we can construct the
signature manifold for a polynomial of any degree in any number of variables.
In practice however we are confronted with the complexity of Gr\"obner basis
computation which in many cases exhaust available  computer resources.
This limitation  significantly affects  the practical implementation of the moving
frame method described below.

\section{Symmetries and Equivalence of  Polynomials.}
\label{symp}
We consider the action of the general linear group on the space of polynomials in $m$
 variables.
The standard action of the general linear group on $\C^m$ induces a representation
on the ring of polynomials $\C[\x]$:
\begin{equation}\label{AF2}
\bar{F}(A\cdot \x)=F(\x),
\end{equation}
where
$A\in GL(m,\C)$,\, $\x\in \C^m$  and  $F\in \C[\x]$. Equivalently:
\[
\bar{F}(\x)=F(A^{-1}\x)
\]
Since this action preserves the grading on  the ring of polynomials we  can restrict it
to the  homogeneous
polynomials of a certain degree, which are called {\it forms}
in the classical literature.

\begin{df}
A polynomial $F$ is said to be {\it equivalent} to a polynomial ${\bar F}$ if there
exists $ A\in GL(m,\C)$ such that $\bar{F}(\x)=F(A\x)$.
\end{df}
\begin{ex} The binary form
$5x^2-2xy+2y^2$ is equivalent to  $ x^2+y^2 $
  under the change of variables
$$x\mapsto x+y ;\quad y \mapsto y-2x.$$
\end{ex}

To each form we associate a unique inhomogeneous polynomial by the formula:
\[
f(p_1,\dots ,p_{m-1})=F(p_1,\dots ,p_{m-1},1).
\]
We call $\p=\{p_1,\dots,p_{m-1}\}$  {\em projective variables} and we  also refer
to $f(\p)$  as a form.
We can restore the homogeneous polynomial from its inhomogeneous version by the formula
\begin{equation}\label{Ff}
F(x_1,\dots ,x_m)=x_m^nf\left
(\frac{x_1}{x_m},\dots ,\frac{x_{m-1}}{x_m}\right).
\end{equation}
\begin{rem} The formula above shows that it is important to remember the degree $n$ of
the
homogeneous form if we want to restore it from its inhomogeneous version. For example
the quartic form $x^2y^2+y^4$ and the qudratic $x^2+y^2$ have the same inhomogeneous version:
$p^2+1$.
\end{rem}

Let $ A=\left(
\begin{array}{cc}
B & \mathbf{t} \\
\mathbf{s} & c
\end{array}
\right) \in GL(m,\C),
$
be a linear transformation, where $B$ is an $\left( m-1\right) \times \left( m-1\right) $ matrix,
 $\mathbf{s}^T,\  \mathbf{t}\in \C^{m-1}$,  and $c$ is a scalar.
Let $\bar F(\x)$ be the image of $F(\x)$ under transformation (\ref{AF}) and
let $f(\p)$ be the inhomogeneous version of $F(\x)$,  then the induced
transformation
of $f(\p)$ follows  from  formula (\ref{Ff}):
\begin{equation}\label{Af}
\left(  \mathbf{s\cdot p} +c\right) ^n\bar f\left(\frac{B \mathbf{p+t}}
{\mathbf{s\cdot p} +c}\right)=f(\p).
\end{equation}
We call the corresponding transformation of $f(\p)$ {\it projective}.
The form $F(\x)$ is equivalent to ${\bar F(\x)}$ under the transformation (\ref{AF}) if
and only if $f(\p)$
is equivalent to ${\bar f(\p)}$ under the corresponding projective transformation (\ref{Af}).

In the homogeneous version one can consider the graph $u=F(\x)$ of
the polynomial as a submanifold in $\C^m$ under the transformation
\begin{equation}\label{lx}
\x \mapsto A\cdot\x,\qquad u\mapsto u.
\end{equation}
For  the inhomogeneous version of the problem we consider the graph $u=f(\p)$ in $\C^m$ under the transformation:
\begin{equation}\label{Aup}
{\p} \mapsto \frac{B \mathbf{p+t}}{\mathbf{s\cdot p} +c},\qquad u\mapsto
\left(  \mathbf{s\cdot p} +c\right) ^{-n} u.
\end {equation}
The transformation of  projective variables $\p$ in (\ref{Aup}) is {\it
linear fractional}:
\begin{equation}\label{lfp}
A\cdot\p= \frac{B \mathbf{p+t}}{\mathbf{s\cdot p} +c}
\end{equation}
Two matrices which are scalar multiples of each other, ${\tilde A} = \lambda A$,
induce
the same linear fractional transformation, and so
(\ref{lfp}) defines an action of the
projective
group $PSL(m,\C) = GL(m,\C)/\{\lambda \, \Id\}$ on $\C^{m-1}$.
Let $\pi \colon GL(m,\C) \mapsto  PSL(m,\C)$
 denote the standard
projection.

\begin{df}\label{prsymm}  The {\it symmetry group} of  $F$ is the subgroup
$G_F\subset GL(m,\C)$ consisting\ of all linear transformations   that map $F$ to itself.
 It coincides with the group  $G_f\subset GL(m,\C)$ which maps inhomogeneous version $f$ of $F$ to itself under transformation (\ref{Af}).
The {\it projective symmetry group} of $f$ is the subgroup
$\Gamma_f =\pi (G) \subset PSL(m,\C)$ consisting of all
linear fractional transformations  of $\p$ that give rise to symmetries of $f$.
In the real case $G_F\subset GL(m,\R)$ and $\Gamma_f\subset PSL(m,\R)$.
\end{df}

\begin{ex}
The form  $F(x,y)= x^2+y^2 $ is symmetric under any orthogonal map:
   \[(x,y)\mapsto\left\{\begin{array}{c} (\cos(\alpha)x+\sin(\alpha)y, \,
-\sin(\alpha)x+\cos(\alpha)y) \\
 (-x,\, -y )
\end{array}\right.\]
The inhomogeneous version of $F$ is $f(p)=p^2+1$ and the corresponding projective group of symmetries $\Gamma_f$ consists of linear fractional transformations
\[p\mapsto \frac{\cos(\alpha)p+\sin(\alpha)}{-\sin(\alpha)p+\cos(\alpha)}\]
\end{ex}
We notice that each projective symmetry in the preceding example corresponds to
two genuine symmetries of $F(x,y)$.
In general if the  form $F(\x)$ has   degree $n$ then $F(\lambda \x)=\lambda^n F(\x)$ and so if $\omega$ is any root of unity, $\omega^n=1$, then the diagonal matrix $\omega\Id\in G_F$, on the other hand $\pi(\omega \Id)=e\in\Gamma_f$ and so each projective symmetry                         gives rise to $n$ genuine symmetries  of the form $f(\p)$.
\begin{prop}
A transformation $A\in GL(m,\C)$ maps a form $F$ to  some scalar multiple of itself, say
$\mu F$ is and only if  $\pi(A)\in\Gamma_f$.
\end{prop}
{\em Proof.}
By substitution  $\bar F=\mu F$ into (\ref{AF2}) one obtains:
$$F(\x)=\mu F(A\cdot \x)=F( \root n\of\mu\, A\cdot \x).$$
and thus $\hat A =\lambda A$, where $\lambda=\root n\of\mu\, A$ belongs to $G_F$, and so by definition  $\pi(A)=\pi(\hat A)\in\Gamma_f$.
\qed
\begin{rem}
If $F$ is mapped to $\mu F$ by $A$ then its inhomogeneous version $f$ is mapped
 to $\mu f$, that is:
\begin{equation}
f(\p) = \left(\mathbf{s\cdot p} +c\right) ^{n}\,\mu\, f\left(\frac{B \mathbf{p+t}}
{\mathbf{s\cdot p} +c}\right).
\end{equation}
\end{rem}
\begin{rem} The original transformation rules (\ref{AF}), (\ref{Af}) apply to  forms of
 weight zero.  One can, more generally, consider forms of nonzero
weight $k$, with transformation rules
\begin{equation}
F(\x) = (\det A)^k \, \bar F(A\cdot \x),\quad
f(\p) = (\det A)^k \left(\mathbf{s\cdot p} +c\right) ^{n} \bar f\left(\frac{B \mathbf{p+t}}
{\mathbf{s\cdot p} +c}\right).
\end{equation}
If $n + mk\neq 0$, then the projective symmetry group of a weight $k$  form is the
same as that of its weight $0$ counterpart. However, the full symmetry group does not
have the same cardinality, and so are not isomorphic.
 Indeed, let   $A\in GL(m,\C)$ be any matrix whose associated linear fractional
transformation  belongs to the
projective symmetry group of $F$ of  weight $0$, and let $\det A=\Delta$.
Then $A$ maps $F$ to a
scalar multiple
of itself, say $\mu \,F$ and so $F(\x)=\mu\, F(A\cdot\x)$.
Consequently, the scalar multiple
\begin{equation}\label{scal}
\hat A = \lambda A ,\quad \mbox{where}\quad \lambda^{n+mk}\Delta^k  =\mu ,
\end{equation}
 is a symmetry of the weight $k$ form $F$.
Therefore, when $n+mk\neq 0$, each projective symmetry gives rise to $n+mk$
matrix symmetries.

In the  exceptional case when $n + mk= 0$,  if $A\in GL(m,\C)$
is any symmetry, so is {\em any} scalar multiple $\lambda A$.
Thus  each projective symmetry gives rise to a one-parameter family of symmetries in
$GL(m,\C)$.
On the other hand, the projective group of symmetries in this  case is smaller in
general than
for the other weights. Indeed, given a projective symmetry for a  weight zero form,  one can
always find a matrix representative $A\in SL(m,\C)$. Due to (\ref{scal}), this
representative is a
symmetry for the exceptional weight  $k=-\frac n m$ if and only if the corresponding
$\mu=1$.
On the other hand, if  $A\in GL(m,\C)$ is a symmetry for weight $k=-\frac n m$ then       a unimodular matrix    $\Delta^{-\frac 1 m}A$ is a
symmetry for the
corresponding zero weight form.

\begin{ex}\label{nonzw}
Let $f(p)=p^4+3p^2+1$ with $n=4,\, m=2$ correspond to a binary form
$f(x,y)=x^4+3x^2y^2+y^4$.
 Then the projective symmetry group in the case   of weight zero
(as well as for  any other general weight) consists of the transformations mapping
$p$ to  $p,-p,\frac 1
p,-\frac 1 p$. Any $GL(2,\C)$ representative of these projective maps is a symmetry
for the exceptional
weight $k = -2$.

On the other hand,  for the form $f(p)=p^4+1$ of with $n=4,m=2$
 the projective symmetries
for a general weight are
$p,-p,\frac 1 p,-\frac 1 p, ip,-ip,\frac i p,-\frac i p$.
In the exceptional case, however, only the first four, namely $p,-p,\frac 1 p,-\frac 1 p$, are  symmetries.
\end{ex}
\end{rem}


\begin{df} A homogeneous  form is called {\it nonsingular} if its symmetry group
$G_f$ is
finite. The {\it index} of a nonsingular form $f(\p)$ is the cardinality
$\#G_f$ of its symmetry group. The {\it projective index} of $f(\p)$ is the cardinality $\#\Gamma_f$ of
its projective symmetry group $\Gamma_f =\pi (G_f)$.
\label{nonsingular}
\end{df}
Thus, for nonsingular forms, the 
indices are simply related by

\begin{equation}\label{indices}
\#G = l\cdot \#\Gamma ,\, \mbox{where}\quad
  l = \left\{\begin{array}{lcl}n&\,&\mbox {for complex forms of degree $n$,}\\
2&\,&\mbox {for real forms of even degree,}\\
1&\,&\mbox{for real forms of odd degree.}\end{array}
\right.
\end{equation}

In what follows we  address the problem of classification  of polynomials under linear transformation in its inhomogeneous version (\ref{Af}).  We  reduce  the problem to  the problem of equivalence of  the graphs of polynomials  under the transformation (\ref{lfp}) so we can make a full use of  the results described in the previous two chapters.

\begin{rem}\label{svariety}
As it has been mentioned in the introduction to this chapter the signature manifold of
the form $f(\p)$ can be parameterized by a set of rational differential invariants
restricted to $f$. Elimination of the variables $\p$ produces polynomial relations among invariants  which define the smallest  variety containing
the signature manifold of $f$. We call this variety
{\em the signature variety} and notice  that it is irreducible in both the real and the complex case (Proposition~6, ch.~4, \S~6 in \cite{CLO96}).
The relations among invariants, which we obtain for a  polynomial with real
coefficients, are  the same  whether the real or complex equivalence problem is
considered. The real
classification however includes more equivalence classes!
It can be explained by the fact that the signature manifold,
defined by parametric equations, does not necessarily fill up entire signature variety
(see \cite{CLO96}) neither in complex nor in real case. In the real case
 two different parts of the signature variety may correspond to two different
signature manifolds of the same dimension (see Example~8.69 in \cite{O99}).
In the complex case, however, different signature manifolds are included in different
 signature varieties.
Indeed assume $\CC(f)$ and $\CC(\tf)$ are two signature manifolds and $V$ is the smallest signature variety containing them.
Assume that   $\CC(f) \cap \CC(\tf)\neq \emptyset$ then $f$ and $\tf$ are locally equivalent by Theorem~\ref{equiv}  and so they are globally equivalent (for any
analytic function  local equivalence implies global). Otherwise, if
$\CC(f) \cap \CC(\tf)= \emptyset$
then  in the complex case, it follows from Theorem~3, ch.~3, \S~2 \cite{CLO96}  that
there is a subvariety $W\subsetneq V$ such that $V-W\subset \CC(f)$, and so
$\CC(\tf)\subset W$. This contradicts to
the assertion that $V$ is the minimal variety  containing $\CC(\tf)$.
Assume that $W$ contains  the signature variety of $\CC(\tf)$, then
$\dim\, \CC(\tf)<\dim\, \CC(f)$ since
$V$ is irreducible and $W\subsetneq V$.

We conclude that  the polynomial relations among invariants provide a solution
for  the problem of equivalence over the complex numbers. In the real case,
these relations produce necessary but not sufficient conditions of  equivalence, and
so more  detailed analysis is required to complete the classification over  reals.
\end{rem}

The symmetry groups  of equivalent polynomials are related by matrix conjugation:
\[\bar{F}(\mathbf{x)=}F(A\mathbf{x}) \Longrightarrow G_{\bar{F}}=
 AG_{{F}}A^{-1}.\]
Thus the problem of the classification of polynomials is closely related to the problem of the classification of their   symmetry groups up to  matrix conjugation.
Theorems~\ref{symm},~\ref{sfinite} of Chapter~\ref{prelim} provide the foundation to an  algorithm
which determines  the dimension of the symmetry group of a given polynomials and, in
the case when the  cardinality of the symmetry group is finite,
 explicitly computes all transformations that belong to it.  We start with the simplest case  of binary forms.

\section {Binary forms.}
\label{bf}
The general linear group
$$
GL(2,\C) = \left\{\ A =
\left.\left(\begin{array}{cc} \alpha& \beta \\ \gamma &\delta \end{array}\right)\right|
\alpha \delta - \beta \gamma \ne 0\ \right\}$$
acts on two-dimensional space by invertible linear transformations
\begin{equation}\label{GL2t}
{\bar{x}  = \alpha  x + \beta y, \qquad \bar{y}  = \gamma x + \delta y,}
\end{equation}
and thereby induces an irreducible linear representation on the space of binary forms
$$F(x, y) = \sum_{i = 0}^n \; a_i \,x^i y^{n-i} $$
of the fixed degree $n$.
This corresponds to a linear fractional transformation
\begin{equation}\label{lft}
\bar p = \frac{\alpha p + \beta}{\gamma p + \delta }\,.
\end{equation}
on the projective coordinate  $p= \frac x y$.
 The induced transformation rule for inhomogeneous polynomials of degree $n$ is:
\begin{equation}\label{fbarp}
f(p) = (\gamma  p + \delta )^n\,\bar f (\bar p) = (\gamma  p + \delta )^n\,
\bar f \left(\frac{\alpha p + \beta }{ \gamma p + \delta }\right)
\end{equation}

We reformulate   the symmetry and equivalence problem for polynomials as the symmetry and equivalence problem for
the graph  of a  polynomial $u=f(p)$ considered as a submanifold in  $\C^2$. The
 transformation  rules
\begin{eqnarray*}
p& \mapsto&  P=\frac{\alpha p + \beta }{ \gamma p + \delta },\\
u & \mapsto&v=(\gamma  p + \delta )^{-n}u
\end{eqnarray*}
 for coordinates $p$ and $u$ in $\C^2$ can be prolonged to the $k$-th order jet space:

\begin{eqnarray*}
u_p=u_1 &\mapsto& v_1=\frac 1 {\Delta \sigma^{n-1}}(-n\gamma u+\sigma u_p),\\
u_{pp}=u_2&\mapsto&  v_2=\frac 1 {\Delta^2 \sigma^{n-2}}(n(n-1)\gamma^2u-2(n-1)
\gamma\sigma u_1+\sigma^2u_2)\\
 u_k  &\mapsto&v_k=\frac {1}{\Delta^{k}\sigma^{n-k}}\sum_{j=0}^k (-1)^{k-j} {k\choose j}
\frac {(n-j)!}{(n-k)!}\gamma ^{k-j} \sigma^j u_j ,
\end{eqnarray*}
where $\Delta=\alpha \delta -\beta \gamma$ and $\sigma=\gamma p+\delta$

A moving frame can be defined on the second order by choosing a cross-section
$p=0,u=1,u_1=0,u_2=\frac 1 {n(n-1)}$. By solving the equations:
$$P=0,\, v=1,\, v_1=0,\, v_2=\frac 1 {n(n-1)}$$
one obtains a moving frame:
$$\left(
\begin{array}{cc}
u^{\frac{1-n}n}\sqrt{H} & -pu^{\frac{1-n}n}\sqrt{H} \\
\frac{u^{\frac 1n}u_p}{nu} & u^{\frac 1n}\left( 1-\frac{u_pp}{nu}\right)
\end{array}
\right),
$$
where $H(p)=n(n-1)(uu_{pp}-\frac{n-1}nu_p^2)$ is the inhomogeneous version of the Hessian
$H(x,y)=F_{xx}F_{yy}-F_{xy}^2$ of a form $F(x,y)$.
By substitution of the moving frame into the formulas for  higher order lifted
invariants $v_k$ we obtain the general formula for invariants:
\[
I_k=\frac 1{\sqrt{H}^k} \sum_{j=0}^k(-1)^{k-j}\left(
\begin{array}{c}
k \\
j
\end{array}
\right) \frac{(n-j)!}{(n-k)!}\left( \frac{u_p}n\right) ^{k-j}u^{j-1}u_{^j}.
\]

Due to theorems \ref{equiv}, \ref{symm} from Chapter 1
the symmetry and equivalence properties
of a binary form are entirely determined by just these two differential invariants:
\begin{eqnarray*}\label{I3I4}
I_3&=&\frac 1{\sqrt{H}^3}(u^2u_{ppp}-3\frac{n-2} n uu_pu_{pp}+2\frac{(n-1)(n-2)}{n^2}
u_p^3),\\
I_4&=&\frac 1{H^2}(u^3u_{pppp}-4\frac{n-3} n u^2 u_pu_{ppp}+
6\frac{(n-2)(n-3)}{n^2}uu_p^2u_{pp}\\
&-&3\frac{(n-1)(n-2)(n-3)}{n^3}u_p^4).
\end{eqnarray*}

This differential invariants can be re-expressed in terms of the classical covariants, which can be obtained by omega process:
\begin{equation}\label{HTU}
H(F)=F_{xx}F_{yy}-F_{xy}^2,\quad T(F)=F_xH_y-F_yH_x\quad \mbox{and}\quad U(F)=F_xT_y-F_yT_x.
\end{equation}
Covariants $H$, $T$ and $U$ has weight 2, 3 and 4 respectively and thus the rational functions:
\begin{equation}\label{JK}
J=\frac{T^2}{H^3},\qquad\qquad K=\frac{U}{H^2}
\end{equation}
are absolute rational covariants.
By writing down the inhomogenization of $H$, $T$ and $U$:
\begin{eqnarray*}
H(f)&=& n(n-1) \left [ f f'' - \frac{n-1}{ n}\,(f')^2 \right],\\
T(f)&=&-n^2(n-1) \left [f^2 f''' - 3\,\frac{(n-2)}n \>f f' f'' +
2 \,\frac{(n-1)(n-2)}{n^2} \>(f')^3\right ],\\
U(f)&=&n^3 (n-1) \left[ f^3 f'''' - 4 \,\frac{(n-3)} n \, f^2 f'f''' +
6 \, \frac{(n-2)(n-3)}{ n^2}\> f f'\,^2f''\right.\\
&-&\left.3\,\frac{(n-1)(n-2)(n-3)} {n^3}\, (f')^4
\right]
 - 3\,\frac{(n-2)}{( n-1)} \,H^2.
\end{eqnarray*}
Comparing this formulas with (\ref{I3I4}) we conclude that:
\begin{equation}\label{JandK}
J=-n^4(n-1)^2I_3^2,\qquad\qquad K=n^3(n-1)I_4-3\,\frac{(n-2)}{( n-1)}.
\end{equation}
Thus the signature manifold for  $f(p)$ parameterized by $J$ and $K$ can be used equally well to solve
 the equivalence problem as the one parameterized by $I_3$ and $I_4$ and we will formulate our results in terms of the former, more traditional covariants.
We note that none of the invariants are defined when $H\equiv 0$. However  this happens if and only if
$F(x,y) = (c x + d y)^n$ is the $n$-th power of a linear form and thus is equivalent to polynomial of one variable only. As it pointed out in chapter 10 of \cite{O99}
the generalization of this statement is true when the number of variables $m\leq 4$.
The Hessian of a form of four and less variables is equal to zero if and only if
the form  is equivalent under a linear change of variables to a form of less number
of variables. When the number of variables is greater than four the `if' part of
the statement is still true, but the `only if' part fails \cite{GordanNother}
 for instance when
$$F(x_1,\dots,x_5)=x_1^2x_3+x_1x_2x_4+x_2^2x_5.$$

We can now state the classification theorem for binary forms:
\begin{thm}
\label{fsym}
Let  $f(p)\not\equiv 0$  correspond to  a nonzero binary form of
degree  $n$. The {\it symmetry group}
of $f(p)$ is:
\begin{description}
\item[a)] A two-parameter group if and only if   $H\equiv  0$ if and only if $f(p)$ is
equivalent to a constant.\item[b)] A one-parameter group if and only if $H\not\equiv 0$ and $T^2$
is a constant multiple
of $H^3$ if and only if $f(p)$ is complex-equivalent  to a monomial $p^k$, with $k\ne 0,n$.
 In this case:
$$ J = -\frac{4(n-2k)^2}{k(n-k)(n-1)},\qquad  K =-\frac{6(n-2k)^2}{k(n-k)(n-1)}
$$

\item[c)] A finite group in all other cases.
\end{description}
\end{thm}

\begin{rem}  A real binary form is complex-equivalent to a monomial if and only if it is
real-equivalent to either a real monomial $\pm p^k$ or to the form $\pm (p^2 + 1)^m$, the
latter only occurring in the case of even degree $n=2m$.
\end{rem}
Therefore, a binary form is nonsingular if and only if its
rational covariant $J$ is {\it not} constant if and only if the form is
not complex-equivalent
to a monomial.  The next result is fundamental  for our algorithm for determining the
(finite) symmetry group of a  nonsingular binary form.

\begin{thm} Let $f(p)$ correspond to  a nonsingular complex binary form.
Then $P = \varphi (p)$ is a
complex analytic solution to the {\it rational symmetry equations}
\begin{equation}\label{JKsym}
J (P) =  J(p), \qquad  K (P) =   K(p)
\end{equation}
if and only if
$P = (\alpha p+\beta )/(\gamma p + \delta )$ is a linear fractional
transformation belonging to the projective symmetry group of $f(\p)$.
\label{S}
\end{thm}
Thus {\it all} the solutions to the symmetry equations (\ref{JKsym}) are
necessarily
linear fractional transformations!  As remarked above, given a projective
symmetry, the corresponding symmetry  matrix
$A\in GL(2,\C)$ is uniquely determined  up to multiplication by an
$n$-th root of unity.  Since the linear fractional transformation only determines $A$ up to a scalar
multiple, one must substitute into the transformation rule \ref{fbarp} for the form to
unambiguously specify the symmetry matrix.

In the real case, if the degree of $F$ is odd, $n = 2m+1$, then the basic symmetry
\ref{S} holds
as stated.  Moreover, each real linear fractional solution to the symmetry equations
(\ref{JKsym}) corresponds to a unique matrix symmetry.
On the other hand, if the degree of
$F$ is even, $n = 2m$, then the sign of $F$ is invariant, and a real solution to the
symmetry
equations (\ref{JKsym}) will induce a real projective symmetry, and thereby two
real matrix symmetries of
the form if and only if it preserves the sign of $F$.
\
The explicit computations for of the symmetry group of a nonsingular complex binary form
 relies
on Theorem \ref{S}, and hence requires  solving the fundamental symmetry equations
(\ref{JKsym}) which can be rewritten as polynomial equations in $p$:
\begin{equation}\label{EqJK}
A(p)B(P) =  A(P)B(p), \qquad  C(P)D(p) =  C(p)D(P),
\end{equation}
where:
$$J=\frac{T^2}{H^3}=\frac A B ,\qquad\qquad K=\frac{U}{H^2}=\frac C D .$$
Polynomials  $A$ and $B$ have no common factors, nor do $C$ and $D$.
 Bounds on the index or number of symmetries of a binary form can be determined without explicitly
solving the bivariate symmetry equations (\ref{EqJK}).
The fact that $f(p)$ is not equivalent to a monomial implies that $T^2 $
is a not a constant
multiple of $H^3$, and hence the first equation in (\ref{EqJK}) is  nontrivial.  Therefore, the
projective index of $f(\p)$ is always bounded by the degree of the first equation in  $p$,
which in turn is bounded by
$6n-12$ with equality if and only if $T$ and $H$ have no common factors.
The second bivariate
polynomial  is trivial  if and only if the covariant
$U$ is a constant multiple of $H^2$.
Forms for which $U=cH^2$ will be distinguished as belonging to the
{\em maximal discrete symmetry class}.  Indeed, if $T$ and $H$ have no common
factors and all the roots of the first equation  are simple,
then the projective index of such a form takes
its maximum possible value, namely $6n-12$.
On the other hand, if
$U$ is not a constant multiple of $H^2$, then the projective
index is bounded by the degree of the second polynomial  which is at most $4n-8$.

\begin{thm}\label{indexn} Let $k$ denote the projective index of a binary form $Q$ of degree
$n$ which is not complex-equivalent to a monomial.  Then
$$
k \leq  \left\{ \begin{array}{lcl}
6n-12&\qquad&\mbox{ if }\, U = c H^2\, \mbox{for some constant }\, $c$,\,\mbox {or}\\
4n-8&\qquad&\mbox{ in all other cases.}
\end{array}\right.
$$
\end{thm}
\smallskip\noindent
The real case clearly admits the same bounds on the projective index, since one must
determine the number of common {\em real} solutions to (\ref{EqJK}), and, in the
case of even degree,
whether the sign of $Q$ is the same at each solution.
Consequently, the index of a binary form of degree $n$ is bounded
by either $(6n - 12)l$ or $(4n - 8)l$, where $l=n$ in the complex case, $l=2$ in the case of real
forms of even degree and $l=1$ for real forms of odd degree.

Since the symmetry groups of
equivalent polynomials are related by matrix conjugation in $\GLC2$, a complete list of
possible projective symmetry groups is provided by the following theorem, as presented
in Blichfeldt, (\cite{Bl17} p.~69).
\def\An{{\cal A}_n}
\def\Dn{{\cal D}_n}
\def\T{{\cal T}}
\def\O{{\cal O}}
\def\I{{\cal I}}
\begin{thm}
Up to matrix conjugation there are five different types of finite
subgroups of the projective group $PSL(2,\C)$:
\begin{description}
\item[a)]
The $n$ element abelian group $\An$ is
generated by the transformation $p\mapsto \omega p$,
where $\omega $ is a primitive $n$-th root of unity.

\item[b)]
The $2n$ element {\em dihedral group} $\Dn$ is the group
obtained from $\An$ by adjoining the
transformation $ p\mapsto 1/p$.

\item[c)]
The $12$ element {\em tetrahedral group} $\T$ is the primitive group
generated by the transformations
\begin{equation}\label{st}
\sigma \,\colon\ p\longmapsto -p, \quad \qquad  \tau \,\colon\ p\longmapsto
\frac{i(p+1)}{p-1},
\end{equation}
of respective orders $2$ and $3$.

\item[d])
The $24$ element {\em octahedral group} $\O$ is the primitive group generated by the
transformation $\tau $ in \ref{st} along with
\begin{equation}\label {iota}
\iota\,\colon\ p\longmapsto ip
\end{equation}
of order $4$. Note that $\iota^2=\sigma $, and so $\T\subset \O$.

\item[e)] The $60$ element {\em icosahedral group} $\I$ is the primitive group generated
by the transformations $\sigma ,\tau $ given above, along with the transformation
\begin{equation}\label{rho}
\rho \,\colon\ p\longmapsto\
\frac{2p-(1-\sqrt{5})i-(1+\sqrt{5})}{\left[ (1-\sqrt{5})i-(1+ \sqrt{5})\right] p-2}
\end{equation}
of order $2$.
The tetrahedral group  is also a subgroup of the icosahedral group: $\T\subset \I$.
\end{description}
\label{finitesubgroup}
\end{thm}

Since the maximal number of elements in the projective symmetry group of a form of
degree $n$ is bounded by $6n-12,$ then the tetrahedral group can appear as a symmetry
group only when $n\geq 4$, the octahedral group  is a possible symmetry group only if
$n\geq 6$
and the icosahedral group  is possible only if $n\geq 12$.

\def\pf{K}
\def\pfb{\bar \pf}
\def\pg{L}
\def\pgt{\widetilde \pg}

We can describe the invariants of the three primitive groups using
the following polynomials:
\begin{equation}
\label{phipsi}
\begin{array}{ll}
\pf _4=x^4-2\sqrt{3}\,i\,x^2y^2+y^4,& \pfb_4=x^4+2\sqrt{3}\,i\,x^2y^2+y^4,\\
\pf _6=x^5y-xy^5&
\pf _8=x^8+14x^4y^4+y^8=\pf _4\pfb_4,\\
\pf _{12}=x^{12}-33\left( x^8y^4+y^8x^4\right) +y^{12},\\
\pg _{12}=22\sqrt{5}\pf _6^2+5\pf _{12},&
\pgt_{12}=-22\sqrt{5}\pf _6^2+5\pf _{12},\\
\pg _{20}=3\pf _8\pf _{12}-38\sqrt{5}\pf _6^2\pf _8,&
\pgt_{20}=3\pf _8\pf _{12}+38\sqrt{5}\pf _6^2\pf _8,\\
\pg _{30}=6696\pf _6^5+225\pf _6\pf _8^3-580\sqrt{5}\pf _6^3\pf _{12},&
\pgt_{30}=6696\pf _6^5+225\pf _6\pf _8^3\\
 &+580\sqrt{5}\pf_6^3\pf _{12}.\hskip-20pt
\end{array}
\end{equation}
Huffman, (\cite{Huffman} Theorem 4.1), provides the complete characterization of polynomials
whose symmetry groups contain one of these primitive groups.

\begin{prop}\label{phiinvs}
The symmetry group of a binary form  $F$ contains:
\begin{description}
\item [a)] An icosahedral group if and only if  it
is equivalent to a polynomial of the one of the two forms
$\Phi (\pg _{12},\pg_{20}) + \pg _{30}\Psi (\pg _{12},\pg _{20})$
or $\Phi (\pgt_{12},\pgt_{20}) + \pgt_{30}\Psi (\pgt_{12},\pgt_{20})$.

\item[c)] An octahedral group if and only if  it is equivalent to
a polynomial of the one of the two forms
$\Phi (\pf _6,\pf _8)$ or $\pf _{12}\Phi (\pf _6,\pf _8)$.\\
\item[d)]
A tetrahedral group if and only if it is equivalent to
a polynomial from the following list:
\begin{eqnarray*}
\Phi (\pf _6,\pf_8) + \pf _{12}\Phi (\pf _6,\pf _8),&
\Phi (\pf_4,\pf _6),& \Phi (\pfb_4,\pf _6),\\
\pf _4\Phi (\pf _6,\pf_8) + \pf _4^2\,\Phi (\pf _6,\pf _8),&
\pf _4\Phi (\pfb_4,\pf _6),& \pfb_4\Phi (\pf_4,\pf _6),\\
\pfb_4\Phi (\pf_6,\pf _8) +\pfb_4^2\,\Phi (\pf _6,\pf _8),&
\pf _4^2\,\Phi (\pfb_4,\pf _6),& \pfb_4^2\,\Phi (\pf _4,\pf _6) .
\end{eqnarray*}
\end{description}

Note in particular that only forms of even degree can admit a primitive symmetry group.
\end{prop}
{\sc Maple} code was written to explicitly compute the
symmetries of binary forms.   Details of the programs and some of the difficulties we experienced in the
implementation are discussed in the appendix~A. The program {\tt symm}   computes the fundamental invariants
$J$ and $K$, determines the dimension of the symmetry group, and, in the case of a finite
symmetry group, solves the two equations \ref{EqJK} to find explicit form of
the projective symmetries. The actual matrix symmetries are then computed by the program {\tt matrices}
by substituting the linear fractional transformations in the projective symmetry group into the form
in order to determine the appropriate scalar multiple.   We now present some typical examples
resulting from our computations.

\begin{ex}\label{cubics}
{\em Cubic forms}. All binary cubics with discrete symmetries are
equivalent to $x^3+y^3$, or, in inhomogeneous form, to $p^3+1$. Therefore, the symmetry group
of a nonsingular cubic is isomorphic to the symmetry group
of $p^3+1$.  Applying our algorithm, we find a complete solution to the symmetry equations
\ref{JKsym} is the projective symmetry group $\Gamma $ given by the six linear fractional transformations
taking
$p$ to
$$
p,\quad \frac 1p,\quad\omega p,\quad \omega ^2p,\quad \frac \omega p,\quad
\frac{\omega ^2}p,
$$
where $\omega =-\frac 12+i\frac{\sqrt{3}}2$ is the primitive cube root of unity.
Since the covariants of any cubic form satisfy the syzygy  $U = -\frac 32 H^2$, all non-degenerate cubics
have maximal discrete symmetry groups of projective index $6 $, which equals the number of different
permutations of the three roots.
The full matrix symmetry group $G$ of this cubic has $18$ elements, since we can also multiply by
a cube root of unity, and is generated by the three matrices
$$\left(\begin{array}{cc} \omega &0\\ 0& \omega \end{array}\right),\qquad
\left(\begin{array}{cc} 0 &1\\ 1& 0 \end{array}\right),\qquad
\left(\begin{array}{cc} \omega &0\\ 0& \omega^2 \end{array}\right).
$$
In this case,  $G \simeq \Gamma \times \Z_3$ is a Cartesian product group.
  In the real case, one requires real solutions
to (\ref{JKsym}, and hence $f$ has (projective) index $6$ if its discriminant $\Delta <0$, but (projective)
index $2$ if $\Delta >0$.

The \Maple\ code can be used to compute the explicit symmetries of other cubics.  For example,
the cubic $f(p) = p^3+p$ leads to the following six element
group of linear fractional transformations
$$
p,\qquad -p,\qquad \frac{ip+1}{3p+i},\qquad \frac{ip-1}{-3p+i}\>,\qquad
\frac{-ip+1}{-3p+i}\>,\qquad \frac{-ip+1}{3p+i}\>.$$
The matrix generators of the symmetry group are
$$\matrixb \omega 00{\omega },\qquad \matrixb 100{-1},\qquad
\frac 12 \matrixb 1{-i}{-3i}1.
$$
The second and third matrices correspond, respectively, to the second and third linear fractional
transformations.  Note that one must, in accordance with the general procedure, rescale the
matrices as required by the condition that $f$ must be mapped to itself.
Difficulties arise when \Maple\ gives the solutions of equations \ref{JKsym} not as
 rational functions, but involving roots of
polynomials. An example is the cubic $f(p) = p^3+p+1$, which is discussed in Appendix A.
\end{ex}

\begin{ex}\label{Ex4} {\em Quartic forms}. A polynomial of degree $4$ has a finite symmetry group
if it is equivalent to either
$$
p^4+\mu\, p^2+1, \quad \mbox{or}\quad  p^2+1,
$$
where $\mu \neq \pm 2$.  The former has all simple roots; the latter has a double root at
$\infty $.

In the first situation, the symmetry group will depend on the value of $\mu $.
For general $\mu $, the projective symmetry group is a dihedral group ${\cal D}_2$, generated by $-p$ and
$1/p$.  When $\mu =0$ it becomes a dihedral group ${\cal D}_4$, generated by $ip$ and
$1/p$.  The associated matrices are the obvious ones, namely $\matrixb{-1}001$, $\matrixb 0110$ in the
first case, and $\matrixb{i}001$, $\matrixb 0110$ in the second.

The cases $\mu =\pm 2i\sqrt{3}$ corresponds to the polynomials $\pf _4$ and $\pfb_4$
listed in \ref{phipsi} above, and so the projective symmetry group is the $12$
element octahedral group $\O$. This case has the maximal size discrete symmetry  group. The linear
fractional transformations are generated by
$-p$ and
${i\left( p-1\right) }/{(p+1)}$.  These correspond to different
matrices in each case:
$$\begin{array}{lclc}\pf _4:&
\matrixb {-1}001,&
\frac 1{\left( 2-2i\sqrt{3}\right) ^{1/4}}\
\matrixb i{-i}11,&
\qquad\mbox{when}\qquad  \mu = 2i\sqrt 3,\\
\pfb_4:&
\matrixb {-1}001,&
\frac 1{\left( 2+2i\sqrt{3}\right) ^{1/4}}\
\matrixb i{-i}11,&
\qquad\mbox{when}\qquad  \mu = -2i\sqrt 3.
\end{array}$$
 The transformations and their matrices are given in the
form they were computed by \Maple.
\end{ex}
Finally, the projective symmetry group of the quartic form $p^2+1$ consists of just two elements:
identity and
$p\rightarrow -p$.

\begin{ex}\label{Ex5} {\em Quintic forms}.
For polynomials of degree $5$, the projective symmetry group is either
cyclic, of type $\An$, or dihedral, of type $\Dn$.  Some representative examples  are listed in the
following table.

\begin{center}Projective Symmetry Groups of Representative Quintics
\hrule
\begin{tabbing}
aaaaaaaaaaaaaaaaaa\=aaaaaaaaaaa\=aaaaaaaaaaaaaaaaaaaaaa\=\kill
\>i.\>$p^5+1$\> ${\cal D}_5$\\
\>ii.\>$p^5+p$ \>${\cal A}_4$\\
\>iii. \>$p^5+p^2$\>${\cal A}_3$\\
\>iv.\>$p^5+p^3$\>${\cal A}_2$\\
\> v.\>$p^5+p^2+1$\>$\{e\}$\\
\>vi.\>$p^5-4p-2$\>$\{e\}$
\end{tabbing}
\vskip5mm
\hrule
\end{center}
\end{ex}

\begin{ex} \label{Ex6} {\em Higher degree forms}.  At the sixth degree, we first encounter a polynomial with an octahedral projective symmetry group:  the sextic $Q(p)=p^5+p$
which
corresponds to the form $Q(x,y)=x^5y+xy^5$,
compare with (\ref{phipsi}). The
inhomogeneous form looks
like the the second quintic polynomial listed in the preceding table, but we are now considering it
as a sextic with an additional root at $\infty $, and so the symmetry group is quite different.
Initially
\Maple\ produces symmetries which involve square roots and so do not initially look like linear
fractional transformations. However, after some simplifications under the radical we obtain the group of
linear fractional transformations generated by
$$i\,p, \qquad \qquad  \frac{\sqrt{2}(1+i)p-2}{\sqrt{2}(1-i)+2p}\>,$$
with corresponding matrices
$$\matrixb{i^{5/6}}00{i^{-1/6}},\qquad
\matrixb{\frac 1 2(1+i)}{-\frac 1 2\sqrt{2}}{\frac 1 2\sqrt{2}}{\frac 1 2(1-i)}.$$

The next time we meet this group is the octavic (degree $8$) form $Q(p)=p^8+14p^4+1$.
The octahedral generators are now
$$p \longmapsto i\, p,\qquad \qquad  p \longmapsto i\, \left ({p+1\over p-1}\right ),$$
which correspond to the matrix symmetries
$$\matrixb{i}001,\\
\frac{\sqrt 2}2 \matrixb ii1{-1}.$$
\end{ex}
\section{Differential Invariants for Ternary Forms.}
\label{ternary}
In this section  we  derive a fundamental set of differential invariants for
the projective action  (\ref{Af}) on  the inhomogeneous version $f(p,q)$ of
a ternary form $F(x,y,z)$.
Let us consider the corresponding (local) action of $GL(3,\C)$ on the graph of polynomial
$u=f(p,q)$:
 \begin{eqnarray}
\nonumber p&\mapsto &P= \frac{\alpha  p+\beta q +\gamma}{\delta p+\epsilon q+\zeta};\\
\label{gl3}q&\mapsto& Q=\frac {\lambda  p+ \nu  q+\tau}{\delta p+\epsilon q+\zeta};\\
\nonumber u&\mapsto&v=(\delta p+\epsilon q+\zeta)^{-n}u
\end {eqnarray}
The direct construction of the moving frame is computationally difficult so we
apply the recursive algorithm~\ref{AlgHi} from Chapter~\ref{chInd}.

\begin{rem} If Gr\"obner basis computation were more feasible in practice the following lengthy derivation of  moving frames and differential invariants would not be necessary. In theory, the  signature manifold could be derived from the lifted invariants $v_{k,l}$, obtained by the prolongation of the action (\ref{gl3}). Indeed,  restricted to the graph of a polynomial, lifted invariants $v_{k,l}$ are some rational functions of the group parameters and variables $p$ and $q$. The action (\ref{gl3}) if locally free and transitive on $J^3$ and we can choose a cross-section
$$\CS=\{ p=q=0,\, u=1,\, u_{1,0}=u_{0,1}= u_{2,0}=u_{0,2}=0,\,u_{3,0}=u_{0,3}=1\}.$$
Using an  algorithm based on Gr\"obner basis computations (see \cite{CLO96}, Theorem~2,
\S~3, ch. 3)
 we can eliminate group
parameters and variables $p$ and $q$ from the equations:
\begin{eqnarray*}
 &P= Q=0,\, v=1,\, v_{1,0}=v_{0,1}= v_{2,0}=v_{0,2}=0,\,v_{3,0}=v_{0,3}=1,\\
 &I_{1,1}=v_{1,1},\, I_{2,1}=v_{2,1},\,I_{1,2}= v_{1,2},\, I_{k,l}=v_{k,l},
\end{eqnarray*}
where $3<k+l<s+1$ and $s$ is the differential invariant  order of the polynomial defined in Chapter~\ref{prelim}, Definition~\ref{dRankOrder}.
As the result we obtain  polynomial relations among $I_{k,l}$ that  generate a prime
ideal.
The corresponding irreducible signature variety contains the signature manifold for
the graph of the  polynomials (see Remark~\ref{svariety}).
This procedure can be easily generalized
for an arbitrarily algebraic action on polynomials.
\end{rem}
In practice however we were  not able to carry out this straightforward algorithm and so we start the recursion  by
choosing
a  cross-section $\CS_0=\{p=0,q=0,u=1\}$ to the local action (\ref{gl3}) of $GL(3,\C)$
on $\C^3$
The isotropy group $H_1$
 of $\CS_0$  consists of matrices:
 $$\left(
\begin{array}{ccc}
\alpha & \beta &0 \\
\lambda&\nu& 0\\
\delta & \epsilon &1
\end{array}
\right) ,
$$
The action of the subgroup

 $$T=\left\{\left(
\begin{array}{ccc}
1 & 0 &t_1 \\
0&1& t_2\\
0 &0 &c\\
\end{array}
\right) \right\}
$$
is locally free and transitive and thus $G=H_1T$.
The prolongation of the $T$-action:
$$p\mapsto \frac{p+t_1} c,\qquad q\mapsto \frac{q+t_2} c, \qquad u\mapsto c^{-n}u$$
to $J^k$ is given by the formulas
\begin{equation}\label{Tu}
u_{k,l}\mapsto c^{k+l-n}u_{k,l},
\end{equation}
where the variable $u_{k,l}$ correspond to the derivative
$\frac{\partial^{k+l} u}{\partial^kp\partial^lq}$.
The condition
$\{p=0,q=0,u=1\}$ normalizes the group parameters of T:
\[t_1=-p,\quad t_2=-q,\quad c=u^{1/n}.\]
The  set $\CS_0^\infty$ such that $\pi^\infty_0(\CS_0^\infty)=\CS_0$ is parameterized by functions
$U_{k,l}=u^{\frac{k+l-n} n} u_{k,l} $ obtained from (\ref{Tu}) by substitution of the
 normalization for $c$.
The group $H_1$ acts on the functions $U_{k,l}$ restricted to  $\CS_0^\infty$ in the same manner as on the corresponding coordinate  functions.
We  note  that the  action of the subgroup $H_1$ on $\CS_0^1$ is transitive and
we choose  the cross-section  $\CS^1_1\subset \CS^1_0$ defined  by the condition
$U_{1,0}=0,\, U_{0,1}=0$.
The isotropy group $H_2\subset H_1$ for this cross-section is isomorphic to $GL(2,\C)$:
 $$H_2=\left\{\left(
\begin{array}{ccc}
\alpha& \beta &0   \\
\gamma&\delta& 0 \\
0 &0 &1\\
\end{array}
\right)\right\}.
$$
We can write $H_1$ is a product of $H_2$ and $R$, where the group $R$ consists of
inversions:
$$R=\left\{\left(
\begin{array}{ccc}
 1    & 0     &0   \\
0     &1     & 0  \\
a &b &1\\
\end{array}
\right)\right\}.
$$
To normalize the parameters of $R$ we  need to  prolong its action:
\begin{eqnarray*} \label{L}
p &\mapsto &\bar p=\frac p{ap+bq+1}, \\
q &\mapsto &\bar q=\frac q{ap+bq+1} ,\\
u &\mapsto &\bar  u=\left( ap+bq+1\right) ^{-n}u
\end{eqnarray*}
to $J^\infty$ and then restrict it to $\CS_0^\infty$.
The lifted invariant differential operators are dual to the lifted contact invariant forms:
\begin{eqnarray*}
d\bar{p} &=&\frac{dp+b(qdp-pdq)}{\sigma ^2}, \\
d\bar{q} &=&\frac{dq+a(pdq-qdp)}{\sigma ^2},
\end{eqnarray*}
where $\sigma=ap+bq+1$ and hence they are equal to:
\begin{eqnarray*}
\bar \CD_p &=&\sigma ^2D_p+\sigma q\left( aD_q-bD_p\right) , \\
\bar \CD_q &=&\sigma ^2D_q+\sigma p\left( bD_p-aD_q\right) .
\end{eqnarray*}
Thus the  first prolongation of the action of $R$ is:
\begin{eqnarray}\label{t1}
u_{1,0} &\mapsto&\bar u_{1,0}=\sigma ^{1-n}\left( -nau+\sigma u_p\right)
+\sigma ^{1-n}q\left( au_q-bu_p\right) , \\
u_{0,1} &\mapsto &\bar u_{0,1}=\sigma ^{1-n}\left( -nbu+\sigma u_q\right) +
\sigma ^{1-n}p\left( bu_p-au_p\right) .
\end{eqnarray}
By restricting  the above  transformations  to $\CS_0^1$, where $p=q=0$ and $u=1$
one obtains the transformations of $U_{1,0}$ and $U_{0,1}$:
\begin{equation}
U_{1,0} \mapsto \bar U_{1,0}= -na+ U_{1,0},  \qquad
U_{1,0} \mapsto \bar U_{0,1}= -nb+ U_{0,1}.
\end{equation}
Normalization $\bar U_{1,0}=0$ and $\bar U_{0,1}=0$ defines the group parameters of $R$:
\begin{equation}\label{ab}
a=\frac {U_{1,0}}n=\frac {u^{\frac{1-n}n}u_{1,0} }{n},\qquad
b=\frac {U_{0,1}}n=\frac {u^{\frac{1-n}n}u_{0,1} }{n}.
\end{equation}
To proceed further we need to determine the invariantization of higher order
functions    $U_{k,l}$ under the action
of  $R$ on $\CS_0^\infty$. In order to do so we need to compute the lifted invariants
 $\bar u_{k,l}=\bar \CD_p^k\bar \CD_q^l\,\bar u$, restrict
 them to $\CS_0^\infty$ and then substitute the normalization (\ref{ab}).
It is not hard to observe (see also (\ref{t1}))  that
$$\bar u_{1,0}=\bar \CD_p \bar u=\sigma^2D_p \bar u+q\times (\dots),$$
where the second term in the sum equals to  $q$ multiplied by an expression
and hence it  is equal to
zero when $q=0$. In general
$$\bar u_{k,0} (\bar\CD_p)^k v=\left(\sigma^2D_p)\right)^k \bar u+q\times (\dots).$$
The differential operator in the first part of the formula  is similar to
the one that  produces the prolongation  formulas  in the
case of binary forms and so  we obtain the familiar expressions:
$$\bar u_{k,0}=\sigma ^{k-n}\sum_{j=0}^k(-1)^{k-j}
{k \choose j} \frac{(n-j)!}{(n-k)!}a^{k-j}\sigma ^j u_{j,0} +q\times(\dots).$$
We restrict these functions to $\CS_0^\infty$, where $p=0,\,q=0,\, u=1$ and
$ u_{k,l}=U_{k,l}$
and substitute the normalization (\ref{ab}) to obtain functions:
\begin{eqnarray}
\nonumber Q_{k,0}&=&\sum_{j=0}^k(-1)^{k-j}
{ k\choose  j}
 \frac{(n-j)!}{(n-k)!}\left( \frac{ U_{1,0}}n\right) ^{k-j}U_{j,0}\\
\label{Qk0}&=&u^{\frac{k(1-n)}n}\sum_{j=0}^k(-1)^{k-j}
{k \choose  j}
 \frac{(n-j)!}{(n-k)!}u^{j-1}\left( \frac{ u_{1,0}}n\right) ^{k-j}u_{j,0},
\end{eqnarray}
In the notation of the Algorithm~\ref{AlgHi}
 the last expression corresponds  to the pull back of the coordinate function
$u_{k,0}$ under the map $\iota_1\colon J^k \rightarrow \CS_1^k$.
In the same manner we derive that
\begin{eqnarray}
\nonumber Q_{0,k}&=&\sum_{j=0}^k(-1)^{k-j}
{k \choose j}
 \frac{(n-j)!}{(n-k)!}\left( \frac{ U_{0,1}}n\right) ^{k-j}U_{0,j}\\
\label{Q0k} &=&u^{\frac{k(1-n)}n}\sum_{j=0}^k(-1)^{k-j}
{k \choose j}
 \frac{(n-j)!}{(n-k)!}u^{j-1}\left( \frac{ u_{0,1}}n\right) ^{k-j}u_{0,j}.
\end{eqnarray}
The similar straightforward (but more complicated) derivations produces the formulas for  $Q_{k,l}\, k\neq 0,\, l\neq 0$ corresponding to the mixed derivatives. We list the ones which we will use later:
\begin{eqnarray}\label{Qkl}
\nonumber  Q_{1,1}&=&u^{\frac{2(1-n)}n}\left[u_{1,1}\,u-\frac{n-1} n
u_{1,0}u_{0,1}\right],\\
\nonumber Q_{1,2}&=&u^{\frac{3(1-n)}n}\left[u_{1,2}\,u^2-\frac{n-2} n
(u_{0,2}u_{1,0}+2\,u_{1,1}u_{0,1})u
+ 2\frac{(n-1)(n-2)}{n^2}u_{0,1}^2u_{1,0}\right],\\
\nonumber Q_{2,1}&=&u^{\frac{3(1-n)}n}\left[u_{2,1}\,u^2-\frac{n-2} n
(u_{2,0}u_{0,1}+2\,u_{1,1}u_{1,0})\,u
+ 2\frac{(n-1)(n-2)}{n^2}u_{1,0}^2u_{0,1}\right],\\
\nonumber Q_{1,3}&=&u^{\frac{4(1-n)}n}\left[u_{1,3}\,u^3-\frac{n-3} n
(u_{1,0}u_{0,3}+ 3\,u_{1,2}u_{0,1})\,u^2+ 3\frac{(n-3)(n-2)}{n^2}\right.\\
&\times&\left.
(u_{0,1}^2u_{1,1}+u_{1,0}u_{0,1}u_{0,2})\,u-
3\frac{(n-1)(n-2)(n-3)}{n^3}u_{0,1}^3u_{1,0}\right],\\
\nonumber Q_{2,2}&=&u^{\frac{4(1-n)}n}\left[u_{2,2}\,u^3-2\frac{n-3}n
(u_{2,1}u_{0,1}+u_{1,2}u_{1,0})\,u^2 + \frac{(n-3)(n-2)}{n^2}\right.\\
\nonumber &\times&\left.(u_{1,0}^2u_{0,2}+4 u_{1,0}u_{0,1}u_{1,1}+u_{0,1}^2u_{2,0})u-3\frac{(n-1)(n-2)(n-3)}{n^3}u_{1,0}^2u_{0,1}^2\right],\\
\nonumber Q_{1,3}&=&u^{\frac{4(1-n)}n}\left[u_{3,1}\,u^3-\frac{n-3} n
(u_{0,1}u_{3,0}+3\,u_{2,1}u_{1,0})\,u^2+ 3\frac{(n-3)(n-2)}{n^2}\right.\\
\nonumber &\times&\left.
(u_{1,0}^2u_{1,1}+u_{0,1}u_{1,0}u_{2,0})\,u-
3\frac{(n-1)(n-2)(n-3)}{n^3}u_{1,0}^3u_{0,1}\right].
\end{eqnarray}
The function $Q_{k,l}$ are invariant  under the transformation  (\ref{gl3}) up to
action of $H_2\backsimeq GL(2,\C)$ (see Proposition~\ref{invupto}
from Chapter~\ref{chInd}).
Moreover restricted to $\CS_1^k$ they are
transformed by $H_2$ by the same formulas as the coordinate functions
 $u_{k,l}$. We  prolong the linear transformation:
$$p\mapsto\alpha p+\beta q,\quad q\mapsto \gamma p+\delta q,\quad u\mapsto u.$$
to $J^\infty$ and then substituting $u_{k,l}$ with $Q_{k,l}$ to obtain:
\begin{eqnarray}\label{Q20}
\bar Q_{2,0}&=&\frac 1 {\Delta^2}\left( \delta ^2Q_{2,0}-2\gamma \delta Q_{1,1}
+\gamma ^2Q_{0,2}\right);\\
\label{Q11} \bar Q_{1,1} &=&\frac 1 {\Delta^2}\left(-\delta \beta Q_{2,0}+
(\gamma \beta +\alpha
\delta)Q_{1,1}-\alpha \gamma Q_{0,2}\right); \\
\label{Q02}\bar Q_{0,2} &= &\frac 1 {\Delta^2}\left(\beta ^2Q_{2,0}-2\alpha \beta Q_{1,1}
+\alpha ^2Q_{0,2}\right); \\
\label{Q30}\bar Q_{3,0} &=&\frac 1 {\Delta^3}
\left(\delta ^3Q_{3,0}-3\gamma \delta ^2Q_{2,1}+3\gamma^2
\delta Q_{1,2}-\gamma ^3Q_{0,3}\right); \\
\bar Q_{2,1} &= &\frac 1 {\Delta^3}
\left(-\delta ^2\beta Q_{3,0}+\delta (2\gamma \beta +\alpha
\delta )Q_{2,1}-\gamma (\gamma \beta +2\alpha \delta )Q_{1,2}+\alpha \gamma
^2Q_{0,3}\right) \\
\bar Q_{1,2} &= &\frac 1 {\Delta^3}
\left(\delta \beta ^2Q_{3,0}-\beta (\gamma \beta +2\alpha
\delta )Q_{2,1}+\alpha (2\gamma \beta +\alpha \delta )Q_{1,2}-\alpha
^2\gamma Q_{0,3}\right); \\
\label{Q03}\bar Q_{0,3} &= &\frac 1 {\Delta^3}
\left(-\beta ^3Q_{3,0}+3\alpha \beta ^2Q_{2,1}-3\alpha
^2\beta Q_{1,2}+\alpha ^3Q_{0,3}\right);\\
\nonumber && etc .
\end{eqnarray}
We can normalize the remaining group parameters by setting
$$
\bar Q_{2,0}=\bar Q_{0,2}=0 \quad \mbox{and} \quad \bar Q_{3,0}=\bar Q_{0,3}=1.
$$
From the first pair of normalizations it follows that
$\frac \delta \gamma $ and $\frac \beta \alpha $ are two
roots of the same quadratic equation (see (\ref{Q20}) and (\ref{Q02})) so we can write that
 \begin{equation}\label{r1r2}
\frac \delta \gamma=r_1 =\frac{Q_{1,1}+\sqrt{Q_{1,1}^2-Q_{2,0}Q_{0,2}}}{Q_{2,0}},\qquad
\frac \beta \alpha=r_2 =\frac{Q_{1,1}-\sqrt{Q_{1,1}^2-Q_{2,0}Q_{0,2}}}{Q_{2,0}}.
\end{equation}
By subtracting these expressions one  obtains that
$$
r_1-r_2=\frac {\alpha\delta-\beta\gamma}{\alpha\gamma}=\frac
{2\sqrt{d}}{Q_{2,0}} \Longrightarrow \ \Delta=\alpha\gamma\,\frac
{2\sqrt{d}}{Q_{2,0}},$$ where $d=Q_{1,1}^2-Q_{2,0}Q_{0,2} $. From
the second pair of normalizations we obtain that:
\begin{eqnarray}\label{alphagamma}
\alpha&=&\frac{Q_{2,0}}{2\sqrt{d}}\left(r_1^3Q_{3,0}-3r_1 ^2Q_{2,1}+3r_1
Q_{1,2}-Q_{0,3}\right)^{1/3},\\
\gamma&=&\frac{Q_{2,0}}{2\sqrt{d}}\left(-r_2^3Q_{3,0}+3r_2 ^2Q_{2,1}-3r_2
\nonumber Q_{1,2}+Q_{0,3}\right)^{1/3}.
\end{eqnarray}
The substitution  of this normalization into $\bar{Q}_{1,1},\bar{Q}_{2,1},\bar{Q}_{1,2}$ and higher
order `derivatives' $\bar{Q}_{k,l},\, k+l>3$ produces  all the fundamental invariants of the action (\ref{gl3}).
These invariants however are not rational expressions in $p$ and $q$ and so we can not
use  Gr\"obner basis elimination algorithm to describe the corresponding signature
 manifold. Fortunately, by turning back  to the classical  invariant theory processes we are able
to derive a complete set of fundamental rational invariants.

We first note that $\bar{Q}_{k,l}$
are lifted invariants under the action of $H_2\backsimeq GL(2,\C)$ on $H_2\times J^{\infty}$
defined by the prolongation of the map:
\begin{equation}\label{Apalpha}
\left(\begin{array}{c}p\\q \end{array}\right)\mapsto
\left(\begin{array}{cc} a_{11}&a_{12}\\a_{21}& a_{22} \end{array}\right)\left(\begin{array}{c}p\\q \end{array}\right),\qquad
\left(\begin{array}{cc}\alpha& \beta \\ \gamma&\delta\end{array}\right)\mapsto
\left(\begin{array}{cc}\alpha& \beta \\ \gamma&\delta\end{array}\right)\left(\begin{array}{cc} a_{11}&a_{12}\\a_{21}&a_{22} \end{array}\right)^{-1}.
\end{equation}
We can also look at  $\bar{Q}_{k,l}$ as on polynomials in
$\{\alpha,\, \beta, \,\gamma,\, \delta\}$ with coefficients $Q_{i,j}$. We note that the coefficients $Q_{i,j}$ are transformed in the same manner as the coefficients of the binary forms $\sum_{i,j} a_{i,j}p^iq^j$, where $i+j=n$, and $p$ and $q$ are transformed as in (\ref{Apalpha}).  From the second formula in (\ref{Apalpha}) we can see that vectors $\left(\begin{array}{c}\alpha\\  \beta \end{array}\right)$ and
$\left(\begin{array}{c} \gamma\\  \delta\end{array}\right)$ are transformed by the matrix
$$\left(\begin{array}{cc} a_{11}&a_{12}\\a_{21}&a_{22} \end{array}\right)^{-t}=
\frac 1 {\det(A)}\left(\begin{array}{cc} a_{22}&-a_{21}\\-a_{12}&a_{11} \end{array}\right).$$
Such vectors are called  {\em covariant} in contrast with {\em contravariant} vectors which are transformed by the left multiplication by $A$
(for instance, $p$ and $q$ form a contravariant vector).  The forms
$$P_2=\Delta^2\,\bar Q_{2,0},\quad
P_3=\Delta^3\,\bar Q_{3,0},\quad  P_4=\Delta^4\,\bar Q_{4,0},\quad\dots$$
are relatively invariant under the simultaneous transformations of $\left(\begin{array}{c} \gamma\\  \delta\end{array}\right)$ and $Q_{k,l}$ with  weights $2,\, 3$ and $4$ respectively. Because of  their dependence  on a covariant vector forms $P_i$  are called
{\em contravariants} (see \cite{Gur64} \S~13 for more details on classical terminology).
In the case of binary form there is a duality between covariant vectors and contravariant vectors. Indeed let $\mu=-\delta$ and  $\eta=\gamma$.  Then $\mu$ and $\eta$ form a
contravariant vector of weight $-1$:
$$\left(\begin{array}{c}\mu\\ \eta \end{array}\right)\mapsto
\frac 1 {\det(A)}\left(\begin{array}{cc} a_{11}&a_{12}\\a_{21}& a_{22} \end{array}\right)\left(\begin{array}{c}\mu\\\eta \end{array}\right).$$
Thus the  polynomials $P_i$ rewritten  in terms of the new variables $\mu$ and $\eta$ are covariants of weight zero:
\begin{eqnarray}\label{P2}
P_2&=& Q_{2,0}\,\mu ^2+2 Q_{1,1}\,\mu\eta+Q_{0,2}\,\mu^2;\\
\label{P3} P_3 &= &
Q_{3,0}\,\mu^3+3 Q_{2,1}\,\mu^2\eta^2+3\,Q_{1,2}\,\mu\eta^2+Q_{0,3}\,\eta^3;\\
\label{P4} P_4 &= &
Q_{4,0}\mu^4+4 Q_{3,1}\mu^3\eta +6\,Q_{2,2}\mu^2\eta^2+4\,Q_{1,3}\mu\eta^3+Q_{0,4}\eta^4;\\
\nonumber && etc .
\end{eqnarray}
Let $H_\kappa(\dots, Q_{k,l},\dots)$ be a rational invariant of  these forms.
We recall that $Q_{k,l}$ are  differential functions (see formulas (\ref{Qk0}),
(\ref{Q0k}) and (\ref{Qkl})) and so functions   $H_\kappa(\dots, Q_{k,l},\dots)$ provide differential invariants of the initial action (\ref{gl3}). We would like to find sufficiently many of such invariants in order to parameterize the signature manifold.
We start with   the following classical  relative invariants \cite{GraceYoung03}.
\begin{description}
\item[] The discriminant of quadratic $P_2$:
$$d_2= Q_{2,0}Q_{0,2}-Q_{1,1}^2$$
 which is a relative invariant of weight 2.
\item[]
The discriminant of the cubic $P_3$:
$$
d_3=Q_{3,0}^2\,Q_{0,3}^2-3\,Q_{2,1}^2Q_{1,2}^2-6\,Q_{3,0}Q_{2,1}Q_{1,2}Q_{0,3}+
4\,Q_{3,0}Q_{1,2}^3+4\,Q_{0,3}Q_{2,1}^3
$$
which is a relative invariant of weight 6.
\item[] Two invariants of the quartic $P_4$:
$$
i:=Q_{4,0}Q_{0,4}-4\,Q_{3,1}Q_{1,3}+3\,Q_{2,2}^2
$$
of weight 4, and
$$
j=\det\left(\begin{array}{ccc}
Q_{0,4}& Q_{1,3}&Q_{2,2}\\
Q_{1,3}& Q_{2,2}&Q_{3,1}\\
Q_{2,2}& Q_{3,1}&Q_{4,0}
\end{array}\right)
$$
of weight 6.
\end{description}

To obtain joint relative invariants of forms $P_i$ we can apply the  omega process,
or transvection
(as described in chapter~3, \S~48 of  \cite{GraceYoung03} and chapter~5 of \cite{O99})
to these forms.
 Let $\Phi$ and $\Psi$ be covariants of  weight $k_1$ and $k_2$ respectively then
their  $r$-th transvectant  $(\Phi,\Psi)^{(r)}$  is a covariant of weight $k_1+k_2+r$.

Let
\begin{eqnarray*}
H_3=(P_4,P_4)^{(2)},\quad H_4=(P_4,P_4)^{(2)}, \\
T_3=(H_3,P_3)^{(1)},\quad T_4=(H_4,P_4)^{(1)},\\
 S=(H_4,P_3^2)^{(3)}.
\end{eqnarray*}
We note that $H$'s  have weight $2$ and they are the Hessians of the corresponding
 forms, $T$'s  have weight $3$ and their explicit formulas (\ref{HTU}) are given
in the preceding section on binary forms.
We also note that the  invariants of single forms  can be expressed as transvectants:
\begin{eqnarray*}
d_2=\frac 1 8 (P_2,P_2)^{(2)}&d_3= \frac 1 {10368} (H_3,H_3)^{(2)}\\
i=\frac 1 {1152}(P_4,P_4)^{(4)}&j= \frac 1 {497664}(H_4,P_4)^{(4)}
\end{eqnarray*}
We  complete this list with the following joint covariants:
\begin{description}
\item[] Joint invariants of cubic $P_3$ and quadratic $P_2$:
$$M_1=\frac 1 {288}(H_3^2,P_2)^{(2)},\qquad M_2=\frac 1 {103680}(P_3^2,P_2^3)^{(6)}$$
of weights $4$ and  $6$ respectively.
\item[]Joint invariants of quadratic $P_2$ and quartic $P_4$:
$$M_3=\frac 1 {576}\,(P_4,P_2^2)^{(4)},\quad
M_4=\frac 1 {1194393600}\,(T_4,P_2^3)^{(6)}$$
of weights $4$,   and $9$ respectively.
\item[] Joint invariant of cubic $P_3$ and quartic $P_4$:
$$M_5=\frac 1 {238878720}(S,P_4)^{(4)}$$
of weight 9.
\end{description}
Taking into account the  weights of the relative invariants above we define three
absolute rational invariants of the third order:
\begin{equation}
\label{third} I_1=\frac{M_1}{{d_2}^2},\qquad I_2=\frac {M_2} {{d_2}^3},\qquad
I_3=\frac {d_3}{{d_2}^3}
\end{equation}
and five invariants of the fourth order:
\begin{equation}
\label{fourth}I_4=\frac j {{d_2}^3},\quad I_5=\frac{i}{{d_2}^2},\quad
I_6=\frac {{M_4}^2} {{d_2}^9},\quad I_7=\frac {M_3} {{d_2}^2},\quad  I_8=\frac {{M_5}^2}{{d_2}^9}.
\end{equation}
Thus we have found eight differential invariants of order four or less.
Since most of the explicit formulas are long we place them in the Appendix~B.
Using the Thomas replacement theorem \cite{FO99} we can rewrite these invariants
in terms of eight independent invariants obtained by invariantizations
$I_{k,l}=\inv(u_{k,l})$ (see Section~\ref{smf} of Chapter~\ref{prelim}) of the
`derivative' coordinates $u_{k,l}$:

\begin{eqnarray*}
I_1&=&{\displaystyle\frac {  I_{2, \,1}\,I_{1, \,2}-1}{{I_{1, \,1}}^{3}}} ,\\
I_2&=& - 4\,{\displaystyle \frac {1 + 9\,{I_{2, \,1}}\,{I_{1, \,2}}}
{{I_{1, \,1}}^{3}}} ,\\
I_3&=&{\displaystyle \frac {1 - 3\,{I_{2, \,1}}^{2}\,{I_{1, \,2}}^{2}
 - 6\,{I_{2, \,1}}\,{I_{1, \,2}} + 4\,{I_{1, \,2}}^{3} + 4\,{I_{2
, \,1}}^{3}}{{I_{1, \,1}}^{6}}} ,  \\
I_4&=&{\displaystyle \frac {{I_{4, \,0}}\,{I_{0, \,4}} - 4\,{I_{3
, \,1}}\,{I_{1, \,3}} + 3\,{I_{2, \,2}}^{2}}{{I_{1, \,1}}^{4}}}
,  \\
I_5&=&{\displaystyle \frac { + {I_{0, \,4}}\,{I_{2, \,2}}\,{I_{
4, \,0}} - {I_{0, \,4}}\,{I_{3, \,1}}^{2} - {I_{1, \,3}}^{2}\,{I
_{4, \,0}} + 2\,{I_{1, \,3}}\,{I_{2, \,2}}\,{I_{3, \,1}} - {I_{2
, \,2}}^{3}}{{I_{1, \,1}}^{6}}} ,  \\
I_6&=& 16\,{\displaystyle \frac {({I_{1, \,3}}^{2}\,{I_{4, \,0}} -
{I_{0, \,4}}\,{I_{3, \,1}}^{2})^{2}}{{I_{1, \,1}}^{12}}},\\
I_7&=&4\,{\displaystyle \frac {{I_{2, \,2}}}{{I_{1, \,1}}^{2}}},\\
I_8&=& {\displaystyle \frac 1 {{I_{1, \,1}}^{18}}}\,
 ( - 9\,{
I_{3, \,1}}\,{I_{2, \,2}}\,{I_{2, \,1}}^{2}\,{I_{0, \,4}} - 6\,{I
_{1, \,3}}^{2}\,{I_{2, \,1}}\,{I_{2, \,2}},
+ 6\,{I_{3, \,1}}^{2}\,{I_{1, \,2}}\,{I_{2, \,2}} -
2\,{I_{3, \,1}}^{3} + 2\,{I_{1, \,3}}^{3}\\
 & & + 2\,{I_{4, \,0}}\,{I_{1, \,3}}\,{I_{1, \,2}}\,{I_{0, \,4}}
- 6\,{I_{3, \,1}}\,{I_{2, \,2}}\,{I_{1, \,2}}\,{I_{0, \,4}}  \\
 & &  + 9\,{I_{2, \,2}}^{2}\,{I_{2, \,1}}\,{I_{0, \,4}} +
{I_{4, \,0}}^{2}\,{I_{0, \,4}}\,{I_{1, \,2}} - {I_{4, \,0}}\,{I_{
0, \,4}}^{2}\,{I_{2, \,1}} + 3\,{I_{4, \,0}}\,{I_{1, \,3}}\,{I_{2
, \,1}}^{2}\,{I_{0, \,4}} \\
 & &  + 2\,{I_{1, \,3}}\,{I_{3, \,1}}\,{I_{1, \,2}}\,{I_{4
, \,0}} - 6\,{I_{3, \,1}}^{2}\,{I_{1, \,2}}^{2}\,{I_{1, \,3}} - 4
\,{I_{3, \,1}}^{2}\,{I_{2, \,1}}\,{I_{1, \,3}} - 9\,{I_{4, \,0}}
\,{I_{2, \,2}}^{2}\,{I_{1, \,2}} \\
 & &  + 3\,{I_{4, \,0}}\,{I_{2, \,2}}\,{I_{3, \,1}} + 4\,{
I_{1, \,3}}^{2}\,{I_{3, \,1}}\,{I_{1, \,2}} - 9\,{I_{4, \,0}}\,{I
_{1, \,3}}^{2}\,{I_{1, \,2}}\,{I_{2, \,1}} + 9\,{I_{3, \,1}}^{2}
\,{I_{2, \,1}}\,{I_{1, \,2}}\,{I_{0, \,4}} \\
 & &  + 9\,{I_{4, \,0}}\,{I_{2, \,2}}\,{I_{1, \,2}}^{2}\,{
I_{1, \,3}} - 3\,{I_{2, \,2}}\,{I_{1, \,3}}\,{I_{0, \,4}} - 3\,{I
_{4, \,0}}\,{I_{0, \,4}}\,{I_{1, \,2}}^{2}\,{I_{3, \,1}} + 6\,{I
_{4, \,0}}\,{I_{2, \,2}}\,{I_{2, \,1}}\,{I_{1, \,3}} \\
 & &  - 2\,{I_{4, \,0}}\,{I_{0, \,4}}\,{I_{2, \,1}}\,{I_{3
, \,1}} - 2\,{I_{1, \,3}}\,{I_{3, \,1}}\,{I_{2, \,1}}\,{I_{0, \,4
}} + {I_{0, \,4}}\,{I_{3, \,1}}^{2} - {I_{1, \,3}}^{2}\,{I_{4, \,
0}} + {I_{3, \,1}}\,{I_{0, \,4}}^{2} \\
 & &  + 6\,{I_{1, \,3}}^{2}\,{I_{3, \,1}}\,{I_{2, \,1}}^{2
} - {I_{4, \,0}}^{2}\,{I_{1, \,3}})^{2}
\end{eqnarray*}
By computing (with the help of a computer) the corresponding Jacobian
we conclude that  invariants $\{I_1,\dots,I_8\}$ are functionally independent.
On the other hand, since $GL(3,\C)$ acts freely  on $J^4$ and $\dim (J^4)=17$,
 there could be
 no more than eight functionally independent invariants and thus $\{I_1,\dots,I_8\}$
form a complete set of differential invariants of order four or less.
\begin{rem}\label{notdefined} None of the invariants is defined when $I_{1,1}\equiv 0$ (or equivalently $d_2\equiv 0$). By substitution of the group parameters (\ref{r1r2}) and (\ref{alphagamma}) into
(\ref{Q11}) we conclude that it happens if and only if the inhomogenization of the
 Hessian:
$$n\,f (f_{pq}^2-f_{pp}\,f_{qq})+(n-1)\,(f_{pp}\,f_q^2+2\,f_p\,f_q\,f_{pq}+f_{qq}\,f_p^2)   $$
is identically zero, and hence if and only if the ternary form can be transformed into a binary form.
\end{rem}

In the next section we will use the first three  invariants  to obtain a classification of ternary cubics and their group of symmetries. In the last section we use all eight invariants to construct the signature manifold for the forms $x^n+y^n+z^n$, therefore obtaining
 necessary and sufficient condition for the equivalence of an arbitrary ternary form of degree $n$   to the sum of $n$-th powers.


\section{Classification of Ternary Cubics}

In this section we reproduce  known results on classification of ternary cubics up
to a linear transformation (\cite{Gur64}, \cite{Kraft87}) and then obtain a classification of their symmetry groups, which we believe is new. We achieve this
by  restricting invariants $I_1,I_2$ and $I_3$, obtained in the previous  section to each if the  canonical forms. Note that the fourth order invariants, restricted  to a cubic, are zero.  By Gr\"obner basis computation we find the  ideal, whose zero set
 defines the corresponding signature manifold  and determine its dimension.
In non-trivial situations we find the dimension using {\sc Maple} function
{\tt hilbertdim},  based on computing  the degree of Hilbert  polynomial \cite{CLO96},
\cite{Froberg97}. Using Theorem~\ref{symm}, from Chapter~\ref{prelim} we make a conclusion about the dimension of the
 symmetry group, which helps us to determine the group explicitly.
In the case of a finite symmetry group we find its cardinality.
We  start with

{\large {\bf Reducible Cubics.}}

A reducible cubic is either a product of three linear factors or a product of linear and quadratic factors. We state the following classification theorem:
\begin{thm} \label{redc} A reducible cubic $F(x,y,z)$ is equivalent under a linear change of
variables to one of the following forms:
\begin{description}
 \item[1.] If it  is a product of there linear factors and
\begin{description}
\item[a)] all three factors are the same, then it is equivalent to $x^3$ and its symmetry group, is conjugate to a four-dimensional group isomorphic to $GL(2,\C)$, of linear transformations on the variables $y$ and $z$;
\item[b)] two factors are the same, then the cubic is equivalent  to  $x^2y$
and its symmetry group is conjugate to four-dimensional group of matrices
$$
\left( \begin{array}{ccc} \alpha&&\\
&\frac 1 {\alpha^2}&\\
\beta&\gamma&\delta
\end{array}\right);
$$
\item[c)] three  factors are linearly dependent, but
any pair of them is linearly independent, then the cubic is
equivalent  to  $xy(x+y)$ and its
 symmetry group is conjugate to the three-dimensional direct
 product of
 arbitrary linear transformations  $z \mapsto \alpha x+\beta y+\gamma
 z$ and a finite subgroup  of  order  $6\times 3$   of linear
 transformations on the $(x,y)$-plane, which preserves $xy(x+y)$ (see Section~\ref{bf}) .

\item[d)] there factors are linearly independent  factors,  then the cubic is equivalent  to  $xyz$ and its
 symmetry group is conjugate to two-dimensional group of matrices
\begin{equation}\label{r1c}
 \left(\begin{array}{ccc} \alpha&&\\
&\beta&\\
&&\frac 1 {\alpha\beta}
\end{array}\right);
\end{equation}
\end{description}
 \item[2.] If it  is a product of quadratic and linear factors then there are two
canonical forms:
 \begin{description}
\item[a)] $F(x,y,z)\sim  z\,(x^2+yz)$. In this case the symmetry group is conjugate to a two-dimensional group generated by
\begin{equation}\label{r2a}
\left(\begin{array}{ccc} 1&0&\alpha\\
-2\alpha &1&-\alpha^2\\
     0&0&1
\end{array}\right) \quad
\mbox{and}\quad  \left(\begin{array}{ccc} \beta&&\\
&\beta^4&\\
&&\frac 1 {\beta^{-2}}
\end{array}\right);
\end{equation}
\item[b)] $F(x,y,z)\sim  z\,(x^2+y^2+z^2)$. In this case the symmetry group is conjugate
 to  a one-dimensional group of rotations around the $z$-axis which is
isomorphic to  $O(2,\C)$.
\end{description}
\end{description}
\end{thm}
 {\em Proof.}

1. The classifications of the cubics reducible into linear factors is obvious.
We note only that a cubic has  repeated factors  (cases 1.~a) and 1.~b))
if and only if its   Hessian is identically  zero and so the invariants $I_1, I_2$
and $I_3$ are not defined (see Remark~\ref{notdefined}).  The graph $S$   of  such
cubic is a  totally singular
submanifold in $\C^4$, in a sense  that there is no  order of
prolongation  $n$ at which the  prolonged
action of the group  becomes locally free  on $j^n(S)$.
Thus, none of the constructions from Chapter~\ref{prelim} can be applied to this case.
The graphs of all  other cubics define regular submanifolds and so, due to
Theorem~\ref{symm}, the dimensions  of their  symmetry groups equal to
$2-\dim{\cal C}(S_f)$
where $f(p,q)$ is the inhomogeneous version of  $F(x,y,z)$.  The signature manifold
 ${\cal C}(S_f)$,  is  parameterized by $I_1(f),\, I_2(f)$ and $I_3(f)$.
When the cubic is equivalent to $xyz$ (case 1.~c)) the signature manifold consists of
a single point:
\[I_1=\frac 4 3,\qquad I_2=\frac {16} 3,\qquad I_3=\frac {16} 9.\]
and so the symmetry group is two-dimensional.
It is easy to check that the group (\ref{r1c})  leaves $xyz$  unchanged.

2. The general form of a cubic reducible into linear and quadratic factors is:
\[
F(x,y,z)=(k_1x+k_2y+k_3z)(k_4z^2+k_5zx+k_6zy+K(x,y))
,\]
where $k_i$ are some complex coefficients and $K(x,y)$ is a non-zero quadratic form
in two
variables.
By taking the first factor as a new variable $z^{\prime }$, one obtains  an
equivalent
cubic
\[
z\,(k_4z^2+k_5zx+k_6zy+K(x,y))
\]
where $k_i$ are some new coefficients and $K(x,y)$ is a new non-zero  quadratic
form. It is
known that by a linear change of variables the quadratic $K(x,y)$ can be
transformed to  either $x^2+y^2$ or $x^2$.  Thus $F(x,y,z)$ is
equivalent to either:
\begin{equation}\label{2a}
z\,(k_2z^2+k_3zx+k_4zy+x^2+y^2)
\end{equation}
or
\begin{equation}\label{2b}
z\,(k_2z^2+k_3zx+k_4zy+x^2),
\end{equation}
where $k_i$ are again new coefficients.
\newline
In the first case (\ref{2a}) we make the transformation
$$
x^{\prime } =(x+\frac{k_3}2z);\qquad y^{\prime } =(y+\frac{k_4}2z); \qquad
z^{\prime } =z
$$
to obtain  an equivalent form
\[
z\,(kz^2+x^2+y^2).
\]
Finally the  scaling:
$$k^{1/6}x^{\prime }=x, \qquad k^{1/6}y^{\prime}=y\qquad k^{-1/3}z^{\prime }=z$$
leads to the  canonical form:
\[
z\,(z^2+x^2+y^2).
\]
In the second case (\ref{2b}) we make the transformation
\begin{equation*}
x^{\prime } =(x+\sqrt{k_2}z); \qquad
y^{\prime } =k_4y+(k_3-2\sqrt{k_2})x; \qquad
z^{\prime } =z/k_4.
\end{equation*}
to obtain  an equivalent cubic:
\[
kz\,(zy+x^2).
\]
Finally, we make the transformation: $z^{\prime }=kz,$ $y^{\prime }=\frac
1 k y$, which leads  to
the canonical form:
\[
z\,(zy+x^2)
\]
We compute the corresponding signature manifolds. The signature manifold of a cubic
  equivalent to $z\,(x^2+yz)$, consists of a single  point:
  \[I_1=-\frac 1 6,\qquad  I_2=\frac {41}6,\qquad  I_3=-\frac 2 9,\]
and hence the symmetry group is two-dimensional. Using Lie's criterion~\ref{Lie}
we found infinitesimal symmetries:
$$ \left( \begin{array}{ccc} 0&0&1\\
-2&0&0\\
0&0&0
\end{array}\right) \quad
\mbox{and}\quad  \left(\begin{array}{ccc} 1&&\\
&4&\\
&&-2
\end{array}\right);
$$
that give rise to the group (\ref{r2a}).

The signature manifold of a cubic  equivalent to $z\,(x^2+y^2+z^2)$, is defined by three equations found by Gr\"obner basis computations:
\begin{eqnarray*}
&&
-1482\,I_3\,I_2+8865\,{I_3}^2+40\,{I_2}^2-1296\,I_1-36\,I_2+17280\,I_3-18522\,I_3\,I_1,\\
&&  40\,I_1\,I_2+582\,I_3\,I_1+42\,I_3\,I_2-315\,I_3^2-144\,I_1-4\,I_2-480\,I_3,\\
&& 360\,I_1^2-378\, I_3 I_1-18\, I_3\,I_2+135\,{I_3}^2-144\,I_1-4\,I_2+120\,I_3
\end{eqnarray*}
The the corresponding  variety is one-dimensional and thus the symmetry group is also one-dimensional. We notice  that the rotations with respect to the $z$-axis  leave the canonical form unchanged.
\qed

We proceed with the classification of the

{\large {\bf Irreducible Ternary Cubics:}}

It is well known that any homogeneous irreducible cubic $F(x,y,z)$ over $\C$
can be transformed by a linear map  into Weierstrass normal  form \cite{Knapp92}.
Let $V$ be
 the set of zeros of the inhomogenization $f(p,q)$  in
$\C P^2$. The normal form is obtained by transforming  one of the
inflection points of $V$ to the infinite point $(0,1,0)$,
and the  tangent line at this point to the line $(k,1,0)$ at infinity.
We state the following classification theorem:
\begin{thm} \label{irredc} An irreducible cubic $F(x,y,z)$ can be transformed
 under a linear change of  variables to one of the following forms:
\begin{description}
 \item[1] If f(p,q) defines a singular variety $V$ then it is equivalent to
either
\begin{description}
\item[a)] $p^3-q^2$ and the it has one-dimensional symmetry group given by:
\begin{equation}\label{1a}
\left( \begin{array}{ccc} 1&&\\
&{\alpha}&\\
&&\frac 1 {\alpha^{2}}
\end{array}\right);
\end{equation}
or
\item[b)]  $p^2(p+1)-q^2$, which has a discrete symmetry group,
consisting of $6$ projective symmetries
(see Definition~\ref{prsymm} of Section~\ref{symp}) which correspond to $18$ genuine symmetries.
\end{description}
 \item[2]  If $f(p,q)$ defines a nonsingular variety $V$ then it either equivalent to:
 \begin{description}
\item[a)] a cubic  in  one-parametric family:
$$p^3+ap+1-q^2$$
and then it has $18$ projective symmetries.

or
\item[b)] it is equivalent to $p^3+p-q^2$. In this case it has $36$ projective  symmetries,

or
\item[c)] it is equivalent to $p^3+1-q^2$. In this case it has $54$ projective  symmetries. We note that the cubic $p^3+q^3+1$ belongs to this class.
\end{description}
\end{description}
\end{thm}

For the proofs of the classification theorems we refer the reader to
 \cite{Knapp92}, \cite{Kraft87}, and restrict ourselves to the
discussion of the signature manifolds and the symmetry groups.
In the case 1.~a) the signature manifold is defined by two equations:
\[I_1=-\frac 1 6,\qquad 6\,I_2+45\,I_3-31=0 \]
and so the symmetries form a  one-dimensional group conjugate to (\ref{1a}). We note that the number of unimodular symmetries is finite.

 In all other cases the symmetry group is finite.
In theory, the projective symmetries can be found explicitly by solving the equations
\begin{equation}\label{symmN}
I_1(p,q)=I_1(P,Q),\quad I_2(p,q)=I_2(P,Q), \quad I_3(p,q)=I_3(P,Q),
\end{equation}
for $P$ and $Q$ in terms of $p$ and $q$. All solutions must be linear fractional
expressions in $p$ and $q$.  In practice however,  we were not able to carry out these
 computations. Nevertheless we can find the cardinality of the symmetry group using
 a well known algebraic geometry result (\cite{CLO96}, Proposition~8, ch.~5, \S~3).
\begin{prop}
Let ${\bf I}\subset \C [x_1,\dots, x_n]$ be an ideal such that its zero  set $V$  is finite then
\begin{description}
\item[(i)] The dimension of $(\C [x_1,\dots, x_n]/{\bf I})$ (as a vector space over $\C$) is finite and greater or equal to the number of points in $V$.
\item[(i)] If ${\bf I}$ is a radical ideal then equality holds, i.e., the number of
points in $V$ is exactly $(\C [x_1,\dots, x_n]/{\bf I})$.
\end{description}\end{prop}
For generic values of $P$ and $Q$ we find the number of the  solutions for
(\ref{symmN}) by computing the dimension of
the quotient space of  $\C[p,q]$ by the corresponding ideal. We use two algorithms
presented in the exercises for  ch.~2, \S~2  of \cite{CLO97}, first to check that the ideal defined by (\ref{symmN})  is radical, and then to compute the dimension of the quotient.

The nonsingular irreducible ternary cubics are known as {\em elliptic curves}
 and play an important role in number theory. The number of symmetries for these curves
 has a natural
explanation. First let us  consider only the symmetries fixing the point
$(0,1,0)$
and mapping  the line at infinity to itself. Following Knapp, \cite{Knapp92} we  call
these
symmetries  {\it  admissible}. It is not difficult to prove \cite{Knapp92}
that there are only $2$ such symmetries in case 2.a), $4$ in case 2.b) and $6$ in
case 2.c). Each of  nonsingular irreducible cubics has 9 inflection points. Each of the inflection point can be mapped to  $(0,1,0)$ with the
corresponding tangent line mapped to the line at infinity.
{\em We observe that the number of the projective symmetries of a nonsingular irreducible cubic equals to the number of its inflection point times the number of its
admissible symmetries.}

Ternary cubics form a  ten-dimensional linear space. So by dimensional consideration we expect to have one absolute  invariant depending on the  coefficients of cubic.
Indeed such invariant is known
\cite{Gur64}, \cite{Knapp92} and has been used to obtain classification results.
This invariant must be expressible in terms of the invariants $I_1,I_2$ and $I_3$ but we have not tried
to obtain the explicit formula. Not all  of the orbits of $GL(3,\C)$ acting on
the space of
cubics  are closed. In fact only the orbits of elliptic curves, the orbit of $xyz$ and
the orbit of $x^3$ are. See  Kraft \cite{Kraft87} for the proofs and the description
of how some orbit is included in the closure of the others.

We conclude this section with a simple corollary from  Theorems~\ref{redc} and \ref{irredc}:
\begin{cor}
A cubic in three variables splits into a linear factor and an irreducible
quadratic factor if and only if its $ SL(3,\C)$  symmetries form a
one-dimensional  Lie group.
\end{cor}
{\sc Maple} computations for ternary cubics can be found in Appendix~C.

\section{The Signature Manifold for $x^n+y^n+z^n$.}

In this section we construct the signature manifold for $f(p,q)=p^n+q^n+1$
therefore determining the necessarily and sufficient condition for a ternary form
to be complex equivalent to the sum of $n$-th powers.
We first  recall that the action (\ref{gl3}) becomes free on the third order jet space $J^3$. We restrict the first three invariants to $f(p,q)=p^n+q^n+1$:

\begin{eqnarray*}
I_1(f)&=&-\frac{(n-2)^2}{n(n-1)}\\
I_2(f)&=&-\frac {(n-2)^2}{n(n-1)}\\
&\times &\frac{5\,q^n\,(p^n)^2+5\,(p^n)^2-26\,p^n\,q^n+5\,(q^n)^2\,p^n+5\,p^n+5\,q^n
+5\,(q^n)^2}{p^n\,q^n}\\
I_3(f)&=&-\frac {(n-2)^4}{n^2(n-1)^2}\,\frac{(p^n+1-q^n)(p^n-1+q^n)(p^n-1-q^n)}
{p^n\,q^n}
\end{eqnarray*}
We observe that the first invariant is constant and we check that the  last two
are functionally
independent.
Thus the third order signature manifold has maximal possible dimension and the symmetry group of $f(p,q)$ is finite. Moreover, we can conclude that the differential invariant order of the graph $u=p^n+q^n+1$ equals to three and so both equivalence and symmetry problems can be solved by construction the fourth order signature manifold.

For $n>3$ we need to  restrict
the remaining five (forth order) invariants  to the graph  $u=p^n+q^n+1$:
\begin{eqnarray*}
I_4(f)&=&-\frac{(n-2)^2(n-3)^2}{n^2(n-1)^2}\,
\frac{\left((p^n)^3+1+(q^n)^3\right)}{p^n\,q^n}\\
I_5(f)&=&-\frac{(n-2)^3(n-3)^3}{n^3(n-1)^3}\\
I_6(f)&=&-\frac {(n-2)^6(n-3)^6}{n^6(n-1)^6}\,\frac{(q^n-1)^2\,(p^n-1)^2\,(p^n-q^n)^2
\,(p^n+q^n+1)^3}{(p^n)^3\,(q^n)^3}\\
I_7(f)&=&\frac {(n-2)(n-3)}{n(n-1)}\,\frac{(p^n)^2+q^n(p^n)^2+p^n-2\,p^n\,q^n+(q^n)^2\,p^n+q^n+(q^n)^2}{p^n\,q^n}\\
I_8(f)&=&-\frac {(n-2)^{10}(n-3)^6}{n^8(n-1)^8} \frac{(q^n-1)^2\,(p^n-1)^2\,(p^n-q^n)^2
\,(p^n+q^n+1)^3}{(p^n)^3\,(q^n)^3}.
\end{eqnarray*}
We note that $I_5(f)$ is constant and $I_8(f)$ is a constant multiple of $I_6(f)$.
We need however three other relations to define the signature manifold. We observe  that all invariants are functions of $p^n$ and $q^n$ and  denote  $P=p^n$ and  $Q=q^n$, then:
\begin{eqnarray*}
i_2 &=& -\frac{5\,Q^2+5\,P\,Q^2+5\,P^2\,Q-26\,P\,Q+5\,Q+5\,P+5\,P^2}{PQ}=
\frac{n(n-1)}{(n-2)^2}\,I_2,\\
i_3 &=& -\frac{(Q+1-P)(Q-1+P)(Q-1-P)}{PQ}=\frac{n^2(n-1)^2}{(n-2)^4}\, I_3,\\
i_4 &=& \frac{Q^3+P^3+1}{PQ}=\frac{n^2(n-1)^2}{(n-2)^2(n-3)^2}\,I_4,\\
i_6 &=& -\frac{(P-1)^2(Q-1)^2(Q-P)^2(P+Q+1)^3}{P^3Q^3}=\frac{n^6(n-1)^6}{(n-2)^6(n-3)^6}\,I_6,\\
i_7 &=&\frac{PQ^2+Q^2-2\,PQ+P^2Q+Q+P^2+P}{PQ}=\frac{n(n-1)}{(n-2)(n-3)}\,I_7,\\
i_8 &=& -\frac{(P-1)^2(Q-1)^2(Q-P)^2(P+Q+1)^3}{P^3Q^3}=\frac{n^8(n-1)^8}
{(n-2)^{10}(n-3)^6}\,I_8\\
\end{eqnarray*}

The following four relations were computed by {\sc Macaulay~2}\hfill\break
(remarkably {\sc Maple~5} was not be able to handle these computations):
\begin{eqnarray*}
i_6-i_8=0,\qquad i_3+i_4-i_7=0,\qquad i_2+5\,i_7-16=0, \\
i_4\,i_7^2-i_7^3-4\,i_4^2+4\,i_4\,i_7-8\,i_4+12\,i_7+i_8-20=0 .
\end{eqnarray*}
The corresponding relations among invariants $I_2(f)$, $I_3(f)$, $I_4(f)$, $I_6(f)$,
$I_7(f)$ and  $I_8(f)$ follow immediately.
Combining them with the conditions:
$$ I_1(f)=-\frac{(n-2)^2}{n(n-1)}\    \quad\mbox{and}  \quad I_5(f)=-\frac{(n-2)^3(n-3)^3}{n^3(n-1)^3}$$
and demanding that  $\dim\CC(f)=2$ (equivalently we demand that the symmetry group of
$f$ is discrete),  we thus obtain necessary and sufficient conditions for a ternary
$n$-form to be equivalent to the sum of $n$-th powers.

\chapter*{Appendices.}
\addcontentsline{toc}{chapter}{\protect\numberline{}{Appendices}}
\pagestyle{myheadings}
\markboth{}{}
\pagestyle{headings}
\appendix
\chapter{Computations on Binary Forms.}

\section{The {\sc Maple} Code for computing Symmetries of Binary Forms}

The \Maple\ code consists of two main programs --- {\tt symm} and {\tt matrices} --- 
and two auxiliary  functions --- {\tt simple}  and {\tt l\_f}.  
The program {\tt symm} is  the main function. 
The input consists of a complex-valued polynomial $ f(p)$  considered as the
projective version of homogeneous binary form $F(x,y)$, and
the degree  $n=deg(F)$.   The program
computes the invariants $J$ and $K$ in reduced form, determines the dimension of the symmetry
group, and, in the case of a finite symmetry group, applies the
\Maple\ command {\tt solve} to solve the two polynomial symmetry equations (\ref{EqJK}) to find
explicit form of symmetries. The
output of {\tt symm} consists of the projective index of the form and the explicit formulae for its
discrete   projective symmetries.  The program also notifies the user if the
symmetry group is not discrete,  or is in the maximal discrete symmetry class.
The program works well  when applied to very simple forms, but experienced
difficulties simplifying complicated rational algebraic formulae into the basic linear
fractional form.
\def\lone{\moveright 8pt\vbox}
\def\ltwo{\moveright 16pt\vbox}
\def\lthree{\moveright 24pt\vbox}
\def\lfour{\moveright 32 pt\vbox}
\def\lfive{\moveright 40 pt\vbox}
\def\lsix{\moveright 48 pt\vbox}
\def\mapskip{\vskip10mm}
\mapskip
{\tt > with(linalg):}

{\tt > symm:=proc(form,n) 

 global tr,error;

local Q,Qp,Qpp,Qppp,Qpppp,H,T,V,U,J,K,j,k, Eq1,Eq2,i,eqtr,

ans;

\lone{  tr:='tr':
 
  Q:=form(p);

  Qp:=diff(Q,p);

  Qpp:=diff(Qp,p);
 
  Qppp:=diff(Qpp,p);

  Qpppp:=diff(Qppp,p);

 H:=n*(n-1)*(Q*Qpp-(n-1)/n*Qp\^{}2);

 if H=0 then }

\ltwo {  ans:=`Hessian is zero: two-dimensional symmetry group`}

\lone {  else }

\ltwo  
{
T:=-n\^{}2*(n-1)*(Q\^{}2*Qppp-3*(n-2)/n*Q*Qp*Qpp

+2*(n-1)*(n-2)/n\^{}2*Qp\^{}3);

V:=Q\^{}3*Qpppp-4*(n-3)/n*Q\^{}2*Qp*Qppp+6*(n-2)*(n-3)/n\^{}2

*Q*Qp\^{}2*Qpp-3*(n-1)*(n-2)*(n-3)/n\^{}3*Qp\^{}4;

  U:=n\^{}3*(n-1)*V-3*(n-2)/(n-1)*H\^{}2;

  J:=simple(T\^{}2/H\^{}3); K:=simple(U/H\^{}2);

  j:=subs(p=P,J);k:=subs(p=P,K);

  Eq1:=simplify(numer(J)*denom(j)-numer(j)*denom(J));

  Eq2:=simplify(numer(K)*denom(k)-numer(k)*denom(K));

  if Eq1=0  then}

\lthree { ans:=`Form has a one-dimensional symmetry group`;}

\ltwo{    else }

\lthree{  if Eq2=0  then }

\lfour{ print (` Form has the maximal possible discrete 

symmetry group`);

        eqtr:= [solve(Eq1,P)];

        tr:=map(radsimp,map(allvalues,eqtr));}
     
\lthree{  else}

\lfour{ eqtr:=[solve(\{Eq1,Eq2\},P)];

        tr:= [];

        for i from 1 to nops(eqtr) do} 

\lfive  { tr:=map(radsimp,[op(tr),allvalues(rhs(eqtr[i][1]))]);}

\lfour {        od}

\lthree {     fi;

    print(`The number of elements in the symmetry group`

   =nops(tr));

    ans:=map(l\_{}f,tr);

    if error=1 then }

\lfour{ print(`ERROR: Some of the transformations are not
 
linear-fractional`) }

\lthree{    else }

\lfour{     if error=2 then }

\lfive {      print(`WARNING: Some of the transformations are not

written   in the form polynomial over polynomial`)}

\lfour{    fi;}

\lthree{    fi;}

\ltwo{   fi;}

\lone{ fi;

 ans }
end:
}

\mapskip
\noindent
The program {\tt matrices} determines the matrix
symmetry corresponding to a given (list of) projective symmetries.  As discussed in the text, this only
requires determining an overall scalar multiple, which can be found by substituting the projective
symmetry into the form.  The output consists of each projective symmetry, the scalar factor $\mu $,
and the resulting matrix symmetry.

\mapskip
{\tt 
 > matrices:=proc(form,n,L::list)
 
  local Q,ks,ksi,i,Sf,M;
 
 \lone {ksi:='ksi'; 

  for i from 1 to nops(L) do }

  \ltwo {Sf:=simplify(denom(L[i])\^{}n*form(L[i]));
 
  ks:=quo(Sf,form(p),p); 

  ksi:=simplify(ks\^{}(1/n),radical,symbolic); 

  M[i]:=matrix(2,2,[coeff(numer(L[i]),p)/ksi, 

  coeff(numer(L[i]),p,0)/ksi,coeff(denom(L[i]),p)/ksi, 

  coeff(denom(L[i]),p,0)/ksi]); 

  print(L[i],    mu=ksi,   map(simplify,M[i]))}

\lone { od;}

end:
}

\mapskip
\noindent
The auxiliary function {\tt simple}  helps to simplify rational
 expressions by manipulating the numerator and denominator separately.
The simplified rational expression is returned.

\mapskip

{\tt 
> simple:=proc(x) 

local nu,de,num,den;

\lone{ nu:=numer(x);

de:=denom(x);

num:=(simplify((nu,radical,symbolic))); 

den:=(simplify((de,radical,symbolic)));

simplify(num/den);}

end:
}

\mapskip
\noindent
The auxiliary function {\tt l\_f} uses polynomial division  to
reduce rational expressions to linear fractional form (when possible).

\mapskip

{ \tt 
> l\_{}f:=proc(x) 

 local A,B,C,S,de,nu,r,R;

 global error;error:='error';

\lone { nu:=numer(x);  

 de:=denom(x);

 if type(nu,polynom(anything,p)) 

and type(de,polynom(anything,p)) then}
 
\ltwo{  if degree(nu,p)+1=degree(de,p) then} 

\lthree{   A:=quo(de,nu,p,'B'); 

   S:=1/A; R:=0} 

\ltwo { else }
 
\lthree { A:=quo(nu,de,p,'B');

   if B=0 then }

\lfour {    S:=A;  R:=0;  }

\lthree{   else }

\lfour {     C:=quo(de,B,p,'r'); R:=simple(r);

     S:=simplify(A+1/C) }
  
\lthree{  fi;}
 
\ltwo{ fi;

  if R=0 then }  

\lthree  {  collect(S,p)}

\ltwo {  else }

\lthree  {  error:=1; x }

\ltwo { fi;}

\lone { else }

\ltwo {  error:=2; x } 
 
\lone {fi;}

end:}

\section{ Cubic Forms.}

We now present the results of applying the function {\tt symm} and {\tt matrices} to cubic forms.
We begin with simple cases, ending with a cubic whose formulae required extensive manipulation.

\medskip

\noindent  {\bf Cubics with  one triple root:}

\smallskip
\begin{maplegroup}
\begin{mapleinput}
\mapleinline{active}{1d}{f:=p->p^3;}{%
}
\end{mapleinput}

\mapleresult
\begin{maplelatex}
\[
f := p\rightarrow p^{3}
\]
\end{maplelatex}

\end{maplegroup}

\begin{maplegroup}
\begin{mapleinput}
\mapleinline{active}{1d}{symm(f,3);}{%
}
\end{mapleinput}

\mapleresult
\begin{maplelatex}
\[
\mathit{Hessian\ is\ zero:\ two-dimensional\ symmetry\ group}
\]
\end{maplelatex}

\end{maplegroup}

\smallskip

\noindent  {\bf Cubics with one double root and one single root:}

\smallskip
\begin{maplegroup}
\begin{mapleinput}
\mapleinline{active}{1d}{f:=p->p;}{%
}
\end{mapleinput}

\mapleresult
\begin{maplelatex}
\[
f := p\rightarrow p
\]
\end{maplelatex}

\end{maplegroup}
\begin{maplegroup}
\begin{mapleinput}
\mapleinline{active}{1d}{symm(f,3);}{%
}
\end{mapleinput}

\mapleresult
\begin{maplelatex}
\[
\mathit{Form\ has\ a\ one-dimensional\ symmetry\ group}
\]
\end{maplelatex}

\end{maplegroup}

\smallskip

\noindent  {\bf Cubics with three simple roots:}

\smallskip
\begin{maplegroup}
\begin{mapleinput}
\mapleinline{active}{1d}{f:=p->p^3+1;}{%
}
\end{mapleinput}

\mapleresult
\begin{maplelatex}
\[
f := p\rightarrow p^{3} + 1
\]
\end{maplelatex}

\end{maplegroup}

\begin{maplegroup}
\begin{mapleinput}
\mapleinline{active}{1d}{S:=symm(f,3);}{%
}
\end{mapleinput}

\mapleresult
\begin{maplelatex}
\[
\mathit{\ Form\ has\ the\ maximal\ possible\ discrete\ symmetry\ 
group}
\]
\end{maplelatex}

\begin{maplelatex}
\[
\mathit{The\ number\ of\ elements\ in\ the\ symmetry\ group}=6
\]
\end{maplelatex}

\begin{maplelatex}
\[
S :=  \left[  \; p, \,{\displaystyle \frac {1}{p}} , \,
{\displaystyle \frac { - {\displaystyle \frac {1}{2}}  + 
{\displaystyle \frac {1}{2}} \,I\,\sqrt{3}}{p}} , \,
{\displaystyle \frac { - {\displaystyle \frac {1}{2}}  - 
{\displaystyle \frac {1}{2}} \,I\,\sqrt{3}}{p}} , \,( - 
{\displaystyle \frac {1}{2}}  + {\displaystyle \frac {1}{2}} \,I
\,\sqrt{3})\,p, \,( - {\displaystyle \frac {1}{2}}  - 
{\displaystyle \frac {1}{2}} \,I\,\sqrt{3})\,p \;  \right] 
\]
\end{maplelatex}

\end{maplegroup}

\begin{maplegroup}
\begin{mapleinput}
\mapleinline{active}{1d}{matrices(f,3,[S[2],S[4]]);}{%
}
\end{mapleinput}

\mapleresult
\begin{maplelatex}
\[
{\displaystyle \frac {1}{p}} , \quad \mu =1, \quad \left[ 
{\begin{array}{rr}
0 & 1 \\
1 & 0
\end{array}}
 \right] 
\]
\end{maplelatex}

\begin{maplelatex}
\[
{\displaystyle \frac { - {\displaystyle \frac {1}{2}}  - 
{\displaystyle \frac {1}{2}} \,I\,\sqrt{3}}{p}} , \quad \mu =2, \quad
 \left[ 
{\begin{array}{rc}
0 &  - {\displaystyle \frac {1}{2}}  - {\displaystyle \frac {1}{2
}} \,I\,\sqrt{3} \\ [2ex]
1 & 0
\end{array}}
 \right] 
\]
\end{maplelatex}

\end{maplegroup}

\smallskip

\noindent {\bf  A more complicated cubic example.}

\smallskip

All cubics with a discrete symmetry group are complex
equivalent to $x^3+y^3$ and have projective index $6$. 
However, when we apply the same code to a cubic not in canonical form. The initial \Maple\ result
is not in the correct linear fractional form.  We must simplify the rational algebraic expressions ``by
hand'' to put them in the form of a projective linear fractional transformation. 

\mapskip
\begin{maplegroup}
\begin{mapleinput}
\mapleinline{active}{1d}{f:=p->p^3+p+1;}{%
}
\end{mapleinput}

\mapleresult
\begin{maplelatex}
\[
f := p\rightarrow p^{3} + p + 1
\]
\end{maplelatex}

\end{maplegroup}
\begin{maplegroup}
\begin{mapleinput}
\mapleinline{active}{1d}{S:=symm(f,3);}{%
}
\end{mapleinput}

\mapleresult
\begin{maplelatex}
\[
\mathit{\ Form\ has\ the\ maximal\ possible\ discrete\ symmetry\ 
group}
\]
\end{maplelatex}

\begin{maplelatex}
\[
\mathit{The\ number\ of\ elements\ in\ the\ symmetry\ group}=6
\]
\end{maplelatex}

\begin{maplelatex}
\begin{eqnarray*}
& &\mathit{WARNING:\ Some\ of\ the\ transformations\ are\ 
not\ written\  in\ the\ form}\\
& & \mathit{polynomial\ over\ polynomial}
\end{eqnarray*}
\end{maplelatex}

\begin{maplelatex}
\begin{eqnarray*}
& & S := [p, \,{\displaystyle \frac {( - 9 + I\,\sqrt{31})\,
p + 2}{9+ I\,\sqrt{31} + 6\,p}} , \, - {\displaystyle \frac {(9
 + I\,\sqrt{31})\,p - 2}{9 - I\,\sqrt{31} + 6\,p}} ,\\
& & 
{\displaystyle \frac {1}{18}} ((54\,p^{4} + 9\,2^{(2/3)}\,3^{(1/3)}\,
\mathrm{\%1}^{(1/3)}\,p^{3} + 324\,p^{3} + 3\,2^{(1/3)}\,3^{(2/3)}\,\mathrm{\%1}
^{(2/3)}\,p^{2} + 450\,p^{2}\\
& & - 108\,p + 9\,2^{(1/3)}\,3^{(2/3)}\,
\mathrm{\%1}^{(2/3)}\,p  + 9\,2^{(2/3)}\,3^{(1/3)}\,\mathrm{\%1}^{(1/3)}\,p
 - 2^{(1/3)}\,3^{(2/3)}\,\mathrm{\%1}^{(2/3)}\\
& & + 6 + 9\,2^{(2/3)}
\,3^{(1/3)}\,\mathrm{\%1}^{(1/3)})3^{(2/3)} 
2^{(1/3)}) \left/{\vrule height0.56em width0em depth0.56em}
 \right. \!  \!\\
& & (\mathrm{\%1}^{(1/3)}\,(27\,p^{3} - 9\,p^{2} - 1)
),  - {\displaystyle \frac {1}{36}} ((54\,I\,\sqrt{3}\,p^{4} + 54
\,p^{4} + 324\,p^{3} \\
& & \mbox{} + 324\,I\,\sqrt{3}\,p^{3} - 18\,2^{(2/3)}\,3^{(1/3)}
\,\mathrm{\%1}^{(1/3)}\,p^{3} + 450\,p^{2} + 450\,I\,\sqrt{3}\,p^{2}
\end{eqnarray*}
\end{maplelatex}

\begin{maplelatex}
\begin{eqnarray*}
 & & \mbox{} - 9\,I\,3^{(1/6)}\,2^{(1/3)}\,\mathrm{\%1}^{(2/3)}\,
p^{2} + 3\,2^{(1/3)}\,3^{(2/3)}\,\mathrm{\%1}^{(2/3)}\,p^{2} - 27
\,I\,3^{(1/6)}\,2^{(1/3)}\,\mathrm{\%1}^{(2/3)}\,p \\
 & & \mbox{} - 108\,p - 18\,2^{(2/3)}\,3^{(1/3)}\,\mathrm{\%1}^{(
1/3)}\,p + 9\,2^{(1/3)}\,3^{(2/3)}\,\mathrm{\%1}^{(2/3)}\,p - 108
\,I\,\sqrt{3}\,p \\
 & & \mbox{} + 3\,I\,3^{(1/6)}\,2^{(1/3)}\,\mathrm{\%1}^{(2/3)}
 + 6\,I\,\sqrt{3} - 18\,2^{(2/3)}\,3^{(1/3)}\,\mathrm{\%1}^{(1/3)
} - 2^{(1/3)}\,3^{(2/3)}\,\mathrm{\%1}^{(2/3)} + 6 \\
 & & )3^{(2/3)}\,2^{(1/3)}) \left/ {\vrule 
height0.56em width0em depth0.56em} \right. \!  \! (\mathrm{\%1}^{
(1/3)}\,(27\,p^{3} - 9\,p^{2} - 1)), {\displaystyle \frac {1}{36}
} ((54\,I\,\sqrt{3}\,p^{4} - 54\,p^{4} - 324\,p^{3} \\
 & & \mbox{} + 324\,I\,\sqrt{3}\,p^{3} + 18\,2^{(2/3)}\,3^{(1/3)}
\,\mathrm{\%1}^{(1/3)}\,p^{3} - 450\,p^{2} + 450\,I\,\sqrt{3}\,p
^{2} \\
 & & \mbox{} - 9\,I\,3^{(1/6)}\,2^{(1/3)}\,\mathrm{\%1}^{(2/3)}\,
p^{2} - 3\,2^{(1/3)}\,3^{(2/3)}\,\mathrm{\%1}^{(2/3)}\,p^{2} - 27
\,I\,3^{(1/6)}\,2^{(1/3)}\,\mathrm{\%1}^{(2/3)}\,p \\
 & & \mbox{} + 108\,p + 18\,2^{(2/3)}\,3^{(1/3)}\,\mathrm{\%1}^{(
1/3)}\,p - 9\,2^{(1/3)}\,3^{(2/3)}\,\mathrm{\%1}^{(2/3)}\,p - 108
\,I\,\sqrt{3}\,p \\
 & & \mbox{} + 3\,I\,3^{(1/6)}\,2^{(1/3)}\,\mathrm{\%1}^{(2/3)}
 + 6\,I\,\sqrt{3} + 18\,2^{(2/3)}\,3^{(1/3)}\,\mathrm{\%1}^{(1/3)
} + 2^{(1/3)}\,3^{(2/3)}\,\mathrm{\%1}^{(2/3)}\\
& & - 6)3^{(2/3)}\,2^{(1/3)}) \left/ {\vrule 
height0.56em width0em depth0.56em} \right. \!  \! (\mathrm{\%1}^{
(1/3)}\,(27\,p^{3} - 9\,p^{2} - 1))] \\
 & & \mathrm{\%1} := 9 + 18\,p - 81\,p^{2} + 261\,p^{3} + 27\,
\sqrt{31}\,\sqrt{3}\,p^{3} - 9\,\sqrt{31}\,\sqrt{3}\,p^{2} - 
\sqrt{31}\,\sqrt{3}
\end{eqnarray*}
\end{maplelatex}

\end{maplegroup}

\mapskip
\begin{maplegroup}
The first three components of {\tt S} are in the proper linear fractional form.  The problem with the
other expressions is that \Maple\ does not automatically factor polynomials under a radical.
One approach to simplification is to first do the required factorization:

\end{maplegroup}
\mapskip
\begin{maplegroup}
\begin{mapleinput}
\mapleinline{active}{1d}{
n1:=factor(9+18*p-81*p^2+261*p^3+27*sqrt(31)*sqrt(3)*p^3\break 
-9*sqrt(31)*sqrt(3)*p^2-sqrt(31)*sqrt(3));
}{%
}
\end{mapleinput}

\mapleresult
\begin{maplelatex}
\[
\mathit{n1} :=  - {\displaystyle \frac {1}{24}} \,(29 + 3\,\sqrt{
31}\,\sqrt{3})\,( - 6\,p - 9 + \sqrt{31}\,\sqrt{3})^{3}
\]
\end{maplelatex}

\end{maplegroup}
\mapskip
\begin{maplegroup}
Substituting {\tt n1} into the fourth rational algebraic expression in {\tt S} above,
we can now force \Maple\ to take the cube root and obtain the actual linear fractional
formula for this symmetry:

\end{maplegroup}
\mapskip
\begin{maplegroup}

\begin{mapleinput}
\mapleinline{active}{1d}{
simp1:=radsimp(1/18*((54*p^4+9*2^(2/3)*3^(1/3)*(n1)^(1/3)*p^3\break
+324*p^3+3*2^(1/3)*3^(2/3)*(n1)^(2/3)*p^2+450*p^2\break
+9*2^(1/3)*3^(2/3)*(n1)^(2/3)*p+9*2^(2/3)*3^(1/3)*(n1)^(1/3)*p\break
-108*p-2^(1/3)*3^(2/3)*(n1)^(2/3)+6\break
+9*2^(2/3)*3^(1/3)*(n1)^(1/3))*3^(2/3)*2^(1/3))/((n1)^(1/3)\break
*(27*p^3-9*p^2-1)));
}{%
}

\end{mapleinput}
\mapleresult
\begin{maplelatex}
\end{maplelatex}
\end{maplegroup}

\begin{maplegroup}
\begin{mapleinput}
\mapleinline{active}{1d}{simp2:=l_f(simp1);}{%
}

\end{mapleinput}
\mapleresult
\begin{maplelatex}
\end{maplelatex}
\end{maplegroup}

\begin{maplegroup}

\begin{mapleinput}
\mapleinline{active}{1d}
{simp3:=collect(expand(numer(simp2))/expand(denom(simp2)),p);}{%
}

\end{mapleinput}

\mapleresult
\begin{maplelatex}
\begin{eqnarray*}
\lefteqn{\mathit{simp3} := (( - 226\,2^{(1/3)} - 20\,\mathrm{\%2}
 - 20\,2^{(2/3)}\,\mathrm{\%1}^{(2/3)} - 22\,\mathrm{\%3} - 208\,
\mathrm{\%1}^{(1/3)})\,p - 8\,\mathrm{\%3}} \\
 & & \mbox{} - 58\,2^{(2/3)}\,\mathrm{\%1}^{(2/3)} - 566\,
\mathrm{\%1}^{(1/3)} - 58\,\mathrm{\%2} - 116\,2^{(1/3)}) \left/ 
{\vrule height0.56em width0em depth0.56em} \right. \!  \! ( \\
 & & ( - 6\,2^{(2/3)}\,\mathrm{\%1}^{(2/3)} - 174\,\mathrm{\%1}^{
(1/3)} - 1686\,2^{(1/3)} - 174\,\mathrm{\%3} - 6\,\mathrm{\%2})\,
p \\
 & & \mbox{} + 20\,2^{(2/3)}\,\mathrm{\%1}^{(2/3)} + 208\,
\mathrm{\%1}^{(1/3)} + 226\,2^{(1/3)} + 22\,\mathrm{\%3} + 20\,
\mathrm{\%2}) \\
 & & \mathrm{\%1} := 29 + 3\,\sqrt{31}\,\sqrt{3} \\
 & & \mathrm{\%2} := \mathrm{\%1}^{(1/3)}\,\sqrt{31}\,\sqrt{3} \\
 & & \mathrm{\%3} := 2^{(1/3)}\,\sqrt{31}\,\sqrt{3}
\mbox{\hspace{303pt}}
\end{eqnarray*}
\end{maplelatex}

\end{maplegroup}

\mapskip

The linear fractional formulae for the other symmetries are derived in a similar fashion.

\section { The Octahedral Symmetry Group.}

As we remarked in the text, the sextic polynomial $Q(p) = p^5 + p$ has an octahedral symmetry group.
Here we illustrate how the symmetries are computed using our \Maple\ program.

\mapskip

\begin{maplegroup}
\begin{mapleinput}
\mapleinline{active}{1d}{f:=p->p^5+p;
}{%
}
\end{mapleinput}

\mapleresult
\begin{maplelatex}
\[
f := p\rightarrow p^{5} + p
\]
\end{maplelatex}

\end{maplegroup}
\begin{maplegroup}
\begin{mapleinput}
\mapleinline{active}{1d}{symm(f,6);}{%
}
\end{mapleinput}

\mapleresult
\begin{maplelatex}
\[
\mathit{\ Form\ has\ the\ maximal\ possible\ discrete\ symmetry\ 
group}
\]
\end{maplelatex}

\begin{maplelatex}
\[
\mathit{The\ number\ of\ elements\ in\ the\ symmetry\ group}=24
\]
\end{maplelatex}

\begin{maplelatex}
\begin{eqnarray*}
\lefteqn{\mathit{WARNING:\ Some\ of\ the\ transformations\ are\ 
not\ written\ in\ the\ form
}} \\
 & & \mathit{polynomial\ over\ polynomial}\mbox{\hspace{246pt}}
\end{eqnarray*}
\end{maplelatex}

\begin{maplelatex}
\begin{eqnarray*}
\lefteqn{\left [\ {\displaystyle \frac {1}{p}} , \,p, \, - p, \, - 
{\displaystyle \frac {1}{p}} , \,{\displaystyle \frac {I}{p}} , 
\, - {\displaystyle \frac {I}{p}} , \,I\,p, \, - I\,p, \,
{\displaystyle \frac { - 2\,p^{3} + 2\,I\,p + \mathrm{\%3}}{p^{4}
 + 1}} , \,{\displaystyle \frac { - 2\,p^{3} + 2\,I\,p - \mathrm{
\%3}}{p^{4} + 1}} ,\right . } \\
 & & {\displaystyle \frac { - 2\,p^{3} - 2\,I\,p + \mathrm{\%4}}{
p^{4} + 1}} , \,{\displaystyle \frac { - 2\,p^{3} - 2\,I\,p - 
\mathrm{\%4}}{p^{4} + 1}} , \,{\displaystyle \frac {2\,p^{3} + 2
\,I\,p + \mathrm{\%4}}{p^{4} + 1}} , \,{\displaystyle \frac {2\,p
^{3} + 2\,I\,p - \mathrm{\%4}}{p^{4} + 1}} ,  \\
 & & {\displaystyle \frac {2\,p^{3} - 2\,I\,p + \mathrm{\%3}}{p^{
4} + 1}} , \,{\displaystyle \frac {2\,p^{3} - 2\,I\,p - \mathrm{
\%3}}{p^{4} + 1}} , \,{\displaystyle \frac { - 2\,p + 2\,I\,p^{3}
 + \mathrm{\%1}}{p^{4} + 1}} , \,{\displaystyle \frac { - 2\,p + 
2\,I\,p^{3} - \mathrm{\%1}}{p^{4} + 1}} ,  \\
 & & {\displaystyle \frac { - 2\,p - 2\,I\,p^{3} + \mathrm{\%2}}{
p^{4} + 1}} , \,{\displaystyle \frac { - 2\,p - 2\,I\,p^{3} - 
\mathrm{\%2}}{p^{4} + 1}} , \,{\displaystyle \frac {2\,p + 2\,I\,
p^{3} + \mathrm{\%2}}{p^{4} + 1}} , \,{\displaystyle \frac {2\,p
 + 2\,I\,p^{3} - \mathrm{\%2}}{p^{4} + 1}} ,  \\
 & & \left .{\displaystyle \frac {2\,p - 2\,I\,p^{3} + \mathrm{\%1}}{p^{
4} + 1}} , \,{\displaystyle \frac {2\,p - 2\,I\,p^{3} - \mathrm{
\%1}}{p^{4} + 1}} \ \right ] \\
 & & \mathrm{\%1} := \sqrt{ - 4\,p^{6} + 4\,p^{2} + I\,p^{8} - 6
\,I\,p^{4} + I} \\
 & & \mathrm{\%2} := \sqrt{ - 4\,p^{6} + 4\,p^{2} - I\,p^{8} + 6
\,I\,p^{4} - I} \\
 & & \mathrm{\%3} := \sqrt{4\,p^{6} - 4\,p^{2} + I\,p^{8} - 6\,I
\,p^{4} + I} \\
 & & \mathrm{\%4} := \sqrt{4\,p^{6} - 4\,p^{2} - I\,p^{8} + 6\,I
\,p^{4} - I}
\end{eqnarray*}
\end{maplelatex}

\end{maplegroup}
\mapskip
\begin{maplegroup}
Again, \Maple\ has failed to simplify the expressions $\%1,\%2,\%3,\%4$, and we need to make it
take the square  root.  In the case of symmetries numbers $9,11,13,15,17,19$, $21,23$ 
this is done as
follows.  The others are handled in a similar fashion, and, for brevity, we omit the formulae here.

\end{maplegroup}
\mapskip
\begin{maplegroup}
\begin{mapleinput}
\mapleinline{active}{1d}{
for j in [9,11,13,15,17,19,21,23] do \break
sq:=sqrt(factor(op(op(numer(tr[j]))[3])[1],I),symbolic):\break
s[j]:=\break
l_f((op(numer(tr[j]))[1]+op(numer(tr[j]))[2]+sq)/denom(tr[j]));\break
 print(s.j=s[j]);\break      
 od:  }{%
}
\end{mapleinput}

\mapleresult
\begin{maplelatex}
\[
\mathit{s9}= - {\displaystyle \frac {(\sqrt{2} + I\,\sqrt{2})\,p
 - 2}{ - \sqrt{2} + I\,\sqrt{2} - 2\,p}} 
\]
\end{maplelatex}

\begin{maplelatex}
\[
\mathit{s11}= - {\displaystyle \frac {( - \sqrt{2} + I\,\sqrt{2})
\,p + 2}{I\,\sqrt{2} + \sqrt{2} + 2\,p}} 
\]
\end{maplelatex}

\begin{maplelatex}
\[
\mathit{s13}={\displaystyle \frac {( - \sqrt{2} + I\,\sqrt{2})\,p
 - 2}{I\,\sqrt{2} + \sqrt{2} - 2\,p}} 
\]
\end{maplelatex}

\begin{maplelatex}
\[
\mathit{s15}={\displaystyle \frac {(\sqrt{2} + I\,\sqrt{2})\,p + 
2}{ - \sqrt{2} + I\,\sqrt{2} + 2\,p}} 
\]
\end{maplelatex}

\begin{maplelatex}
\[
\mathit{s17}= - {\displaystyle \frac {( - \sqrt{2} + I\,\sqrt{2})
\,p - 2}{ - \sqrt{2} + I\,\sqrt{2} - 2\,I\,p}} 
\]
\end{maplelatex}

\begin{maplelatex}
\[
\mathit{s19}= - {\displaystyle \frac {I\,((\sqrt{2} + I\,\sqrt{2}
)\,p + 2)}{ - \sqrt{2} + I\,\sqrt{2} + 2\,p}} 
\]
\end{maplelatex}

\begin{maplelatex}
\[
\mathit{s21}={\displaystyle \frac {(\sqrt{2} + I\,\sqrt{2})\,p - 
2}{\sqrt{2} + I\,\sqrt{2} + 2\,I\,p}} 
\]
\end{maplelatex}

\begin{maplelatex}
\[
\mathit{s23}={\displaystyle \frac {( - \sqrt{2} + I\,\sqrt{2})\,p
 + 2}{ - \sqrt{2} + I\,\sqrt{2} + 2\,I\,p}} 
\]
\end{maplelatex}

\end{maplegroup}

\mapskip
\begin{maplegroup}
As we remarked in the text, the octahedral symmetry group has two generators.
The matrix form of these generators is computed as follows:

\end{maplegroup}
\begin{maplegroup}
\begin{mapleinput}
\mapleinline{active}{1d}{matrices(f,6,[tr[7],s[9]]);}{%
}
\end{mapleinput}

\mapleresult
\begin{maplelatex}
\[
I\,p, \qquad \mu =(-1)^{(1/12)}, \qquad  \left[ 
{\begin{array}{cc}
(-1)^{(5/12)} & 0 \\
0 &  - (-1)^{(11/12)}
\end{array}}
 \right] 
\]
\end{maplelatex}

\begin{maplelatex}
\[
 - {\displaystyle \frac {(\sqrt{2} + I\,\sqrt{2})\,p - 2}{ - 
\sqrt{2} + I\,\sqrt{2} - 2\,p}} , \qquad \mu =2\,\sqrt{2}, \qquad  \left[ 
{\begin{array}{cc}
 - {\displaystyle \frac {1}{2}}  - {\displaystyle \frac {1}{2}} 
\,I & {\displaystyle \frac {1}{2}} \,\sqrt{2} \\ [2ex]
 - {\displaystyle \frac {1}{2}} \,\sqrt{2} &  - {\displaystyle 
\frac {1}{2}}  + {\displaystyle \frac {1}{2}} \,I
\end{array}}
 \right] 
\]
\end{maplelatex}

\end{maplegroup}
\mapskip
\begin{maplegroup}
We end with two further examples.  We already know that the following octavic polynomial also has an 
octahedral symmetry group.  In this case, {\tt symm} produces the projective symmetries directly:

\end{maplegroup}
\mapskip
\begin{maplegroup}
\begin{mapleinput}
\mapleinline{active}{1d}{f:=p->p^8+14*p^4+1;
}{%
}
\end{mapleinput}

\mapleresult
\begin{maplelatex}
\[
f := p\rightarrow p^{8} + 14\,p^{4} + 1
\]
\end{maplelatex}

\end{maplegroup}
\begin{maplegroup}
\begin{mapleinput}
\mapleinline{active}{1d}{S:=symm(f,8);}{%
}
\end{mapleinput}

\mapleresult
\begin{maplelatex}
\[
\mathit{The\ number\ of\ elements\ in\ the\ symmetry\ group}=24
\]
\end{maplelatex}

\begin{maplelatex}
\begin{eqnarray*}
&&S := [ - {\displaystyle \frac {1}{p}} , \, - 
{\displaystyle \frac {p - 1}{p + 1}} , \, - {\displaystyle 
\frac {p + 1}{p - 1}} , \,p, \, - p, \,{\displaystyle \frac {p + 
1}{p - 1}} , \,{\displaystyle \frac {p - 1}{p + 1}} , \,
{\displaystyle \frac {1}{p}} , \,{\displaystyle \frac {I\,(p - 1)
}{p + 1}} , \, \\
& &- {\displaystyle \frac {I\,(p - 1)}{p + 1}} , \,
{\displaystyle \frac {I\,(p + 1)}{p - 1}} ,\\  
 & &  - {\displaystyle \frac {I\,(p + 1)}{p - 1}} , \,
{\displaystyle \frac {I}{p}} , \, - {\displaystyle \frac {I}{p}} 
, \,I\,p, \, - I\,p, \,{\displaystyle \frac { - 1 + I\,p}{ - I + 
p}} , \, - {\displaystyle \frac {1 + I\,p}{I + p}} , \\
& &{\displaystyle \frac {1 + I\,p}{I + p}} , \, - {\displaystyle 
\frac { - 1 + I\,p}{ - I + p}} , \,{\displaystyle \frac { - 1 + I
\,p}{1 + I\,p}} , \, 
  {\displaystyle \frac {1 + I\,p}{ - 1 + I\,p}} , \, - 
{\displaystyle \frac {1 + I\,p}{ - 1 + I\,p}} , \, - 
{\displaystyle \frac { - 1 + I\,p}{1 + I\,p}} ]
\end{eqnarray*}
\end{maplelatex}

\end{maplegroup}
\begin{maplegroup}
\begin{mapleinput}
\mapleinline{active}{1d}{matrices(f,8,[S[11],S[15]]);}{%
}
\end{mapleinput}

\mapleresult
\begin{maplelatex}
\[
{\displaystyle \frac {I\,(p + 1)}{p - 1}} , \qquad \mu =\sqrt{2}, \qquad
 \left[ 
{\begin{array}{cc}
{\displaystyle \frac {1}{2}} \,I\,\sqrt{2} & {\displaystyle 
\frac {1}{2}} \,I\,\sqrt{2} \\ [2ex]
{\displaystyle \frac {1}{2}} \,\sqrt{2} &  - {\displaystyle 
\frac {1}{2}} \,\sqrt{2}
\end{array}}
 \right] 
\]
\end{maplelatex}

\begin{maplelatex}
\[
I\,p, \qquad \mu =1, \qquad \left[ 
{\begin{array}{cr}
I & 0 \\
0 & 1
\end{array}}
 \right] 
\]
\end{maplelatex}

\end{maplegroup}
\mapskip
\begin{maplegroup}
Finally, for illustrative purposes, we present a higher order example given by a binary form of degree
$12$.

\end{maplegroup}
\mapskip
\begin{maplegroup}
\begin{mapleinput}
\mapleinline{active}{1d}{f:=p->p^12-33*p^8-33*p^4+1;}{%
}
\end{mapleinput}

\mapleresult
\begin{maplelatex}
\[
f := p\rightarrow p^{12} - 33\,p^{8} - 33\,p^{4} + 1
\]
\end{maplelatex}

\end{maplegroup}
\begin{maplegroup}
\begin{mapleinput}
\mapleinline{active}{1d}{S:=symm(f,12);}{%
}
\end{mapleinput}

\mapleresult
\begin{maplelatex}
\[
\mathit{The\ number\ of\ elements\ in\ the\ symmetry\ group}=24
\]
\end{maplelatex}

\begin{maplelatex}
\begin{eqnarray*}
\lefteqn{S := [p, \, - p, \, - {\displaystyle \frac {1}{p}} , \,
{\displaystyle \frac {1}{p}} , \, - {\displaystyle \frac {p - 1}{
p + 1}} , \, - {\displaystyle \frac {p + 1}{p - 1}} , \,
{\displaystyle \frac {p + 1}{p - 1}} , \,{\displaystyle \frac {p
 - 1}{p + 1}} , \,{\displaystyle \frac {I}{p}} , \, - 
{\displaystyle \frac {I}{p}} , \,I\,p, \, - I\,p, \,
{\displaystyle \frac { - 1 + I\,p}{ - I + p}}} , \\
 & &- {\displaystyle \frac {1 + I\,p}{I + p}} , \,
 {\displaystyle \frac {1 + I\,p}{I + p}} , \, 
- {\displaystyle \frac { - 1 + I\,p}{ - I + p}} , \,{\displaystyle 
\frac {I\,(p - 1)}{p + 1}} , \, - {\displaystyle \frac {I\,(p - 1
)}{p + 1}} , \,{\displaystyle \frac {I\,(p + 1)}{p - 1}} , \, - 
{\displaystyle \frac {I\,(p + 1)}{p - 1}} , \\
& &
{\displaystyle \frac { - 1 + I\,p}{1 + I\,p}} , \,
{\displaystyle \frac {1 + I\,p}{ - 1 + I\,p}} ,  \,
 - {\displaystyle \frac {1 + I\,p}{ - 1 + I\,p}} , \, - 
{\displaystyle \frac { - 1 + I\,p}{1 + I\,p}} ]
\end{eqnarray*}
\end{maplelatex}

\end{maplegroup}
\begin{maplegroup}
\begin{mapleinput}
\mapleinline{active}{1d}{matrices(f,12,[S[11],S[19]]);}{%
}
\end{mapleinput}

\mapleresult
\begin{maplelatex}
\[
I\,p, \qquad \mu =1, \qquad  \left[ 
{\begin{array}{cr}
I & 0 \\
0 & 1
\end{array}}
 \right] 
\]
\end{maplelatex}

\begin{maplelatex}
\[
{\displaystyle \frac {I\,(p + 1)}{p - 1}} , \qquad \mu =(-1)^{(1/12)}
 \sqrt{2}, \qquad \left[ 
{\begin{array}{cc}
{\displaystyle \frac {1}{2}} \,(-1)^{(5/12)}\,\sqrt{2} & 
{\displaystyle \frac {1}{2}} \,(-1)^{(5/12)}\,\sqrt{2} \\ [2ex]
 - {\displaystyle \frac {1}{2}} \,(-1)^{(11/12)}\,\sqrt{
2} & {\displaystyle \frac {1}{2}} \,(-1)^{(11/12)}\,
\sqrt{2}
\end{array}}
 \right] 
\]
\end{maplelatex}

\end{maplegroup}
\chapter{Invariants for Ternary Forms.}

In this appendix we list joint relative invariants of the forms (\ref{P2}--\ref{P3})
of  Section~\ref{ternary}, obtained by omega process,
that  we subsequently used  to obtain a fundamental set of differential invariants for ternary forms: 
\begin{small}
\begin{eqnarray*}
&M_1=\frac 1 {288}(H_3^2,P_2)^{(2)}={Q_{3, \,0}}\,{Q_{1, \,2}}\,{Q_{0, \,2}}
 - {Q_{2, \,1}}^{2}\,{Q_{0, \,2}} - {Q_{3, \,0}}\,{Q_{0, \,3}}\,{Q_{1, \,1}} + 
{Q_{2, \,1}}\,{Q_{1, \,2}}\,{Q_{1, \,1}}\\
& + {Q_{2, \,1}}\,{Q_{0, \,3}}\,{Q_{2, \,0}} - {Q_{1, \,2}}^{2}\,{Q_{2, \,0}}
\end{eqnarray*}
of weight 4.
\begin{eqnarray*}
&&M_2=\frac 1 {103680}(P_3^2,P_2^3)^{(6)}=\\
& &\lefteqn{5\,{Q_{3, \,0}}^{2}\,{Q_{0, \,2}}^{3} - 30\,{Q_{3, \,0}}
\,{Q_{0, \,2}}^{2}\,{Q_{2, \,1}}\,{Q_{1, \,1}} + 24\,{Q_{3, \,0}}
\,{Q_{1, \,1}}^{2}\,{Q_{0, \,2}}\,{Q_{1, \,2}} + 36\,{Q_{2, \,1}}
^{2}\,{Q_{1, \,1}}^{2}\,{Q_{0, \,2}}} \\
 & & \mbox{} + 9\,{Q_{2, \,1}}^{2}\,{Q_{0, \,2}}^{2}\,{Q_{2, \,0}
} - 4\,{Q_{3, \,0}}\,{Q_{1, \,1}}^{3}\,{Q_{0, \,3}} - 36\,{Q_{2, 
\,1}}\,{Q_{1, \,1}}^{3}\,{Q_{1, \,2}} \\
 & & \mbox{} + 36\,{Q_{1, \,2}}^{2}\,{Q_{2, \,0}}\,{Q_{1, \,1}}^{
2} + 9\,{Q_{1, \,2}}^{2}\,{Q_{2, \,0}}^{2}\,{Q_{0, \,2}} + 6\,{Q
_{3, \,0}}\,{Q_{0, \,2}}^{2}\,{Q_{1, \,2}}\,{Q_{2, \,0}} \\
 & & \mbox{} - 6\,{Q_{3, \,0}}\,{Q_{2, \,0}}\,{Q_{0, \,2}}\,{Q_{0
, \,3}}\,{Q_{1, \,1}} - 54\,{Q_{2, \,1}}\,{Q_{2, \,0}}\,{Q_{0, \,
2}}\,{Q_{1, \,2}}\,{Q_{1, \,1}} \\
 & & \mbox{} + 24\,{Q_{2, \,1}}\,{Q_{2, \,0}}\,{Q_{1, \,1}}^{2}\,
{Q_{0, \,3}} + 6\,{Q_{2, \,1}}\,{Q_{2, \,0}}^{2}\,{Q_{0, \,2}}\,{
Q_{0, \,3}} - 30\,{Q_{1, \,2}}\,{Q_{2, \,0}}^{2}\,{Q_{0, \,3}}\,{
Q_{1, \,1}} \\
 & & \mbox{} + 5\,{Q_{0, \,3}}^{2}\,{Q_{2, \,0}}^{3}
\mbox{\hspace{312pt}}
\end{eqnarray*}
of weight 6.
\begin{eqnarray*}
&&M_3=\frac 1 {576}\,(P_4,P_2^2)^{(4)}={Q_{4, \,0}}\,{Q_{0, \,2}}^{2} + 
4\,{Q_{2, \,2}}\,{Q_{1, \,1}}^{2} - 4\,{Q_{3, \,1}}\,{Q_{0, \,2}}\,{Q_{1, \,1}}\\
&& + 2\,{Q_{2, \,2}}\,
{Q_{0, \,2}}\,{Q_{2, \,0}} - 4\,{Q_{1, \,3}}\,{Q_{2, \,0}}
\,{Q_{1, \,1}}  + {Q_{0, \,4}}\,{Q_{2, \,0}}^{2}
\end{eqnarray*}
of weight 4.
\begin{eqnarray*}
&& M_4=\frac 1 {1194393600}\,(T_4,P_2^3)^{(6)}\\
&&= 
 4\,{Q_{3, \,1}}^{2}\,{Q_{0, \,4}}\,{Q_{1, \,1}}^{3} - 4
\,{Q_{1, \,3}}^{2}\,{Q_{4, \,0}}\,{Q_{1, \,1}}^{3} - {Q_{4, \,0}}
^{2}\,{Q_{1, \,3}}\,{Q_{0, \,2}}^{3} + {Q_{3, \,1}}\,{Q_{0, \,4}}
^{2}\,{Q_{2, \,0}}^{3} \\
&& + 2\,{Q_{1, \,3}}^{3}\,{Q_{2, \,0}}^{3}
 - 2\,{Q_{3,\,1}}^{3}\,{Q_{0, \,2}}^{3} + 9\,{Q_{2, \,2}}^{2}\,{Q_{0, \,4}}\,
{Q_{2, \,0}}^{2}\,{Q_{1, \,1}}   + 3\,{Q_{4, \,0}}\,{Q_{2, \,2}}\,
{Q_{3, \,1}}\,{Q_{0, \,2}}^{3}\\
&& + {Q_{4, \,0}}^{2}\,{Q_{0, \,4}}\,{Q_{0, \,2}}^{2}\,{
Q_{1, \,1}} + 2\,{Q_{3, \,1}}\,{Q_{1, \,3}}^{2}\,{Q_{2, \,0}}^{2}
\,{Q_{0, \,2}}  - 9\,{Q_{4, \,0}}\,{Q_{2, \,2}}^{2}\,{Q_{0, \,2}}^{2
}\,{Q_{1, \,1}}\\
&& + 6\,{Q_{3, \,1}}^{2}\,{Q_{2, \,2}}\,{Q_{0, \,2}}
^{2}\,{Q_{1, \,1}} - 2\,{Q_{3, \,1}}^{2}\,{Q_{1, \,3}}\,{Q_{0, \,2}}^{2}\,{Q_{2, \,0}}
 - 8\,{Q_{3, \,1}}^{2}\,{Q_{1, \,3}}\,{Q_{1, \,1}}^{2}\,{Q_{0, \,2}}\\
&& - 6\,{Q_{2, \,2}}\,{Q_{1, \,3}}^{2}\,{Q_{2, \,0}}
^{2}\,{Q_{1, \,1}} - {Q_{4, \,0}}\,{Q_{0, \,4}}^{2}\,{Q_{2, \,0}}
^{2}\,{Q_{1, \,1}} - 3\,{Q_{2, \,2}}\,{Q_{1, \,3}}\,{Q_{0, \,4}}\,{Q_{2
, \,0}}^{3}  \\
&&  + 8\,{Q_{3, \,1}}\,{Q_{1, \,3}}^{2}\,{Q_{2, \,0}}\,{Q_{1, \,1}}^{2}- 12\,{Q_{3, \,1}}\,{Q_{0, \,4}}\,{Q_{2, \,2}}\,{Q_{2, \,0}}\,{Q_{1, \,1}}^{2} 
- 2\,{Q_{3, \,1}}\,{Q_{1, \,3}}\,{Q_{0, \,4}}\,{Q_{2, \,0}}^{2}\,{Q_{1, \,1}}
 \\
&& + 3\,{Q_{4, \,0}}\,{Q_{2, \,2}}\,{Q_{1, \,3}}\,{Q_{0, \,2}}^{2}\,{Q_{2, \,0}}
+ 12\,{Q_{4, \,0}}\,{Q_{2, \,2}}\,{Q_{1, \,3}}\,{Q_{1, \,1}}^{2}\,{Q_{0, \,2}} 
 - 3\,{Q_{3, \,1}}\,{Q_{0, \,4}}\,{Q_{2, \,2}}\,{Q_{2, \,0}}^{2}
\,{Q_{0, \,2}}   \\
&&+4\,{Q_{4, \,0}}\,{Q_{0, \,4}}\,{Q_{1, \,3}}\,{Q_{2, \,0}}
\,{Q_{1, \,1}}^{2} + 2\,{Q_{4, \,0}}\,{Q_{1, \,3}}\,{Q_{3, \,1}}\,{Q_{0
, \,2}}^{2}\,{Q_{1, \,1}} + {Q_{4, \,0}}\,{Q_{0, \,4}}\,{Q_{1, \,
3}}\,{Q_{2, \,0}}^{2}\,{Q_{0, \,2}} \\
&& - 4\,{Q_{4, \,0}}\,{Q_{0, \,4}}
\,{Q_{3, \,1}}\,{Q_{1, \,1}}^{2}\,{Q_{0, \,2}} - {Q_{4, \,0}}\,{Q_{0, \,4}}\,{Q_{3, \,1}}\,{Q_{0, \,2}}^{2}\,{Q_{2, \,0}}
  - 6\,{Q_{1, \,3}}^{2}\,{Q_{4, \,0}}\,{Q_{2, \,0}}\,{Q_{0, \,2}}
\,{Q_{1, \,1}}\\
&& + 6\,{Q_{3, \,1}}^{2}\,{Q_{0, \,4}}\,{Q_{2, \,0}}\,
{Q_{0, \,2}}\,{Q_{1, \,1}},
\end{eqnarray*}
of weight 9.
\begin{eqnarray*}
&&M_5=\frac 1 {238878720}(S,P_4)^{(4)}=\\
\lefteqn{ - 6\,{Q_{1, \,3}}^{2}\,{Q_{3, \,1}}\,{Q_{2, \,1}}^{2}
 + 6\,{Q_{3, \,1}}^{2}\,{Q_{1, \,2}}^{2}\,{Q_{1, \,3}} + {Q_{1, 
\,3}}\,{Q_{4, \,0}}^{2}\,{Q_{0, \,3}}^{2} - 6\,{Q_{4, \,0}}\,{Q_{
2, \,2}}\,{Q_{2, \,1}}\,{Q_{0, \,3}}\,{Q_{1, \,3}}} \\
 & & \mbox{} + 2\,{Q_{0, \,4}}\,{Q_{4, \,0}}\,{Q_{0, \,3}}\,{Q_{2
, \,1}}\,{Q_{3, \,1}} - 2\,{Q_{1, \,3}}\,{Q_{3, \,1}}\,{Q_{0, \,3
}}\,{Q_{1, \,2}}\,{Q_{4, \,0}} + 2\,{Q_{3, \,1}}^{3}\,{Q_{0, \,3}
}^{2} \\
 & & \mbox{} - {Q_{0, \,4}}\,{Q_{4, \,0}}^{2}\,{Q_{0, \,3}}\,{Q_{
1, \,2}} - {Q_{3, \,1}}^{2}\,{Q_{3, \,0}}\,{Q_{0, \,3}}\,{Q_{0, 
\,4}} + 4\,{Q_{3, \,1}}^{2}\,{Q_{2, \,1}}\,{Q_{0, \,3}}\,{Q_{1, 
\,3}} \\
 & & \mbox{} + {Q_{1, \,3}}^{2}\,{Q_{4, \,0}}\,{Q_{0, \,3}}\,{Q_{
3, \,0}} - 3\,{Q_{4, \,0}}\,{Q_{2, \,2}}\,{Q_{0, \,3}}^{2}\,{Q_{3
, \,1}} + 9\,{Q_{4, \,0}}\,{Q_{2, \,2}}^{2}\,{Q_{1, \,2}}\,{Q_{0
, \,3}} \\
 & & \mbox{} - 6\,{Q_{3, \,1}}^{2}\,{Q_{1, \,2}}\,{Q_{0, \,3}}\,{
Q_{2, \,2}} + 6\,{Q_{1, \,3}}^{2}\,{Q_{2, \,1}}\,{Q_{3, \,0}}\,{Q
_{2, \,2}} - 2\,{Q_{1, \,3}}^{3}\,{Q_{3, \,0}}^{2} \\
 & & \mbox{} + 3\,{Q_{0, \,4}}\,{Q_{4, \,0}}\,{Q_{1, \,2}}^{2}\,{
Q_{3, \,1}} + 9\,{Q_{1, \,3}}^{2}\,{Q_{4, \,0}}\,{Q_{2, \,1}}\,{Q
_{1, \,2}} - 9\,{Q_{3, \,1}}^{2}\,{Q_{2, \,1}}\,{Q_{1, \,2}}\,{Q
_{0, \,4}} \\
 & & \mbox{} + 9\,{Q_{2, \,2}}\,{Q_{3, \,1}}\,{Q_{2, \,1}}^{2}\,{
Q_{0, \,4}} + 3\,{Q_{1, \,3}}\,{Q_{2, \,2}}\,{Q_{3, \,0}}^{2}\,{Q
_{0, \,4}} - 4\,{Q_{1, \,3}}^{2}\,{Q_{3, \,1}}\,{Q_{3, \,0}}\,{Q
_{1, \,2}} \\
 & & \mbox{} - 9\,{Q_{4, \,0}}\,{Q_{2, \,2}}\,{Q_{1, \,2}}^{2}\,{
Q_{1, \,3}} - 3\,{Q_{1, \,3}}\,{Q_{4, \,0}}\,{Q_{2, \,1}}^{2}\,{Q
_{0, \,4}} - 9\,{Q_{2, \,2}}^{2}\,{Q_{3, \,0}}\,{Q_{2, \,1}}\,{Q
_{0, \,4}} \\
 & & \mbox{} + {Q_{0, \,4}}^{2}\,{Q_{4, \,0}}\,{Q_{3, \,0}}\,{Q_{
2, \,1}} - {Q_{0, \,4}}^{2}\,{Q_{3, \,1}}\,{Q_{3, \,0}}^{2} + 2\,
{Q_{1, \,3}}\,{Q_{3, \,1}}\,{Q_{3, \,0}}\,{Q_{2, \,1}}\,{Q_{0, \,
4}} \\
 & & \mbox{} + 6\,{Q_{2, \,2}}\,{Q_{3, \,1}}\,{Q_{3, \,0}}\,{Q_{1
, \,2}}\,{Q_{0, \,4}} - 2\,{Q_{1, \,3}}\,{Q_{4, \,0}}\,{Q_{3, \,0
}}\,{Q_{1, \,2}}\,{Q_{0, \,4}}\mbox{\hspace{108pt}}
\end{eqnarray*}
of weight 9.
\end{small}

We remind the reader that 
\begin{eqnarray*}
H_3=(P_4,P_4)^{(2)},\quad H_4=(P_4,P_4)^{(2)}, \\
T_3=(H_3,P_3)^{(1)},\quad T_4=(H_4,P_4)^{(1)},\\
 S=(H_4,P_3^2)^{(3)}.
\end{eqnarray*}

The following  {\sc Maple} code was used to  compute eight 
fundamental invariants: 
\begin{equation*}
 I_1=\frac{M_1}{{d_2}^2},\qquad I_2=\frac {M_2} {{d_2}^3},\qquad
I_3=\frac {d_3}{{d_2}^3}, 
\end{equation*}
\begin{equation*}
I_4=\frac j {{d_2}^3},\quad I_5=\frac{i}{{d_2}^2},\quad
I_6=\frac {{M_4}^2} {{d_2}^9},\quad I_7=\frac {M_3} {{d_2}^2},\quad  
I_8=\frac {{M_5}^2}{{d_2}^9}.
\end{equation*}
restricted to a given polynomial.

\chapter{Computations  on Ternary Cubics.}

\DefineParaStyle{Maple Output}
\DefineCharStyle{2D Math}
\DefineCharStyle{2D Output}

\begin{maplegroup}
\begin{center}
REDUCIBLE CUBICS IN THREE VARIABLES
\end{center}

\end{maplegroup}
\begin{maplegroup}
The following standard packages are used:

\end{maplegroup}
\begin{maplegroup}
\begin{mapleinput}
\mapleinline{active}{1d}{with(linalg):with(Groebner):}{%
}
\end{mapleinput}
\end{maplegroup}
\begin{maplegroup}
Our code includes the following programs:

\end{maplegroup}

\begin{maplegroup}
 \underline{Pinv} computes  fundamental  invariants.

\end{maplegroup}
\begin{maplegroup}
\begin{mapleinput}
\mapleinline{active}{1d}{read ternary3;}{%
}
\end{mapleinput}

\end{maplegroup}
\begin{maplegroup}
\underline{Psignature }computes syzygies between fundamental 
invariants.

\end{maplegroup}
\begin{maplegroup}
\begin{mapleinput}
\mapleinline{active}{1d}{read Psignature;}{%
}
\end{mapleinput}

\end{maplegroup}
\begin{maplegroup}
\begin{center}
\underline{Two-dimensional unimodular group of  symmetries}
\end{center}

\end{maplegroup}
\begin{maplegroup}
\begin{mapleinput}
\mapleinline{active}{1d}{F:=(x,y,z)->x*y*z;}{%
}
\end{mapleinput}

\mapleresult
\begin{maplelatex}
\[
F := (x, \,y, \,z)\rightarrow x\,y\,z
\]
\end{maplelatex}

\end{maplegroup}
\begin{maplegroup}
\begin{mapleinput}
\mapleinline{active}{1d}{f:=(p,q)->p*q;}{%
}
\end{mapleinput}

\mapleresult
\begin{maplelatex}
\[
f := (p, \,q)\rightarrow p\,q
\]
\end{maplelatex}

\end{maplegroup}
\begin{maplegroup}
\begin{mapleinput}
\mapleinline{active}{1d}{Pinv(f,3);}{%
}
\end{mapleinput}

\mapleresult
\begin{maplelatex}
\[
[{\displaystyle \frac {4}{3}} , \,{\displaystyle \frac {16}{3}} 
, \,{\displaystyle \frac {16}{9}} ]
\]
\end{maplelatex}

\end{maplegroup}
\begin{maplegroup}
an equivalent polynomial:

\begin{mapleinput}
\mapleinline{active}{1d}{f:=(p,q)->1/2*p*(q^2-1);}{%
}
\end{mapleinput}

\mapleresult
\begin{maplelatex}
\[
f := (p, \,q)\rightarrow {\displaystyle \frac {1}{2}} \,p\,(q^{2}
 - 1)
\]
\end{maplelatex}

\end{maplegroup}
\begin{maplegroup}
\begin{mapleinput}
\mapleinline{active}{1d}{Pinv(f,3);}{%
}
\end{mapleinput}

\begin{mapleinput}
\end{mapleinput}

\mapleresult
\begin{maplelatex}
\[
[{\displaystyle \frac {4}{3}} , \,{\displaystyle \frac {16}{3}} 
, \,{\displaystyle \frac {16}{9}} ]
\]
\end{maplelatex}

\end{maplegroup}
\begin{maplegroup}
\begin{center}
\underline{Two-dimensional symmetry group with a  one-dimensional unimodular subgroup.
}
\end{center}

\end{maplegroup}
\begin{maplegroup}
\begin{mapleinput}
\mapleinline{active}{1d}{F:=(x,y,z)->(x^2+y*z)*z:}{%
}
\end{mapleinput}

\end{maplegroup}
\begin{maplegroup}
\begin{mapleinput}
\mapleinline{active}{1d}{f:=(p,q)->p^2+q;}{%
}
\end{mapleinput}

\mapleresult
\begin{maplelatex}
\[
f := (p, \,q)\rightarrow p^{2} + q
\]
\end{maplelatex}

\end{maplegroup}
\begin{maplegroup}
\begin{mapleinput}
\mapleinline{active}{1d}{Pinv(f,3);}{%
}
\end{mapleinput}

\mapleresult
\begin{maplelatex}
\[
[{\displaystyle \frac {-1}{6}} , \,{\displaystyle \frac {41}{6}} 
, \,{\displaystyle \frac {-2}{9}} ]
\]
\end{maplelatex}

\end{maplegroup}
\begin{maplegroup}
\mapleresult
\begin{maplettyout}
an equivalent polynomial:
\end{maplettyout}

\end{maplegroup}
\begin{maplegroup}
\begin{mapleinput}
\mapleinline{active}{1d}{F:=(x,y,z)->(x^2-y*z-z^2)*z:                 
   f:=(p,q)->p^2-q-1;}{%
}
\end{mapleinput}

\mapleresult
\begin{maplelatex}
\[
f := (p, \,q)\rightarrow p^{2} - q - 1
\]
\end{maplelatex}

\end{maplegroup}
\begin{maplegroup}
\begin{mapleinput}
\mapleinline{active}{1d}{Pinv(f,3);}{%
}
\end{mapleinput}

\mapleresult
\begin{maplelatex}
\[
[{\displaystyle \frac {-1}{6}} , \,{\displaystyle \frac {41}{6}} 
, \,{\displaystyle \frac {-2}{9}} ]
\]
\end{maplelatex}

\end{maplegroup}
\begin{maplegroup}
\begin{center}
{\underline{One-dimensional unimodular group of symmetries}}
\end{center}
\end{maplegroup}

\begin{maplegroup}
\begin{mapleinput}
\mapleinline{active}{1d}{F:=(x,y,z)->(x^2+y^2+z^2)*z;}{%
}
\end{mapleinput}

\mapleresult
\begin{maplelatex}
\[
F := (x, \,y, \,z)\rightarrow (x^{2} + y^{2} + z^{2})\,z
\]
\end{maplelatex}

\end{maplegroup}
\begin{maplegroup}
\begin{mapleinput}
\mapleinline{active}{1d}{f:=(p,q)->(p^2+q^2+1);}{%
}
\end{mapleinput}

\mapleresult
\begin{maplelatex}
\[
f := (p, \,q)\rightarrow p^{2} + q^{2} + 1
\]
\end{maplelatex}

\end{maplegroup}
\begin{maplegroup}
\begin{mapleinput}
\mapleinline{active}{1d}{Pinv(f,3);}{%
}
\end{mapleinput}

\mapleresult
\begin{maplelatex}
\begin{eqnarray*}
&&[{\displaystyle \frac {4}{3}} \,{\displaystyle \frac {(p^{2} + 3
 + q^{2})\,(p^{2} + q^{2})}{(p^{2} - 3 + q^{2})^{2}}} , \,
{\displaystyle \frac {16}{3}} \,{\displaystyle \frac {(p^{2} + q
^{2})\,(p^{4} + 2\,p^{2}\,q^{2} + 81 + q^{4})}{(p^{2} - 3 + q^{2}
)^{3}}} , \\
& &{\displaystyle \frac {16}{9}} \,{\displaystyle 
\frac {(9 + p^{2} + q^{2})\,(p^{2} + q^{2})^{2}}{(p^{2} - 3 + q^{
2})^{3}}} ]
\end{eqnarray*}
\end{maplelatex}

\end{maplegroup}
\begin{maplegroup}
\begin{mapleinput}
\mapleinline{active}{1d}{Psignature(inv);}{%
}
\end{mapleinput}

\begin{mapleinput}
\end{mapleinput}

\mapleresult
\begin{maplelatex}
\begin{eqnarray*}
\lefteqn{\mathit{elimination\ of\ u\ from\ the\ equations:}}\\
&&4\,
(p^{2} + 3 + q^{2})\,(p^{2} + q^{2}) - 3\,(p^{2} - 3 + q^{2})^{2}
\,{I_{1}},  \\
 & & 16\,(p^{2} + q^{2})\,(p^{4} + 2\,p^{2}\,q^{2} + 81 + q^{4})
 - 3\,(p^{2} - 3 + q^{2})^{3}\,{I_{2}},\\
 & & 16\,(9 + p^{2} + q^{2})
\,(p^{2} + q^{2})^{2} - 9\,(p^{2} - 3 + q^{2})^{3}\,{I_{3}},  \\
 & & 1 - 243\,(p^{2} - 3 + q^{2})^{3}\,w
\end{eqnarray*}
\end{maplelatex}

\begin{maplelatex}
\[
\mathit{dimension\ of\ the\ signature\ manifold}=1
\]
\end{maplelatex}

\begin{maplelatex}
\[
\mathit{the\ signature\ manifold\ is\ defined\ by:}
\]
\end{maplelatex}

\begin{maplelatex}
\begin{eqnarray*}
\lefteqn{[ - 1482\,{I_{3}}\,{I_{2}} + 8865\,{I_{3}}^{2} + 40\,{I
_{2}}^{2} - 1296\,{I_{1}} - 36\,{I_{2}} + 17280\,{I_{3}} - 18522
\,{I_{3}}\,{I_{1}}, } \\
 & & 40\,{I_{1}}\,{I_{2}} + 582\,{I_{3}}\,{I_{1}} + 42\,{I_{3}}\,
{I_{2}} - 315\,{I_{3}}^{2} - 144\,{I_{1}} - 4\,{I_{2}} - 480\,{I
_{3}},  \\
 & & 360\,{I_{1}}^{2} - 378\,{I_{3}}\,{I_{1}} - 18\,{I_{3}}\,{I_{
2}} + 135\,{I_{3}}^{2} - 144\,{I_{1}} - 4\,{I_{2}} + 120\,{I_{3}}
]\mbox{\hspace{21pt}}
\end{eqnarray*}
\end{maplelatex}

\end{maplegroup}
\begin{maplegroup}
an equivalent polynomial:

\end{maplegroup}
\begin{maplegroup}
\begin{mapleinput}
\mapleinline{active}{1d}{F:=(x,y,z)->x*y*z+x^2*y+z*y^2;}{%
}
\end{mapleinput}

\begin{mapleinput}
\mapleinline{active}{1d}{f:=(p,q)->p*q+p^2*q+q^2;}{%
}
\end{mapleinput}

\mapleresult
\begin{maplelatex}
\[
F := (x, \,y, \,z)\rightarrow x\,y\,z + x^{2}\,y + z\,y^{2}
\]
\end{maplelatex}

\begin{maplelatex}
\[
f := (p, \,q)\rightarrow p\,q + p^{2}\,q + q^{2}
\]
\end{maplelatex}

\end{maplegroup}
\begin{maplegroup}
\begin{mapleinput}
\mapleinline{active}{1d}{Pinv(f,3);}{%
}
\end{mapleinput}

\mapleresult
\begin{maplelatex}
\begin{small}

\begin{eqnarray*}
&&
[{\displaystyle \frac {4}{3}} \,{\displaystyle \frac {(p
 + 1 - q)\,(p + q)\,(p^{2} + p + 2\,q^{2} + q)}{( - 4\,q^{2} + q
 + p + p^{2})^{2}}} , \\
&  & \frac {16}{3} \,
 \frac {(p + 1 - q)\,(p + q)\,(p^{4} + 2\,p^{3} + 2
\,p^{2}\,q - 2\,p^{2}\,q^{2} + p^{2} + 2\,p\,q - 2\,p\,q^{2} - 2
\,q^{3} + 82\,q^{4} + q^{2})}{( - 4\,q^{2} + q + p + p^{2})^{3}}, \\
 & & \,{\displaystyle \frac {16}{9}} \,{\displaystyle \frac {(p
 + p^{2} + q + 8\,q^{2})\,(p + q)^{2}\,(p + 1 - q)^{2}}{( - 4\,q
^{2} + q + p + p^{2})^{3}}} ]
\end{eqnarray*}
\end{small}
\end{maplelatex}

\end{maplegroup}
\begin{maplegroup}
\begin{mapleinput}
\mapleinline{active}{1d}{Psignature(inv);}{%
}
\end{mapleinput}

\mapleresult
\begin{maplelatex}
\begin{eqnarray*}
\lefteqn{\mathit{elimination\ of\ u\ from\ the\ equations:}}\\
& &4\,(p + 1 - q)\,(p + q)\,(p^{2} + p + 2\,q^{2} + q) - 3\,\mathrm{\%1
}^{2}\,{I_{1}}, \\
 & & 16\,(p + 1 - q)\,(p + q)\,(p^{4} + 2\,p^{3} + 2\,p^{2}\,q - 
2\,p^{2}\,q^{2} + p^{2} + 2\,p\,q - 2\,p\,q^{2} - 2\,q^{3}\\
  & &+ 82\,q^{4} + q^{2}) - 3\,\mathrm{\%1}^{3}\,{I_{2}},  \\
 & & 16\,(p + p^{2} + q + 8\,q^{2})\,(p + q)^{2}\,(p + 
1 - q)^{2} - 9\,\mathrm{\%1}^{3}\,{I_{3}}, \,1 - 243\,\mathrm{\%1
}^{3}\,w \\
 & & \mathrm{\%1} :=  - 4\,q^{2} + q + p + p^{2}
\end{eqnarray*}
\end{maplelatex}

\begin{maplelatex}
\[
\mathit{dimension\ of\ the\ signature\ manifold}=1
\]
\end{maplelatex}

\begin{maplelatex}
\[
\mathit{the\ signature\ manifold\ is\ defined\ by:}
\]
\end{maplelatex}

\begin{maplelatex}
\begin{eqnarray*}
\lefteqn{[ - 1482\,{I_{3}}\,{I_{2}} + 8865\,{I_{3}}^{2} + 40\,{I
_{2}}^{2} - 1296\,{I_{1}} - 36\,{I_{2}} + 17280\,{I_{3}} - 18522
\,{I_{3}}\,{I_{1}}, } \\
 & & 40\,{I_{1}}\,{I_{2}} + 582\,{I_{3}}\,{I_{1}} + 42\,{I_{3}}\,
{I_{2}} - 315\,{I_{3}}^{2} - 144\,{I_{1}} - 4\,{I_{2}} - 480\,{I
_{3}},  \\
 & & 360\,{I_{1}}^{2} - 378\,{I_{3}}\,{I_{1}} - 18\,{I_{3}}\,{I_{
2}} + 135\,{I_{3}}^{2} - 144\,{I_{1}} - 4\,{I_{2}} + 120\,{I_{3}}
]\mbox{\hspace{21pt}}
\end{eqnarray*}
\end{maplelatex}

\end{maplegroup}
\begin{maplegroup}
\begin{mapleinput}
\end{mapleinput}

\end{maplegroup}


\newpage

\DefineParaStyle{Error}
\DefineParaStyle{Maple Output}
\DefineParaStyle{Warning}
\DefineCharStyle{2D Math}
\DefineCharStyle{2D Output}

\begin{maplegroup}
\begin{center}
IRREDUCIBLE CUBICS IN THREE VARIABLES
\end{center}
\end{maplegroup}

\begin{maplegroup}
The program \underline{SymN} computes the cardinality of the symmetry group in the case
 it is finite. 
\end{maplegroup}
\begin{maplegroup}
\begin{mapleinput}
\mapleinline{active}{1d}{read symN;}{%
}
\end{mapleinput}
\end{maplegroup}
\begin{maplegroup}
it uses \underline{kbasis5} adopted from \cite{CLO97}
 to compute the number of the solutions
for a system of polynomial equations:
\end{maplegroup}
\begin{maplegroup}
\begin{mapleinput}
\mapleinline{active}{1d}{read kbasis5;}{%
}
\end{mapleinput}

\end{maplegroup}
\begin{maplegroup}
\begin{center}
\textbf{Singular Curves:}
\end{center}

\end{maplegroup}
\begin{maplegroup}
{\bf q\symbol{94}2=p\symbol{94}3:}  one-dimensional symmetry group, finite 
number of unimodular symmetries:
\end{maplegroup}
\begin{maplegroup}
\begin{mapleinput}
\mapleinline{active}{1d}{f:=(p,q)->p^3-q^2;}{%
}
\end{mapleinput}

\mapleresult
\begin{maplelatex}
\[
f := (p, \,q)\rightarrow p^{3} - q^{2}
\]
\end{maplelatex}

\end{maplegroup}
\begin{maplegroup}
\begin{mapleinput}
\end{mapleinput}

\end{maplegroup}
\begin{maplegroup}
\begin{mapleinput}
\mapleinline{active}{1d}{Pinv(f,3);}{%
}
\end{mapleinput}

\mapleresult
\begin{maplelatex}
\[
[{\displaystyle \frac {-1}{6}} , \,{\displaystyle \frac {1}{6}} 
\,{\displaystyle \frac {36\,p^{3} + 5\,q^{2}}{p^{3}}} , \, - 
{\displaystyle \frac {1}{9}} \,{\displaystyle \frac {p^{3} + q^{2
}}{p^{3}}} ]
\]
\end{maplelatex}

\end{maplegroup}
\begin{maplegroup}
\begin{mapleinput}
\mapleinline{active}{1d}{Psignature(inv);}{%
}
\end{mapleinput}

\mapleresult
\begin{maplelatex}
\begin{eqnarray*}
&&\mathit{elimination\ of\ u\ from\ the\ equations:}\\
&&  - 1 - 6\,{I
_{1}}, \,36\,p^{3} + 5\,q^{2} - 6\,{I_{2}}\,p^{3}, \, - p^{3} - q
^{2} - 9\,{I_{3}}\,p^{3}, \,1 - 18\,w\,p^{3}
\end{eqnarray*}
\end{maplelatex}

\begin{maplelatex}
\[
\mathit{dimension\ of\ the\ signature\ manifold}=1
\]
\end{maplelatex}

\begin{maplelatex}
\[
\mathit{the\ signature\ manifold\ is\ defined\ by:}
\]
\end{maplelatex}

\begin{maplelatex}
\[
[6\,{I_{2}} + 45\,{I_{3}} - 31, \,1 + 6\,{I_{1}}]
\]
\end{maplelatex}

\end{maplegroup}
\begin{maplegroup}
\textbf{q\symbol{94}2=p\symbol{94}2*(p+1): } The
projective symmetry group has \textbf{6} elements.

\end{maplegroup}
\begin{maplegroup}
\begin{mapleinput}
\mapleinline{active}{1d}{f:=(p,q)->p^2*(p+1)-q^2;}{%
}
\end{mapleinput}

\mapleresult
\begin{maplelatex}
\[
f := (p, \,q)\rightarrow p^{2}\,(p + 1) - q^{2}
\]
\end{maplelatex}

\end{maplegroup}
\begin{maplegroup}
\begin{mapleinput}
\mapleinline{active}{1d}{Pinv(f,3);}{%
}
\end{mapleinput}

\mapleresult
\begin{maplelatex}
\begin{eqnarray*}
&  &[ - {\displaystyle \frac {1}{6}} \,{\displaystyle 
\frac {(3\,p + 4)\,(q + p)\,(q - p)\,(3\,p\,q^{2} - 2\,q^{2} + 2
\,p^{2} + 3\,p^{3})}{(3\,p\,q^{2} + q^{2} - p^{2})^{2}}} ,\\
& &{\displaystyle \frac {1}{6}} ((q + p)\,(q - p)(135\,q^{6} + 
32\,q^{4} + 972\,p^{3}\,q^{4} + 1107\,p^{2}\,q^{4} + 72\,p\,q
^{4} + 1269\,p^{4}\,q^{2}\\
& & - 144\,p^{3}\,q^{2} - 64\,p^{2}\,q^{2}
 + 972\,q^{2}\,p^{5} + 72\,p^{5} + 81\,p^{6} + 32\,p^{4}))/(3\,p\,q^{2}
 + q^{2} - p^{2})^{3},\\
&   &  - {\displaystyle \frac {1}{9}} ( \\
 & & ( - 16\,q^{2} + 27\,q^{4} + 72\,p\,q^{2} + 16\,p^{2} + 81\,p
^{2}\,q^{2} + 72\,p^{3} + 27\,p^{3}\,q^{2} + 108\,p^{4} + 54\,p^{
5}) \\
 & & (q - p)^{2}\,(q + p)^{2}) \left/ {\vrule 
height0.44em width0em depth0.44em} \right. \!  \! (3\,p\,q^{2} + 
q^{2} - p^{2})^{3}]\mbox{\hspace{207pt}}
\end{eqnarray*}
\end{maplelatex}

\end{maplegroup}
\begin{maplegroup}
\begin{mapleinput}
\mapleinline{active}{1d}{symN(inv,1,2);}{%
}
\end{mapleinput}

\mapleresult
\begin{maplelatex}
\[
\mathit{the\ number\ of\ symmetries}=6
\]
\end{maplelatex}

\end{maplegroup}
\begin{maplegroup}
In {\tt symN} we have  chosen the point  $P=1$ and $Q=2$ to substitute into
equations (\ref{symmN}). We find the number of the solutions using procedure 
{\tt kbasis5},
but in this particular case it is not difficult to solve the equations explicitly:

\end{maplegroup}
\begin{maplegroup}
\begin{mapleinput}
\mapleinline{active}{1d}{E;}{%
}
\end{mapleinput}

\mapleresult
\begin{maplelatex}
\begin{eqnarray*}
\lefteqn{[76\,p^{2} - 49\,q^{2} + 72\,p + 48, \,7543\,p\,q^{2} - 
6633\,q^{2} - 4008\,p + 368, } \\
 & & 20961997\,q^{4} - 83898780\,q^{2} - 13605280\,p + 13808448]
\mbox{\hspace{0pt}}
\end{eqnarray*}
\end{maplelatex}
\end{maplegroup}

\begin{maplegroup}
\begin{mapleinput}
\mapleinline{active}{1d}{map(allvalues,[solve(\{op(E)\},\{p,q\})]);}{%
}
\end{mapleinput}

\mapleresult
\begin{maplelatex}
\begin{eqnarray*}
\lefteqn{[\{p=1, \,q=2\}, \,\{q=-2, \,p=1\},} \\
& & \{p= - 
{\displaystyle \frac {212}{397}}  + {\displaystyle \frac {112}{
397}} \,I\,\sqrt{3}, \,q= - {\displaystyle \frac {208}{397}}  + 
{\displaystyle \frac {20}{397}} \,I\,\sqrt{3}\},  \\
 & & \{p= - {\displaystyle \frac {212}{397}}  - {\displaystyle 
\frac {112}{397}} \,I\,\sqrt{3}, \,q= - {\displaystyle \frac {208
}{397}}  - {\displaystyle \frac {20}{397}} \,I\,\sqrt{3}\},  \\
 & & \{p= - {\displaystyle \frac {212}{397}}  - {\displaystyle 
\frac {112}{397}} \,I\,\sqrt{3}, \,q={\displaystyle \frac {208}{
397}}  + {\displaystyle \frac {20}{397}} \,I\,\sqrt{3}\},  \\
 & & \{p= - {\displaystyle \frac {212}{397}}  + {\displaystyle 
\frac {112}{397}} \,I\,\sqrt{3}, \,q={\displaystyle \frac {208}{
397}}  - {\displaystyle \frac {20}{397}} \,I\,\sqrt{3}\}]
\mbox{\hspace{122pt}}
\end{eqnarray*}
\end{maplelatex}

\end{maplegroup}
\begin{maplegroup}
The list above contains the images of the point (1,2)  under all possible symmetries.
We note, however, that if we put a non-generic point into {\tt symN} we might obtain an incorrect answer or no answer at all:
\end{maplegroup}

\begin{maplegroup}
\begin{mapleinput}
\mapleinline{active}{1d}{symN(inv,1,1);}{%
}
\end{mapleinput}

\mapleresult
\begin{maplettyout}
Error, (in kbasis) Ideal is not zero-dimensional, no finite basis
\end{maplettyout}

\end{maplegroup}
\begin{maplegroup}
Other generic points produce the correct result for the cardinality of the symmetry group.
\end{maplegroup}
\begin{maplegroup}
\begin{mapleinput}
\mapleinline{active}{1d}{symN(inv,0,1);}{%
}
\end{mapleinput}

\mapleresult
\begin{maplelatex}
\[
\mathit{the\ number\ of\ symmetries}=6
\]
\end{maplelatex}

\end{maplegroup}
\begin{maplegroup}
\begin{mapleinput}
\mapleinline{active}{1d}{symN(inv,1,0);}{%
}
\end{mapleinput}

\mapleresult
\begin{maplelatex}
\[
\mathit{the\ number\ of\ symmetries}=6
\]
\end{maplelatex}

\end{maplegroup}
\begin{maplegroup}
\begin{mapleinput}
\end{mapleinput}

\end{maplegroup}

\begin{maplegroup}
An equivalent polynomial:

\end{maplegroup}
\begin{maplegroup}
\begin{mapleinput}
\mapleinline{active}{1d}{f:=(p,q)->p^2*(p+4)-q^2;}{%
}
\end{mapleinput}

\mapleresult
\begin{maplelatex}
\[
f := (p, \,q)\rightarrow p^{2}\,(p + 4) - q^{2}
\]
\end{maplelatex}

\end{maplegroup}
\begin{maplegroup}
\begin{mapleinput}
\mapleinline{active}{1d}{Pinv(f,3):}{%
}
\end{mapleinput}

\end{maplegroup}

\begin{maplegroup}
\begin{mapleinput}
\mapleinline{active}{1d}{symN(inv,1,0);}{%
}
\end{mapleinput}
\end{maplegroup}

\begin{maplegroup}
\begin{maplelatex}
\[
\mathit{the\ number\ of\ symmetries}=6
\]
\end{maplelatex}
\end{maplegroup}

\begin{maplegroup}
\begin{center}
\textbf{Non singular (elliptic) curves:}
\end{center}

\end{maplegroup}
\begin{maplegroup}
\textbf{q\symbol{94}2=p\symbol{94}3+1}. The number of  the 
projective symmetries is \textbf{54.}

\end{maplegroup}
\begin{maplegroup}
\begin{mapleinput}
\mapleinline{active}{1d}{f:=(p,q)->p^3-q^2+1;}{%
}
\end{mapleinput}

\mapleresult
\begin{maplelatex}
\[
f := (p, \,q)\rightarrow p^{3} - q^{2} + 1
\]
\end{maplelatex}

\end{maplegroup}
\begin{maplegroup}
\begin{mapleinput}
\mapleinline{active}{1d}{Pinv(f,3);}{%
}
\end{mapleinput}

\mapleresult
\begin{maplelatex}
\begin{eqnarray*}
&&[{\displaystyle \frac {-1}{6}} , {\displaystyle \frac {1}{6}} 
( - 135 + 432\,p^{3} + 540\,p^{6}\,q^{2} + 90\,q^{4} - 540
\,p^{6} + 40\,q^{6} + 5\,q^{8} + 2052\,p^{3}\,q^{2}\\
&& - 216\,p^{3}
\,q^{4} + 36\,p^{3}\,q^{6})/(p^{3}\,(q^{2} + 3)^{3}), \\
&& - 
{\displaystyle \frac {1}{9}} \,{\displaystyle \frac {(q^{2} - 1
 + p^{3})\,(q^{6} - 18\,q^{4} + 81\,q^{2} + 27\,p^{6})}{p^{3}\,(q
^{2} + 3)^{3}}} ]
\end{eqnarray*}
\end{maplelatex}
\end{maplegroup}

\begin{maplegroup}
\begin{mapleinput}
\mapleinline{active}{1d}{symN(inv,1,0);}{%
}
\end{mapleinput}

\mapleresult
\begin{maplelatex}
\[
\mathit{the\ number\ of\ symmetries}=54
\]
\end{maplelatex}

\end{maplegroup}
\begin{maplegroup}
\begin{mapleinput}
\mapleinline{active}{1d}{symN(inv,1,2);}{%
}
\end{mapleinput}

\mapleresult
\begin{maplelatex}
\[
\mathit{the\ number\ of\ symmetries}=54
\]
\end{maplelatex}

\end{maplegroup}

\begin{maplegroup}
\begin{mapleinput}
\mapleinline{active}{1d}{symN(inv,2,2);}{%
}
\end{mapleinput}

\mapleresult
\begin{maplelatex}
\[
\mathit{the\ number\ of\ symmetries}=54
\]
\end{maplelatex}

\end{maplegroup}

\begin{maplegroup}
An equivalent polynomial:

\end{maplegroup}
\begin{maplegroup}
\begin{mapleinput}
\mapleinline{active}{1d}{f:=(p,q)->p^3+q^3+1;}{%
}
\end{mapleinput}

\mapleresult
\begin{maplelatex}
\[
f := (p, \,q)\rightarrow p^{3} + q^{3} + 1
\]
\end{maplelatex}

\end{maplegroup}
\begin{maplegroup}
\begin{mapleinput}
\mapleinline{active}{1d}{Pinv(f,3);}{%
}
\end{mapleinput}

\mapleresult
\begin{maplelatex}
\begin{eqnarray*}
\lefteqn{[{\displaystyle \frac {-1}{6}} , \, - {\displaystyle 
\frac {1}{6}} \,{\displaystyle \frac {5\,q^{6} + 5\,p^{3}\,q^{6}
 + 5\,q^{3} - 26\,p^{3}\,q^{3} + 5\,p^{6}\,q^{3} + 5\,p^{3} + 5\,
p^{6}}{p^{3}\,q^{3}}} , } \\
 & &  - {\displaystyle \frac {1}{36}} \,{\displaystyle \frac {(q
^{3} + 1 - p^{3})\,(q^{3} - 1 - p^{3})\,(q^{3} - 1 + p^{3})}{p^{3
}\,q^{3}}} ]\mbox{\hspace{47pt}}
\end{eqnarray*}
\end{maplelatex}

\end{maplegroup}
\begin{maplegroup}
\begin{mapleinput}
\mapleinline{active}{1d}{symN(inv,1,2);}{%
}
\end{mapleinput}

\mapleresult
\begin{maplelatex}
\[
\mathit{the\ number\ of\ symmetries}=54
\]
\end{maplelatex}

\end{maplegroup}
\begin{maplegroup}
\begin{mapleinput}
\mapleinline{active}{1d}{symN(inv,1,0);}{%
}
\end{mapleinput}

\mapleresult
\begin{maplettyout}
Error, (in symN) division by zero
\end{maplettyout}

\end{maplegroup}
\begin{maplegroup}
\begin{mapleinput}
\mapleinline{active}{1d}{symN(inv,1,1);}{%
}
\end{mapleinput}

\mapleresult
\begin{maplelatex}
\[
\mathit{the\ number\ of\ symmetries}=54
\]
\end{maplelatex}

\end{maplegroup}
\begin{maplegroup}
The signature manifold for this class of polynomials is defined by:

\end{maplegroup}
\begin{maplegroup}
\begin{mapleinput}
\mapleinline{active}{1d}{Psignature(inv);}{%
}
\end{mapleinput}

\mapleresult
\begin{maplelatex}
\begin{eqnarray*}
\lefteqn{\mathit{elimination\ of\ u\ from\ the\ equations:}}\\
&& - 1 - 6\,{I_{1}},  \\
 & &  - 5\,p^{6}\,q^{3} - 5\,p^{6} - 5\,p^{3}\,q^{6} - 5\,p^{3}
 + 26\,p^{3}\,q^{3} - 5\,q^{6} - 5\,q^{3} - 6\,p^{3}\,q^{3}\,{I_{
2}},  \\
 & & - (p^{3} - 1 - q^{3})\,(p^{3} - 1 + q^{3})\,(p^{3
} + 1 - q^{3}) - 36\,p^{3}\,q^{3}\,{I_{3}},\\
& & \,1 - 36\,w\,p^{3}\,q^{3}
\end{eqnarray*}
\end{maplelatex}

\begin{maplelatex}
\[
\mathit{dimension\ of\ the\ signature\ manifold}=2
\]
\end{maplelatex}

\begin{maplelatex}
\[
\mathit{the\ signature\ manifold\ is\ defined\ by:}
\]
\end{maplelatex}

\begin{maplelatex}
\[
[1 + 6\,{I_{1}}]
\]
\end{maplelatex}

\end{maplegroup}
\begin{maplegroup}
\textbf{q\symbol{94}2=p\symbol{94}3+ap. } They are equivalent for all
\textbf{a}. The number of symmetries is \textbf{36}.

\end{maplegroup}
\begin{maplegroup}
\begin{mapleinput}
\mapleinline{active}{1d}{f:=(p,q)->p^3+p-q^2;}{%
}
\end{mapleinput}

\mapleresult
\begin{maplelatex}
\[
f := (p, \,q)\rightarrow p^{3} + p - q^{2}
\]
\end{maplelatex}

\end{maplegroup}
\begin{maplegroup}
\begin{mapleinput}
\mapleinline{active}{1d}{Pinv(f,3);}{%
}
\end{mapleinput}

\mapleresult
\begin{maplelatex}
\begin{eqnarray*}
\lefteqn{[ - {\displaystyle \frac {1}{6}} \,{\displaystyle 
\frac {27\,q^{4} + 9\,p^{2}\,q^{4} + 1 + 21\,p^{2} + 63\,p^{4} + 
27\,p^{6} - 60\,p\,q^{2} - 36\,p^{3}\,q^{2}}{(3\,p\,q^{2} - 1 + 3
\,p^{2})^{2}}} ,  - {\displaystyle \frac {1}{6}} (729\,p^{8}} \\
 & & \mbox{} - 2916\,p^{7}\,q^{2} - 972\,p^{6} - 2916\,p^{5}\,q^{
2} + 486\,p^{4}\,q^{4} + 1350\,p^{4} + 7236\,p^{3}\,q^{2} - 972\,
p^{3}\,q^{6} \\
 & & \mbox{} - 2052\,p^{2}\,q^{4} + 468\,p^{2} - 828\,p\,q^{2} - 
756\,q^{6}\,p + 41 - 135\,q^{8} - 378\,q^{4})/ \\
 & & (3\,p\,q^{2} - 1 + 3\,p^{2})^{3}, {\displaystyle \frac {1}{9
}} (2 + 36\,p^{2} + 216\,p^{4} + 540\,p^{6} + 486\,p^{8} - 27\,q
^{8} - 189\,q^{4} \\
 & & \mbox{} - 918\,p^{2}\,q^{4} - 81\,p^{4}\,q^{4} + 189\,q^{6}
\,p - 27\,p^{3}\,q^{6} + 117\,p\,q^{2} + 999\,p^{3}\,q^{2} - 81\,
p^{5}\,q^{2} \\
 & & \mbox{} - 243\,p^{7}\,q^{2})/(3\,p\,q^{2} - 1 + 3\,p^{2})^{3
}]
\end{eqnarray*}
\end{maplelatex}

\end{maplegroup}
\begin{maplegroup}
\begin{mapleinput}
\mapleinline{active}{1d}{symN(inv,1,2);}{%
}
\end{mapleinput}

\mapleresult
\begin{maplelatex}
\[
\mathit{the\ number\ of\ symmetries}=36
\]
\end{maplelatex}

\end{maplegroup}
\begin{maplegroup}
\begin{mapleinput}
\mapleinline{active}{1d}{symN(inv,2,3);}{%
}
\end{mapleinput}

\mapleresult
\begin{maplelatex}
\[
\mathit{the\ number\ of\ symmetries}=36
\]
\end{maplelatex}

\end{maplegroup}
\begin{maplegroup}
\begin{mapleinput}
\mapleinline{active}{1d}{f:=(p,q)->-q^2+2*p+p^3;}{%
}
\end{mapleinput}

\mapleresult
\begin{maplelatex}
\[
f := (p, \,q)\rightarrow  - q^{2} + 2\,p + p^{3}
\]
\end{maplelatex}

\end{maplegroup}
\begin{maplegroup}
\begin{mapleinput}
\mapleinline{active}{1d}{Pinv(f,3):}{%
}
\end{mapleinput}

\end{maplegroup}
\begin{maplegroup}
\begin{mapleinput}
\mapleinline{active}{1d}{symN(inv,1,2);}{%
}
\end{mapleinput}

\mapleresult
\begin{maplelatex}
\[
\mathit{the\ number\ of\ symmetries}=36
\]
\end{maplelatex}

\end{maplegroup}
\begin{maplegroup}
\textbf{q\symbol{94}2=p\symbol{94}3+ap+1.} This is a family of equivalence classes. The number    of symmetries is the same:  \textbf{18.}

\end{maplegroup}
\begin{maplegroup}
\begin{mapleinput}
\mapleinline{active}{1d}{f:=(x,y,z)->-q^2+4+2*p+p^3;}{%
}
\end{mapleinput}

\mapleresult
\begin{maplelatex}
\[
f := (x, \,y, \,z)\rightarrow  - q^{2} + 4 + 2\,p + p^{3}
\]
\end{maplelatex}

\end{maplegroup}
\begin{maplegroup}
\begin{mapleinput}
\mapleinline{active}{1d}{Pinv(f,3);}{%
}
\end{mapleinput}

\mapleresult
\begin{maplelatex}
\begin{small}
\begin{eqnarray*}
& &[ - {\displaystyle \frac {1}{6}} ( - 432\,q^{2} - 240\,p
\,q^{2} + 216\,p^{2}\,q^{2} - 72\,p^{3}\,q^{2} + 880 + 576\,p + 
1464\,p^{2} + 864\,p^{3}+\\ 
& &252\,p^{4}  \mbox{} + 54\,p^{6} + 54\,q^{4} + 9\,p^{2}\,q^{4})/(3\,p\,q
^{2} - 4 + 36\,p + 6\,p^{2})^{2},  - {\displaystyle \frac {1}{6}
} (225504\,p^{6} - 972\,p^{3}\,q^{6} \\
& &  + 774144\,p + 46656\,p^{7} + 2916\,p^{8} - 135\,q^{8} - 670464\,
p^{3} + 142272\,p\,q^{2} \\
 & & \mbox{} + 295488\,p^{2}\,q^{2} - 10368\,q^{2} + 512640\,p^{2
} - 58320\,p^{6}\,q^{2} - 5832\,p^{7}\,q^{2} + 972\,p^{4}\,q^{4}
 \\
 & & \mbox{} - 8208\,p^{2}\,q^{4} - 41904\,q^{4} - 318816\,p^{4}
\,q^{2} - 4320\,q^{6} - 211680\,p^{4} + 1025600 \\
 & & \mbox{} - 828576\,p^{3}\,q^{2} - 10368\,p\,q^{4} + 23328\,p
^{3}\,q^{4} - 11664\,p^{5}\,q^{2} - 1512\,q^{6}\,p \\
 & & \mbox{} - 124416\,p^{5})/(3\,p\,q^{2} - 4 + 36\,p + 6\,p^{2}
)^{3},  - {\displaystyle \frac {1}{9}} ( - 15984\,p^{6} + 27\,p^{
3}\,q^{6} - 1728\,p \\
 & & \mbox{} + 2916\,p^{9} - 1944\,p^{8} + 27\,q^{8} - 6912\,p^{3
} - 118512\,p\,q^{2} - 45360\,p^{2}\,q^{2} - 145152\,q^{2} \\
 & & \mbox{} - 16704\,p^{2} + 2916\,p^{6}\,q^{2} + 486\,p^{7}\,q
^{2} + 162\,p^{4}\,q^{4} + 3672\,p^{2}\,q^{4} + 44280\,q^{4} \\
 & & \mbox{} + 15552\,p^{4}\,q^{2} - 2052\,q^{6} - 26784\,p^{4}
 + 27000\,p^{3}\,q^{2} + 16848\,p\,q^{4} - 1944\,p^{3}\,q^{4} \\
 & & \mbox{} + 324\,p^{5}\,q^{2} - 378\,q^{6}\,p - 7776\,p^{5} - 
3584)/(3\,p\,q^{2} - 4 + 36\,p + 6\,p^{2})^{3}]
\end{eqnarray*}
\end{small}
\end{maplelatex}

\end{maplegroup}
\begin{maplegroup}
\begin{mapleinput}
\mapleinline{active}{1d}{symN(inv,1,1);}{%
}
\end{mapleinput}

\mapleresult
\begin{maplelatex}
\[
\mathit{the\ number\ of\ symmetries}=18
\]
\end{maplelatex}

\end{maplegroup}
\begin{maplegroup}
\begin{mapleinput}
\mapleinline{active}{1d}{f:=(p,q)->-q^2+1+3*p+p^3;}{%
}
\end{mapleinput}

\mapleresult
\begin{maplelatex}
\[
f := (p, \,q)\rightarrow  - q^{2} + 1 + 3\,p + p^{3}
\]
\end{maplelatex}

\end{maplegroup}
\begin{maplegroup}
\begin{mapleinput}
\mapleinline{active}{1d}{Pinv(f,3);}{%
}
\end{mapleinput}

\mapleresult
\begin{maplelatex}
\begin{small}
\begin{eqnarray*}
\lefteqn{[ - {\displaystyle \frac {1}{6}} (18 + 36\,p + 72\,p^{2}
 + 36\,p^{3} + 63\,p^{4} + 9\,p^{6} - 18\,q^{2} - 60\,p\,q^{2} + 
6\,p^{2}\,q^{2} - 12\,p^{3}\,q^{2}}\\
& & + 9\,q^{4} + p^{2}\,q^{4})/(p\,q^{2} - 3 + 3\,p + 3\,p^{2})^{2}
,  - {\displaystyle \frac {1}{6}} ( - 432\,p^{6} - 36\,p^{3}\,q^{
6} + 1944\,p + 648\,p^{7} \\
 & & \mbox{} + 243\,p^{8} - 5\,q^{8} + 1944 + 1944\,p^{3} - 1944
\,p\,q^{2} + 6156\,p^{2}\,q^{2} - 324\,q^{2} + 6804\,p^{2} \\
 & & \mbox{} - 540\,p^{6}\,q^{2} - 324\,p^{7}\,q^{2} + 54\,p^{4}
\,q^{4} - 684\,p^{2}\,q^{4} - 468\,q^{4} - 4428\,p^{4}\,q^{2} - 
40\,q^{6} \\
 & & \mbox{} + 3240\,p^{4} + 5184\,p^{3}\,q^{2} - 144\,p\,q^{4}
 + 216\,p^{3}\,q^{4} - 972\,p^{5}\,q^{2} - 84\,q^{6}\,p - 2592\,p
^{5})/ \\
 & & (p\,q^{2} - 3 + 3\,p + 3\,p^{2})^{3},  - {\displaystyle 
\frac {1}{9}} ( - 81 - 567\,p^{6} + p^{3}\,q^{6} - 81\,p + 27\,p
^{9} - 162\,p^{8} + q^{8} \\
 & & \mbox{} - 216\,p^{3} - 756\,p\,q^{2} - 945\,p^{2}\,q^{2} - 
243\,q^{2} - 405\,p^{2} + 27\,p^{6}\,q^{2} + 27\,p^{7}\,q^{2} + 9
\,p^{4}\,q^{4} \\
 & & \mbox{} + 306\,p^{2}\,q^{4} + 288\,q^{4} + 216\,p^{4}\,q^{2}
 - 19\,q^{6} - 729\,p^{4} - 918\,p^{3}\,q^{2} + 234\,p\,q^{4} \\
 & & \mbox{} - 18\,p^{3}\,q^{4} + 27\,p^{5}\,q^{2} - 21\,q^{6}\,p
 - 162\,p^{5})/(p\,q^{2} - 3 + 3\,p + 3\,p^{2})^{3}]
\mbox{\hspace{77pt}}
\end{eqnarray*}
\end{small}
\end{maplelatex}

\end{maplegroup}
\begin{maplegroup}
\begin{mapleinput}
\mapleinline{active}{1d}{symN(inv,1,0);
                                 }{%
}
\end{mapleinput}

\mapleresult
\begin{maplelatex}
\[
\mathit{the\ number\ of\ symmetries}=18
\]
\end{maplelatex}

\end{maplegroup}
\begin{maplegroup}
\begin{mapleinput}
\mapleinline{active}{1d}{f:=(p,q)->-q^2+1-p+p^3;}{%
}
\end{mapleinput}

\mapleresult
\begin{maplelatex}
\[
f := (p, \,q)\rightarrow  - q^{2} + 1 - p + p^{3}
\]
\end{maplelatex}

\end{maplegroup}
\begin{maplegroup}
\begin{mapleinput}
\mapleinline{active}{1d}{Pinv(f,3):}{%
}
\end{mapleinput}

\end{maplegroup}
\begin{maplegroup}
\begin{mapleinput}
\end{mapleinput}

\end{maplegroup}
\begin{maplegroup}
\begin{mapleinput}
\mapleinline{active}{1d}{symN(inv,1,0);}{%
}
\end{mapleinput}

\mapleresult
\begin{maplelatex}
\[
\mathit{the\ number\ of\ symmetries}=18
\]
\end{maplelatex}

\end{maplegroup}
\begin{maplegroup}
\begin{mapleinput}
\end{mapleinput}

\end{maplegroup}



\bibliographystyle{plain}

\end{document}